\documentstyle[fullpage,11pt,psfig,righttag,amssymb]{amsart}

\catcode`\@=11

\long\def\@savemarbox#1#2{\global\setbox#1\vtop{\hsize\marginparwidth 
  \@parboxrestore\tiny\raggedright #2}}
\newcommand\lref[1]{\ref{#1}%
\@ifundefined{r@DisplaY #1}{}{ (#1)}}

\newcommand\fakelabel[2]{\@bsphack\if@filesw {\let\thepage\relax
   \newcommand\protect{\noexpand\noexpand\noexpand}%
\xdef\@gtempa{\write\@auxout{\string
      \newlabel{#1}{{#2}{\thepage}}}}}\@gtempa
   \if@nobreak \ifvmode\nobreak\fi\fi\fi\@esphack}

\catcode`\@=12

\def\Empty{}
\newcommand\oplabel[1]{
  \def\OpArg{#1} \ifx \OpArg\Empty {} \else
  	\label{#1}
  \fi}

\long\def\realfig#1#2#3{
\begin{figure}[htbp]
\centerline{\psfig{figure=#2}}
\caption[#1]{#3}
\oplabel{#1}
\end{figure}}

\newtheorem{theoremSt}{Theorem}[section]

\newtheorem{exampleSt}[theoremSt]{Example}
\newtheorem{exerciseSt}[theoremSt]{Exercise}

\newcommand\MakeStEnv[1]{
  \newenvironment{#1}[1]{
  \begin{#1St} \oplabel{##1}%
  \global\def\CrntSt{\thetheoremSt}%
}{ 
  \end{#1St} }
  \newenvironment{#1+}[1]{
  \begin{#1St} \label{##1}%
  \label{DisplaY ##1}%
  \global\def\CrntSt{\thetheoremSt}%
  \def\Labl{##1}\ifx\Labl\Empty{} \else {\em (\Labl)\,}\fi%
}{ 
  \end{#1St} }
}
\MakeStEnv{theorem}
\MakeStEnv{corollary}
\MakeStEnv{proposition}
\MakeStEnv{lemma}
\MakeStEnv{definition}
\MakeStEnv{conjecture}

\newenvironment{example}[1]{
  \begin{exampleSt} \oplabel{#1}%
  \global\def\CrntSt{\thetheoremSt}%
  \normalshape%
}{ 
  \end{exampleSt} }

\newcommand{\startproof}[1]{%
\medbreak\mbox{}\noindent{\it Proof of #1:}%
}
\newcommand{\finishproof}[1]{ 
  \def\FPArg{#1}
  \ifx\FPArg\Empty
  	\newcommand\FPArg{\CrntSt}  \fi
  \smallbreak\noindent\makebox[\textwidth]{\hfill\fbox{\FPArg}}
  \medbreak\noindent
}

\newcommand{\bfheading}[1]{\par\smallskip\noindent {\bf #1}}

\def\IMSmarkvadjust{0 pt}
\def\IMSmarkhadjust{0 pt}
\def\IMSmarkhpadding{0 pt}
\def\IMSpubltext{Published in modified form:}
\def\SBIMSMark#1#2#3{
 \font\SBF=cmss10 at 10 true pt
 \font\SBI=cmssi10 at 10 true pt
 \setbox0=\hbox{\SBF \hbox to \IMSmarkhpadding{\relax}
                Stony Brook IMS Preprint \##1}
 \setbox2=\hbox to \wd0{\hfil \SBI #2}
 \setbox4=\hbox to \wd0{\hfil \SBI #3}
 \setbox6=\hbox to \wd0{\hss
             \vbox{\hsize=\wd0 \parskip=0pt \baselineskip=10 true pt
                   \copy0 \break%
                   \copy2 \break%
                   \copy4 \break}}
 \dimen0=\ht6   \advance\dimen0 by \vsize \advance\dimen0 by 8 true pt
                \advance\dimen0 by -\pagetotal
	        \advance\dimen0 by \IMSmarkvadjust
 \dimen2=\hsize \advance\dimen2 by .25 true in
	        \advance\dimen2 by \IMSmarkhadjust

  \openin2=publishd.tex
  \ifeof2\setbox0=\hbox to 0pt{}
  \else 
     \setbox0=\hbox to 3.1 true in{
                \vbox to \ht6{\hsize=3 true in \parskip=0pt  \noindent  
                {\SBI \IMSpubltext}\hfil\break
                {\it J. Differential Geom.}~{\bf 47} (1997) 17--94 
                \vfill}}
  \fi
  \closein2
  \ht0=0pt \dp0=0pt
 \ht6=0pt \dp6=0pt
 \setbox8=\vbox to \dimen0{\vfill \hbox to \dimen2{\copy0 \hss \copy6}}
 \ht8=0pt \dp8=0pt \wd8=0pt
 \copy8
 \message{*** Stony Brook IMS Preprint #1, #2. #3 ***}
}

\def\IMSmarkvadjust{-25pt}
\newcommand{\BrkOK}{\discretionary{}{}{}}

\newcommand\AAA{{\cal A}}
\newcommand\BB{{\cal B}}
\newcommand\CC{{\cal C}}
\newcommand\DD{{\cal D}}

\newcommand\FF{{\cal F}}
\newcommand\GG{{\cal G}}
\newcommand\HH{{\cal H}}

\newcommand\JJ{{\cal J}}
\newcommand\KK{{\cal K}}
\newcommand\LL{{\cal L}}
\newcommand\MM{{\cal M}}
\newcommand\NN{{\cal N}}

\newcommand\PP{{\cal P}}
\newcommand\QQ{{\cal Q}}
\newcommand\RR{{\cal R}}
\newcommand\SS{{\cal S}}
\newcommand\TT{{\cal T}}
\newcommand\UU{{\cal U}}
\newcommand\VV{{\cal V}}

\newcommand\XX{{\cal X}}

\newcommand\PMF{{\PP\kern-2pt\MM\FF}}

\newcommand\PML{{\PP\kern-2pt\MM\LL}}

\newcommand\ep{\epsilon}

\newcommand\union{\cup}
\newcommand\intersect{\cap}
\newcommand\bbR{{\mathord{\text{I\kern-2pt R}}}}        
\newcommand\bbH{{\mathord{\text{I\kern-2pt H}}}}        

\newcommand\C{{\Bbb C}}
\newcommand\D{{\Bbb D}}
\newcommand\Z{{\Bbb Z}}
\newcommand\R{{\Bbb R}}

\newcommand\Hyp{{\Bbb H}}

\newcommand\bigrightarrow[1]{\hbox to #1{\rightarrowfill}}
\newcommand\bigleftarrow[1]{\hbox to #1{\leftarrowfill}}

\newcommand\homeo{\cong}
\newcommand\boundary{\partial}
\newcommand\semidir{\mathrel{\hbox{\vrule depth-.03ex height1.1ex\kern-0.15em$\times$}}}

\newcommand\til{\widetilde}

\newcommand{\ssm}{\setminus}
\newcommand{\diam}{\operatorname{diam}}

\numberwithin{equation}{section}

\def\marginpar#1{}

\begin{document}

\def\Draft{$\text{(Revision Draft)}$}
\title[Laminations in holomorphic dynamics ]{Laminations in holomorphic
dynamics}
\author{Mikhail Lyubich}
\author{Yair Minsky}
\address{Mathematics Department and IMS, SUNY Stony Brook}

\maketitle
\thispagestyle{empty}
\SBIMSMark{1994/20}{December 1994}{revised version: June 26, 1998}
\newcommand{\id}{\operatorname{id}}
\newcommand{\dist}{\operatorname{dist}}
\newcommand{\orb}{\operatorname{orb}}
\newcommand{\cl}{\operatorname{cl}}
\newcommand{\inter}{\operatorname{int}}
\newcommand{\Dis}{\operatorname{Dis}}
\newcommand{\Aff}{\operatorname{Aff}}
\def\ra{\rightarrow}
\def\eps{\epsilon}
\def\bz{{\bold{z}}}
\def\bzeta{\text{\boldmath{$\zeta$}}}
\newcommand{\comm}[1]{}
\newcommand{\sm}{\setminus}
\newcommand{\dens}{\operatorname{dens}}
\newcommand{\meas}{\operatorname{meas}}
\newcommand{\T}{{\Bbb T}}

\newcommand{\an}{\AAA^{\bold{n}}} 
\newcommand{\al}{\AAA^{\bold{\ell}}} 
\newcommand{\ac}{\AAA}	 
\newcommand{\hn}{\HH^{\bold{n}}}  
\newcommand{\hl}{\HH^{\bold{\ell}}}	 
\newcommand{\hc}{\HH}	 
\newcommand{\jr}{\JJ^{\bold{r}}} 
\newcommand{\jn}{\JJ^{\bold{n}}}	 
\newcommand{\jl}{\JJ^{\bold{\ell}}}	 
\newcommand{\jc}{\JJ}	 
\newcommand{\fr}{\FF^{\bold{r}}} 
\newcommand{\fn}{\FF^{\bold{n}}}	 
\newcommand{\fl}{\FF^{\bold{\ell}}}	 
\newcommand{\fc}{\FF}	 
\newcommand{\uu}{\UU}	 
\newcommand{\ua}{\UU^{\bold{a}}} 
\newcommand{\uh}{\UU^{\bold{h}}}	 
\newcommand{\k}{\KK}	 
\newcommand{\ka}{\KK^{\bold{a}}}	 
\newcommand{\kh}{\KK^{\bold{h}}}	 

\newcommand{\QED}{\rlap{$\sqcup$}$\sqcap$\smallskip}

\setcounter{tocdepth}{1}
\tableofcontents

\section{A missing line in the dictionary.}
\label{intro}

There is an intriguing dictionary between two branches of conformal dynamics:
the theory of Kleinian groups and dynamics of rational  maps. 
This dictionary was introduced by Sullivan, and led him
in the early 80's to the no wandering domains theorem, deformation theory and
geometric measure theory for holomorphic maps.  Thurston's rigidity
and realization theory, developed at the same time, was also motivated
by this analogy.
More recently, McMullen has made important 
contributions to the renormalization theory motivated by the analogy 
with 3-manifolds which fiber over the circle 
\cite{mcmullen:icm90,mcmullen:renormfiber}.   

However, the translation from one language to
another, as in usual life, 
is not automatic. There are concepts and methods in each of
these fields which only barely allow translation to the other one. And
even when it is possible, the results achieved are often complementary
(see Sullivan's table in \cite{sullivan:QCDII} of the results on the
structural stability and hyperbolicity problems).

In this paper we explore a construction which attempts to provide an
element of the dictionary that has so far been missing: an explicit
object that plays for a rational map the role played by the hyperbolic
3-orbifold quotient of a Kleinian group. To build this object we
replace the notion of manifold by ``lamination'', which is a
topological object whose local structure is the product of Euclidean
space by a (possibly complicated) transverse space. 

Another  goal of this work is to study
 the space of backward orbits of a rational function.
 Since Fatou and Julia,
inverse branches of iterated rational functions have 
 played a crucial role in the theory. 
Unfortunately, the space of such branches, with its natural topology, is 
wild  (should be compared with the H\'enon
attractor), and may deserve to be called a ``turbulation''. 
By imposing a finer topology and completing,
we turn this space  into an  affine  lamination, in the hope that
this will tame it.

Laminations were introduced into conformal dynamics by Sullivan, whose
{\it Riemann surface laminations} play a role similar
to that of Riemann surfaces for Kleinian groups (see
\cite{sullivan:bounds,sullivan:universalities}
 or \S \ref{Julia and Fatou sets} and Appendix 1 of this paper).
These are objects which locally look like a product of a complex disk
times a Cantor set. Sullivan associated such a holomorphic object to 
any $C^2$-smooth expanding circle map. The construction 
involves  ``conformal extension"  of a non-analytic one-dimensional map
(see Appendix 1).

In this paper we go one dimension up and make a ``hyperbolic
three dimensional extension" of a non-M\"{o}bius map.
This object is called a {\em hyperbolic orbifold 3-lamination}
and can be constructed in the following way. 

\smallskip\noindent {\it Step 1: the natural extension.} 
Consider the full natural  extension $\hat f: \NN_f\rightarrow\NN_f$
of a rational map $f$ (points of $\NN_f$ are backward orbits 
$\hat z=(\cdots \mapsto z_{-1}\mapsto z_0)$ of $f$). 

\smallskip\noindent {\it Step 2: the regular leaf space}.
Restrict $\hat f$ to
the ``regular part" $\RR_f$ of $\NN_f$ where the inverse iterates
branch only finitely many times. This space is a union of
leaves which are  non-compact 
Riemann surfaces, simply connected except for Herman rings,
that is, hyperbolic or parabolic planes. It can be viewed as
a Riemann surface with uncountably many sheets where all inverse
iterates $f^{-n}$ live simultaneously.

\smallskip\noindent {\it Step 3: affine orbifold lamination}.
Consider the subset $\an_f$ of $\RR_f$ consisting of parabolic leaves.
The parabolic leaves possess a canonical
affine structure preserved by the map. However this structure is not
necessarily continuous in the transverse direction. 
To make it continuous, we refine the topology on $\an_f$, obtaining	
a space $\al_f$ with a laminar 
structure\footnote{ A different approach to this part was independently
  suggested by Meiyu Su who imposed a laminar topology associated to the
  transversal measure structure \cite{Su}.}.
 We then complete $\al_f$ to obtain a final object
$\ac_f$ some of whose new leaves may be 2-orbifolds.

This step is technically the hardest.

\smallskip\noindent {\it Step 4: three-dimensional extension}. 
Each affine leaf is naturally the boundary of a three-dimensional
hyperbolic space (in the half-space model). The union of these spaces
forms a hyperbolic
orbifold 3-lamination  $\HH_f$ with $\hat f$ acting 
properly discontinuously, and by isometries on the leaves.

\smallskip\noindent {\it Step 5: quotient}. Finally taking
the quotient $\HH_f/\hat f$  of this lamination by $\hat f$ we
obtain the desired hyperbolic orbifold 3-lamination. 

\realfig{z plus epsilon picture}{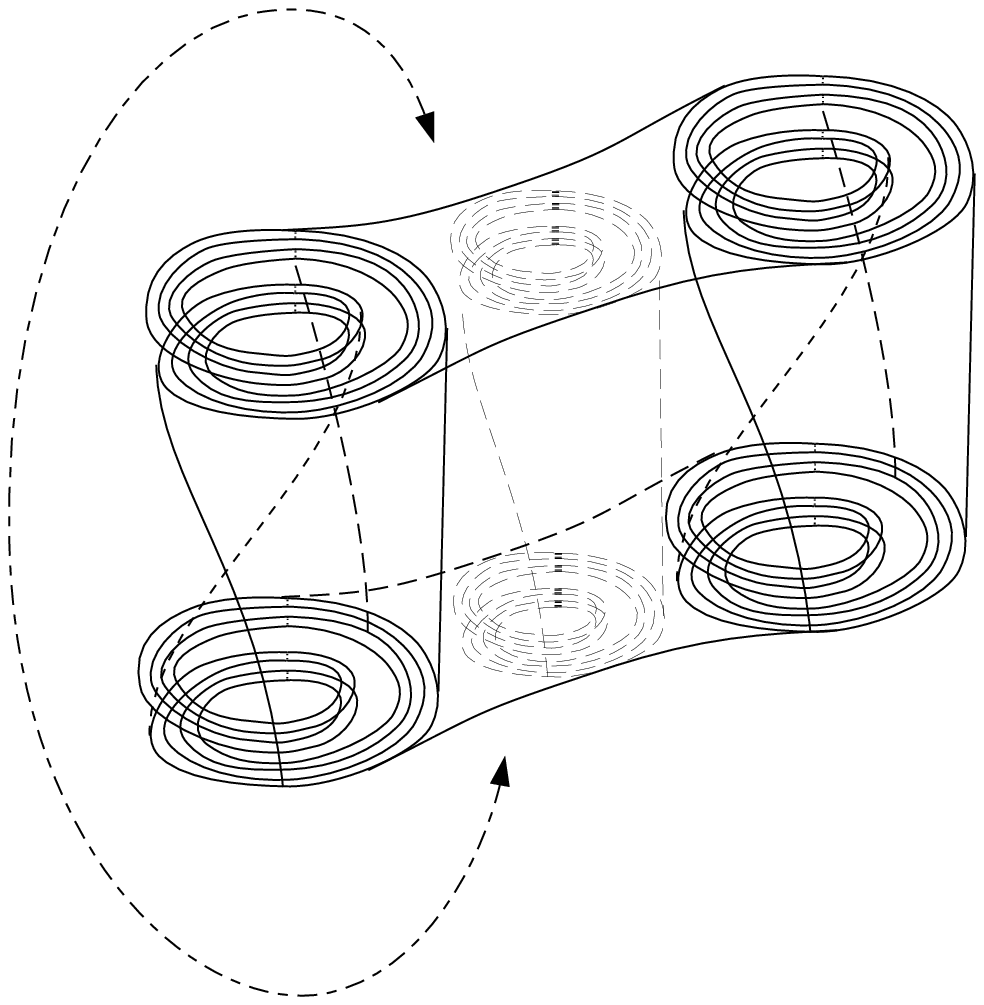}{A caricature of the
hyperbolic 3-lamination for $z^2+\ep$. 
Sullivan's solenoidal Riemann surface lamination is shown as the
mapping torus of the doubling map on the solenoid (final gluing
indicated by the double arrow), and the product of this with the
interval is suggested. See discussion in \S\ref{3d lams}.}

\medskip

We also define the convex core of the lamination
$\HH_f/\hat f$ and prove that it is
compact if and only if $f$ is critically non-recurrent without parabolic
points. 
Using this criterion, we prove a rigidity theorem \ref{pcf rigidity}  for
critically non-recurrent maps which extends 
Thurston's rigidity theorem for post-critically
finite rational maps (see Douady-Hubbard \cite{douady-hubbard:pcf}).
Our three-dimensional proof gives an explicit connection
between Thurston's and Mostow's rigidity theorems.

\medskip

The structure of the paper is as follows: 

{\bf\S\ref{laminations background}.} Basic notions of
laminations and orbifold laminations.  

{\bf\S\ref{natural extension and regular part}.} 
The natural extension $\NN_f$, and its
regular part $\RR_f$. The space $\RR_f$ consists of backward orbits which
have neighborhoods whose pullbacks 
hit the critical points only finitely many times. 
This space can be decomposed into leaves that admit
a natural conformal structure.  
We show that (with the exception of Herman rings) the leaves of $\RR_f$ are
simply connected non-compact Riemann surfaces, i.e. either hyperbolic
or parabolic planes.

We discuss criteria for when $\RR_f$ is all of $\NN_f$ except for a
finite set, and when $\RR_f$ is open in $\NN_f$.
This discussion crucially depends on a theorem by R.~Ma\~n\'e
on the  behaviour of non-recurrent critical points \cite{mane:fatou}. 

{\bf\S\ref{affine structure}.} Here we discuss the affine part $\an_f$
 of $\RR_f$, which leads us to  the type problem for the leaves. 
This problem seems to be intimately related to the geometry of the Julia set.
Parabolicity of leaves reflects ``some'' (but not
necessarily uniform) expansion -- see Lemma \ref{derivative condition}.
We give several simple criteria for parabolicity and apply them to some
special cases. In particular, all leaves of the real Feigenbaum quadratic
are parabolic.  This follows from an expansion property of $f$ with
respect to a hyperbolic metric (compare McMullen 
 \cite{mcmullen:renormalization}). 
The only  examples known to us of hyperbolic leaves are the invariant lifts of
Siegel disks and Herman rings.

We also give an explicit formula for the affine coordinate on a
parabolic leaf. It generalizes the classical formulas for the 
linearizing K\"{o}nigs and Leau-Fatou coordinates near repelling and
parabolic points. From this point of view the affine structures on the
leaves are  just  the linearizing coordinates along the backward 
orbits of $f$.  

{\bf\S\ref{pcf orbifolds}.} Here we carry out Step 3 of the construction 
for the
post-critically finite case, the construction of an affine orbifold
lamination $\ac_f$. We refine the topology on $\an_f$ to separate
leaves which branch inconsistently over the sphere, and 
enlarge $\an_f$ to $\ac_f$ by making several copies
of the post-critical periodic
leaves, and replacing the original affine structure on some of them by 
an orbifold affine
structure. This is the price we pay for having the affine structure
transversally continuous, while keeping the lamination complete (in an
appropriate sense).

{\bf\S\ref{3d lams}.} We define the notion of an
{\em orbifold affine  extension}  $\hat  f: \ac\rightarrow \ac$ of a
rational map $f$, and show that it is naturally the boundary at infinity
for an orbifold  hyperbolic 3D extension 
$\hat f: \hc\rightarrow \hc$. We prove that the action of $\hat f$
on $\hc$ is properly discontinuous, so that the quotient 
$\hc/\hat f$ inherits the structure of a hyperbolic orbifold 3-lamination.

Then we introduce and discuss the notion of the convex core $\CC_f$
in $\hc_f/\hat f$,
which will play a key role in the rigidity argument.

We  also describe the topological structure of the 3-lamination 
associated to quadratics
 $p_{\epsilon}: z\mapsto z^2+\epsilon$ with $\epsilon$
inside of the main cardioid of the Mandelbrot set. We show that it is
homeomorphic to $\SS\times (0,1)$ where $\SS$ is the Sullivan lamination.
So, like in the case of quasi-Fuchsian groups, the 3-lamination
connects the 2-laminations associated
to the attracting basins of $p_{\epsilon}$.

At the end of this section we discuss the ``scenery flow"
introduced by A. Fisher as an analogue of the geodesic flow on
3-manifolds.  The phase space of this flow,
constructed in \cite{bedford-fisher-urbanski}  for rational maps
satisfying axiom A, is loosely speaking the
set of ``pictures", that is all possible rescalings of the
infinitesimal germs of the Julia set.  This scenery flow
is topologically equivalent to the ``vertical geodesic flow" on the
3-lamination over the lifted Julia set.

This vertical geodesic flow is an extra piece of structure which makes
a difference between 3-laminations of rational maps and 3-manifolds
of Kleinian groups. An equivalent way of viewing this structure is by
saying that there  exists  a preferred $\hat f$-invariant
cross-section, ``$\infty$'', at the boundary of the lamination  $\hc_f$.

{\bf\S\ref{universal construction}.}
In this section we give a general construction of the affine and
hyperbolic orbifold laminations associated to a rational map.
The main hurdle is, as in the post-critically finite case, the
fact that a sequence of disks in $\RR_f$  whose projection to the sphere
is branched can limit onto a disk on which the projection is univalent. In the
general case sorting out the different branching types is more
involved since the set of points where this happens is no longer
finite. Thus many copies of a leaf, possibly a continuum, must be
added. One can keep track of this, and define an appropriate topology,
using the affine structures themselves and their limiting behavior.

The self-organizing idea for this construction is to observe that the
natural projection $\pi:\NN_f \to \bar\C$ gives a meromorphic function
on each leaf of $\an_f$, and this family of functions has a natural
topology which induces a topology for $\an_f$. In fact, the space of
non-constant meromorphic functions on $\C$  with the right action
of the affine group serves as a ``universal'' lamination on 
which every rational function acts. For any fixed $f$ the
structure of $\ac_f$ and $\hc_f$ can be extracted from the attractor
of $f$ in this universal space.

In conclusion we prove (using Ahlfors' five islands theorem)
 that every lamination $\hc_f$
is minimal, except  for the Chebyshev and Latt\`es maps.
In these special cases, the lamination becomes minimal
after removing the invariant isolated leaf. This is the
characteristic property of these remarkable maps from the
lamination point of view. 

{\bf\S\ref{convex cocompactness section}.} 
In this section we prove that $f$ is convex co-compact
(that is, its convex core $\CC_f$ is compact) if and only if
it is critically non-recurrent
and does not have parabolic periodic points.
 Note that thus convex
co-compactness differs from hyperbolicity, while these two notions 
are equivalent for Kleinian groups (one more illustration of the
loose nature of the dictionary). 

We also define the conical limit set
 and give in these terms a criterion of  convex cocompactness. 
We then study ergodic properties of the conical limit 
set by means of the blow-up technique (on the lamination level) and
Ahlfors' harmonic extension method. Along the lines we obtain
the lamination insight on the existence of invariant line fields
for the Latt\`es examples: it comes from the existence of the 
isolated leaves. 

{\bf\S\ref{pcf rigidity section}.} 
This section contains the three-dimensional proof of rigidity for
convex co-compact maps. 
 
We start by lifting the topological equivalence between the maps to a
quasi-isometry  $\hat h$ 
between their 3-laminations (using the convex co-compactness).
It follows that $\hat h$ is quasi-conformal on the leaves of the
affine extension.
 This reduces the problem to the existence of invariant
line fields on the Julia set of $\ac_f$, which was analyzed in 
the previous section.

{\bf\S\ref{conjectures}.} 
Conjectures and further program.

{\bf\S\ref{appendix1}.} In the first appendix we outline Sullivan's
costruction 
  of the Riemann surface lamination associated  to an expanding map of
   the circle.  We also give a globalization construction for  the 
 natural extension of polynomial-like maps via the inductive limit procedure.  

{\bf\S\ref{appendix2}.} Appendix 2 fills in some necessary
background, all 
of which is well-known to those who work in either dynamics or geometry,
but not always to both. It also fixes some terminology and notation.

\bigskip{\bf Acknowledgements.} Since the first drafts of this manuscript 
appeared in spring 1994, 
we have received valuable feedback from many people, and
we take this opportunity to thank all of them.
John Milnor read parts of
the manuscript and made many useful suggestions and comments. 
He also gave a simpler proof of Lemma \ref{not elliptic} 
on non-compactness of the
leaves. Alberto Baider simplified the proof of Lemma \ref{leaves are planes}
on simple connectivity of the leaves. We thank Dennis Sullivan for
several enjoyable discussions and stimulating questions, and Curt
McMullen for some useful comments on renormalization.
We  also enjoyed talks  with Alby Fisher about scenery flows, 
as well as discussions with Meiyu Su who is independently
working on the structure of the natural extension of a rational map from
a measure-theoretical point of view. After the preprint was issued
(IMS at Stony Brook, 1994/20) we have
found out that Jeremy Kahn
has independently considered the idea of a 3D lamination and its possible
application to the measure problem.

\section{Laminations: general concepts}
\label{laminations background}

In this paper, a {\it lamination} will be a Hausdorff topological space
$\XX$ equipped with 
a covering $\{U_i\}$ and coordinate charts  $\phi_i:U_i\to T_i\times D_i$,
where $D_i$ is homeomorphic to a domain in $\R^n$ and $T_i$ is  
a topological space. The transition maps
$\phi_{ij}=\phi_i\circ\phi_j^{-1}:\phi_j(U_i\intersect U_j) \to
\phi_i(U_i\intersect U_j)$ are required to be homeomorphisms that
take leaves to leaves (see Sulivan \cite{sullivan:universalities} and
Candel \cite{candel}).

Subsets of the form $\phi_i^{-1}(\{t\}\times D)$ are called {\it local
leaves.} The requirement on the transition maps implies that the local
leaves piece together to form {\it global leaves}, which are
$n$-manifolds immersed injectively in $\XX$. 

As usual we may restrict the class of transition maps to obtain finer
structures on $\XX$. If $D_i$ are taken to lie in $\C$ and $\phi_{ij}$
are conformal maps, we call $\XX$ a {\it Riemann surface lamination}
and note that the global leaves have the structure of Riemann
surfaces. If $\phi_{ij}$ are further restricted to be complex affine
maps $z\mapsto az+b$, then we call $\XX$ a {\it (complex) affine lamination},
and the global leaves have a (complex) affine structure. 
If the leaves of an affine lamination are isomorphic to the complex
plane, we also call it a $\C$-lamination. One can similarly consider
real affine laminations, but as they will not play a role in this paper
we shall assume from now on that ``affine" means ``complex affine".

If $D$ are taken to lie in $\Hyp^n$ and $\phi_{ij}$ are hyperbolic
isometries, then $\XX$ is an {\it $n$-dimensional hyperbolic
lamination}, or {\it hyperbolic $n$-lamination.} 
 In the case when all leaves
of the lamination
are hyperbolic spaces, let us call it an $\Hyp^n$-lamination.

When the laminated space $\XX$ is a (smooth/analytic)
 manifold, the lamination is
usually called a {\it foliation}. It is called {\it smooth/analytic} if
there is a smooth/analytic  atlas of laminar local charts.

We shall need the notion of distance between affine structures on a
Riemann surface.  
Let $S$ be a Riemann surface supplied with two affine structures
$\AAA_1$ and $\AAA_2$. 
Let $\phi_1$ and $\phi_2$ be any two local charts of the
structures $\AAA_1$ and $\AAA_2$ respectively, $\psi=\phi_1\circ \phi_2^{-1}:
U_1\rightarrow U_2 $
be the transition function. 
Then we may define
$$\dist(\AAA_1, \AAA_2)=\sup_{\phi_1, \phi_2} 
\Dis({\phi_1\circ \phi_2^{-1}}), $$
where $\Dis$ stands for the distortion (see Appendix 2, \S\ref{appendix2}).

We will encounter situations where a Riemann surface lamination $\RR$ can be
refined to give an affine lamination. Suppose that the global leaves
of $\RR$ admit affine structure -- that is, each global leaf $L$ admits a
collection of conformal coordinate charts with affine transition maps. 
We say
that these affine structures {\it vary continuously} in $\RR$ if, for
any product box  $U=T\times D$ the induced family of affine structures
on $D$ vary continuously with $T$, in the sense of the above notion of
distance. 

In other words, continuity of affine structure means that
for each coordinate chart $\phi:U\to T\times D$ there
is a choice of coordinate $\psi_t:\phi^{-1}(\{t\}\times D) \to \C$
for each  $t\in T$, so that $\psi_t$ is a restriction of an affine
coordinate chart on a global leaf, and so that the family
$\psi_t\circ\phi^{-1}(t,\cdot) :D\to\C$ varies continuously with $t$. 
The following is easy to check: 

\begin{lemma}{leafwise to global affine}
A continuous family of affine structures on the
global leaves of a Riemann surface lamination $\RR$ induces an affine
lamination structure on $\RR$  compatible with the original
structure. 
\end{lemma}

Similarly the Riemann surface lamination can be viewed as a 
topological lamination with transversally continuous conformal 
structure on the leaves. 

\subsection{Orbifold laminations}
\label{orbifold defs}
In analogy with Thurston's notion of orbifolds
(see Thurston \cite{wpt:textbook},
Scott \cite{scott:survey} and also Satake's similar notion of
V-manifolds, \cite{satake:vmanifolds}), we may define
an orbifold lamination to be a space for which every point has a
neighborhood that is either homeomorphic to a standard product box
neighborhood in a lamination, or to a quotient of such a box
by a finite leaf-preserving group (called an {\em orbifold box}).

If the covering box has an affine (or conformal, or hyperbolic) structure
which is preserved by the finite group, then we say that the orbifold
box inherits an
orbifold affine (or conformal, or hyperbolic) structure.

For example, let $T\times D$ be a product box, with $D$ a two
dimensional disk, let $\sigma:T\to T$ be a finite-order map and let
$\rho:D\to D$ be a finite order rotation of $D$.  Then the map
$\sigma\times\rho$ generates a finite cyclic group action on $T\times
D$ and the quotient is an orbifold box. Cycles of $\sigma$ of order
not divisible by the order of $\rho$ (fixed points, for example) give
rise to quotient leaves with orbifold points.

See also \cite{wpt:ratmaps,douady-hubbard:pcf} 
for the use of (regular 2-dimensional) orbifolds in the context of
post-critically finite maps. 

\begin{example}{deg 2 example}
This example illustrates how orbifold boxes will arise in 
\S\ref{pcf orbifolds}.
Let $K$ be a Cantor set and $K'= K\ssm\{a\}$ for some $a\in K$, and let
$\pi:\hat D \to D$ be a doubly branched map. 
Let $B$ denote $(K'\times \hat D)\union (\{a\}\times D)$, 
topologized so that a sequence $(b_i,z_i)$ in $K'\times \hat D$ converges
to $(a,z)$ in $\{a\}\times D$ if and only if $b_i\to a$ and
$\pi(z_i) \to z$. 

We can then express $B$ as an orbifold box, by letting $T$ be the
double of $K$, with both copies of $a$ identified, and  $\sigma :T\to T$ the
map that interchanges copies. Let $\rho:\hat D \to \hat D$ be the
involution that interchanges pairs of preimages of points in $D$. Then
$(T\times \hat D)/(\sigma\times \rho)$ is exactly $B$.
\end{example}

\section{Natural extension and its regular part.}
\label{natural extension and regular part}

\subsection{Natural extension}
\label{natural extension}
Let 
$f: \bar{\C}\rightarrow \bar{\C}$ be a rational endomorphism of the Riemann 
sphere. Let us consider the space of its backward orbits:
$${\cal N}={\cal N}_f=\{\hat z=(z_0, z_{-1},\ldots): z_0\in \bar{\C},
 f: z_{-(n+1)}\mapsto z_{-n}\},$$
with topology induced by the product topology in 
$\bar\C\times\bar\C\times\ldots$.
This is a compact space projected down to $\bar{\C}$ by 
$\pi: \hat z\mapsto z_0$.
 The endomorphism  $f$  naturally lifts to a
homeomorphism $\hat f: {\cal N}\rightarrow {\cal N}$
 as $\hat f(\hat z)=(fz_0, z_0, z_{-1},\ldots)$.
(The inverse map forgets the first coordinate of the backward orbit).
Moreover, $\pi\circ \hat f=f\circ \pi$. In dynamics the map $\hat f$
is usually called  {\sl the natural extension} of $f$.
In algebra this object is also called
{\sl the projective (or inverse) limit} of
$$\bar{\C}\underset{f}\leftarrow\bar{\C}\underset{f}\leftarrow\bar{\C}
\underset{f}\leftarrow\ldots$$ 

One can also think of a point $\hat z\in \NN$ as a {\it full orbit}
$\{z_n\}_{n=-\infty}^{\infty}$, where $f: z_n\mapsto z_{n+1}$. 
(But don't confuse them with {\it grand
orbits} generated by the equivalence relation $z\sim \zeta$ if there exist
natural $m$ and $n$ such that $f^m z=f^n\zeta$). Along with the projection
$\pi\equiv \pi_0$ let us also consider  projections 
$\pi_n: \NN_f\rightarrow \bar\C$ such that $\pi_n(\hat z) = z_n$.
Clearly $\pi_n=f^{n-m}\circ\pi_m$ for $n\geq m$.

Given a (forward) invariant set  $X\subset \bar\C$,  let $\hat X\subset \NN_f$ 
denote its {\em invariant lift} to $\NN_f$, that is, the set of orbits
$\{z_n\}\subset X$. 
This is nothing but the natural extension of $f|X$. 
Note that it differs from $\pi^{-1} X$, unless $X$ is completely
invariant (that is, $f^{-1}X=X$).

Let $\hat z=(z_0,z_{-1},\ldots)\in \NN_f$,
$D$ be a topological disk containing $z_0$, and $N$ be a natural number. 
Consider  the pullback $D_0,D_{-1},\ldots$ of $D$ along $\hat z$. That
is, $D_{-n}$ is the component of $f^{-n}(D)$ containing $z_{-n}$.
Let us define the following ``boxes":
\begin{align}\label{box}
B(D, \hat z,  N) &= \pi_{-N}^{-1}(D_{-N}) \\\nonumber
&=\{\hat\zeta=(\zeta_0, \zeta_{-1},\ldots)\in \NN_f:\;
  \zeta_{-N}\in D_{-N}\},
\end{align}
which form a basis of the topology in $\NN_f$.
For $N=0$ we will shorten the notation as
$B(D,\hat z)\equiv B(D, \hat z, 0)$.

\subsection{The regular leaf space.}
\label{regular part}
Let us say that a point $\hat z=(z_0, z_{-1},\ldots)\in {\cal N}$ is 
{\em regular} if there is
neighborhood $U$ of $z_0$ in $\bar{\C}$ whose pullback $U_{-n}$
along the backward orbit $(z_0, z_{-1},\ldots)$ is eventually univalent. 
Let ${\cal R}={\cal R}_f$ denote the set of regular
points of the natural extension.
This set is clearly completely invariant. 
Moreover, if $z_0$ is outside the $\omega$-limit set $\omega(C)$
of the critical points, then $\hat z\in \RR_f$
(see Appendix 2). 

The path connected components of  $\RR$ will be called the ${\sl leaves}$
and denoted by $ L(\hat z)$ for $\hat z\in\RR$. 

\begin{lemma}{leaf structure}
The leaves $L(\hat z)$ possess an intrinsic 
topology and analytic structure such that the projection 
$\pi: L(\hat z)\rightarrow \bar{\C}$ is analytic.
The branched points are the backward orbits passing through
critical points. Moreover $\hat f: L(\hat z)\rightarrow L(\hat
f\hat z)$
is a biholomorphic isomorphism.
\end{lemma}

\begin{pf}
Let $\hat z=(z_0, z_{-1},\ldots)\in {\cal R}$. 
Then there is a neighborhood $U\ni z_0$ whose pull-back $U_{-n}$ along
the orbit $z_{-n}$ is eventually univalent. Let us take 
$\hat U=\{\hat\zeta=(\zeta_0,\zeta_{-1},\ldots): \zeta_{-n}\in U_{-n}\}$
as a base neighborhood of $\hat z$ (also called a leafwise neighborhood).

Let $f: U_{-(n+1)}\rightarrow U_{-n}$ be univalent for $n\geq N$. 
Then the map $\pi_{-N}: \hat \zeta\mapsto \zeta_{-N}$ is a homeomorphism between
$\hat U$ and $U_{-N}$. Let it be our local chart. The transition 
functions are just appropriate iterates of $f$, so that this provides
us with a complex structure.

The last two statements are obvious.  
\end{pf}

We may characterize the leaves in dynamical terms via the following
observation. 
\begin{lemma}{backward asymptotic on leaf}
Two points ${\hat z}$ and $\hat \zeta$ in $\RR_f$ belong to the same
leaf iff the following holds. 
There is a sequence of paths $(\gamma_{-n})$ in $\bar\C$ such
that $\gamma_{-n}$ connects $z_{-n}$ to $\zeta_{-n}$, and
$f(\gamma_{-n}) = \gamma_{-n+1}$. Furthermore, for $n$ sufficiently
large there are neighborhoods $U_{-n}$ of $\gamma_{-n}$ such that
there is a branch $g$ of $f^{-1}$ defined on $U_{-n}$  and $f(U_{-n}) =
U_{-n-1}$. In particular 
$\zeta_{-n}$ can be obtained from $z_{-n}$ by analytic continuation of
$f^{-1}$ along $\gamma_{-n+1}$.
\end{lemma}

\begin{pf} Assume that $\hat z$ and $\hat \zeta$ are on the same leaf
and let $\hat\gamma$ be a path connecting them.
We may represent any such path as a sequence
of paths $(\gamma_{-n})$ in $\bar\C$ such
that $\gamma_{-n}$ connects $z_{-n}$ to $\zeta_{-n}$, and
$f(\gamma_{-n}) = \gamma_{-n+1}$. 
Since each point in $\RR_f$ has a neighborhood whose projections are 
eventually univalent, we take a finite covering of $\hat \gamma$ 
and consider its projections by $\pi_{-n}$ for $n$ sufficiently large.
These are the neighborhoods $U_{-n}$. 

Conversely, given the sequence $\gamma_{-n}$ satisfying the
conditions, it is immediate that the
path $\hat\gamma$ in $\NN_f$ that they define in fact lies in $\RR_f$.
\end{pf}
 
By {\it local leaves} in a box $B(D, \hat z,  N)$ 
we will mean the components of intersection of the
global leaves with this box.

Unfortunately these boxes in general don't have a product structure,
so that $\RR_f$ is not always a Riemann surface lamination. 
For this reason $\RR_f$ will be called a {\it  (conformal) leaf space},
 that is,
a space which is decomposed into the union of leaves supplied with
smooth (conformal) structure. 
Actually the leaves behave so wildly
(keep in mind the Henon map) that one might rather call the
space a ``{\it turbulation}". 

However, if the orbits of the critical points don't meet $D_{-N}$ then
$B(D, \hat z,  N)\approx T\times D_{-N}$, 
where $T$ may be identified with the fiber $\pi^{-1}(z_{-N})$.
The local leaves in this box correspond to slices
$\{t\}\times D_{-N}$.

\subsection{Topology of the leaves}
\label{topology of the leaves}
Our main tool in this section will be the Shrinking Lemma (see
Appendix 2), which states roughly that, in a uniform sense, backward
iterates of a region on which the branching of $f$ is bounded have
(spherical) diameters that shrink to 0. This holds except if the
iterates remain 
in a {\em rotation domain} -- a Siegel disk or Herman ring -- for all
time. 

Let us first consider some exceptional cases:
If a component $W$ of the Fatou domain $F_f$ is a rotation domain,
then its invariant lift $\hat
W$, consisting of all orbits which remain in $W$ for all time, is a 
full leaf of $\RR_f$, and $\pi:\hat W\to W$ is a conformal
equivalence. The second part is obvious since
$f|_W:W\to W$ is a 1-1 conformal map.
It only remains to check that $\hat W$ is not properly
contained in a leaf. That is, we must check that  any point on
$\boundary \hat W$ in $\NN_f$ does not lie in $\RR_f$. Such a point
$\hat w$ is an orbit which stays in $\boundary W$ for all time, and in
particular is on the Julia set. If $w_0$ had a neighborhood $D_0$
which pulled back along $\hat w$ eventually univalently, then by the
Shrinking Lemma (after possibly trimming $D_0$ to a slightly smaller
disk), $\diam (D_{-n})\to 0$. However $D_0\intersect W$ is
being pulled back by the univalent map $f|_W$, and so the diameter of
$D_{-n}\intersect W$ cannot shrink.

We shall adopt the convention of using rotation domain, Siegel disk or
Herman ring, to refer also to the leaves of $\RR_f$ which are
invariant lifts of these domains. 

Except in the case of rotation domains, the
structure of a leaf reflects the behavior of $f$ at {\em small scales}
-- this is another consequence of the Shrinking Lemma.  The following
two lemmas show that, barring the obvious exception, all leaves are
topologically trivial.
\begin{lemma}{leaves are planes}
All leaves of $\RR_f$ which are not Herman rings are simply
connected. 
\end{lemma}

\begin{pf}
The invariant lift of a Siegel disk is, by the above discussion, a disk, so we
may from now on consider a leaf $L$ which is not 
either kind of rotation domain. That is, for $\hat z\in L$
there is some $n$ for which $z_{-n}$ is not in a rotation domain. 

Let $\hat\gamma:S^1\to L$ be a simple closed smooth curve on $L$, 
which does not pass through the branched points of $\pi$.  We need to show
that $\hat\gamma$ bounds a disk on the leaf $L$.
Let us consider  the  corresponding sequence of smooth curves on
the Riemann sphere: 
$\gamma_{-n}=\pi_{-n} \circ \hat\gamma: S^1\rightarrow \bar\C$. 
Deforming $\hat\gamma$ slightly, we can get $\gamma_0$ to have only
finitely many  points of self intersection, all of which are double points.
Clearly, the $\gamma_{-n}$ have no more points of self intersection
than $\gamma_0$, since 
if $\gamma_{-n}(a)$ is a simple point for some $a\in S^1$, 
so is $\gamma_{-(n+1)}(a)$. 

Let us  now consider a point of self intersection,
$\gamma_0(a)=\gamma_0(b)$, where
$a, b\in S^1$, and $a\not= b$. Since $\hat\gamma(a)\not=\hat\gamma(b)$,
there is an $n_0$ such that $\gamma_{-n}(a)\not=\gamma_{-n}(b)$
 for $n\geq n_0$, so that
$\gamma_{-n}$ has strictly fewer points of self intersection than $\gamma_0$.
It follows that eventually all the curves $\gamma_{-n}$ are simple. 

Furthermore, by the Shrinking lemma, $\diam \gamma_{-n}\to 0$ as $n\to\infty$.
Let $D_{-n}$ be the component of $\C\ssm \gamma_{-n}$ of small diameter.
Then it contains at most one critical point of $f$ for $n$ sufficiently big.
If $D_{-(n+1)}$ actually contained a critical point,
the curve
$\gamma_{-(n+1)}$ (obtained by analytic continuation of $f^{-1}$ along the
simple curve $\gamma_{-n}$)
would not be closed. Hence the $D_{-n}$ eventually don't contain 
the critical points.

It follows that  the maps $f: D_{-(n+1)}\rightarrow D_{-n}$ are univalent for $n$
 sufficiently big: $n\geq N$.  Hence the set $\hat D$ of backward orbits
$\{(z_{-n}): z_{-n}\in D_{-n}\;\text{for}\; n\ge N\}$
represents a topological   disc in $L$ 
bounded by $\hat\gamma$ (with a homeomorphic projection 
$\pi_{-N}: \hat D\rightarrow D_{-N}$).
\end{pf}

The following lemma excludes elliptic leaves (that is, conformal 
spheres). 
 
\begin{lemma}{not elliptic} If $\deg f >1$, there are no compact
leaves in the lamination $\RR$.
\end{lemma}

\begin{pf} Assume that a leaf $L$ is compact. Then the projection
$\pi: L\rightarrow \bar\C$ is a finite-sheeted branched covering.
However,  we can also express $\pi$ as $f^n\circ \pi \circ \hat
f^{-n}$, so $\deg \pi \ge (\deg f)^n$ for any $n$. This is a contradiction.
\end{pf}

\subsection{Criteria for regularity}
Let us consider some cases in which we can say which part of $\NN_f$
is regular. 

\subsection*{Axiom A case} (See Appendix 2 for definitions) We will
call these functions ``Axiom A" 
instead of the more common ``hyperbolic" in order to avoid sentences
like ``in the hyperbolic case all leaves are parabolic". 

If $f$ satisfies axiom A then 
${\cal R}_f={\cal N}_f\ssm\{\text{finite set of points}\}$, namely the 
attracting cycles of $\hat f$. Note that the backward orbits like
$(\alpha,\ldots, 
\alpha,\beta,\ldots)$, where  $\alpha$ is an attracting fixed point
and $\beta\ne \alpha$ is another preimage, are included into $\RR_f$,
since $\beta\not\in \omega(C)$.

\subsection*{Critically non-recurrent case} 
We will use the notation $\alpha(\hat z)\subset\bar\C$ for the limit set of the
backward orbit  $\hat z= (z_{-n})_{n>0}$. 

\begin{lemma}{nonrec1}
Let $\hat z=(z_{-n})\in \NN_f$ be a backward orbit satisfying the property
that for some $N$, $z_{-N}$ does not belong to  an attracting or parabolic
cycle, nor to the $\omega$-limit set of a recurrent critical point.
Then $\hat z\in \RR_f$.
\end{lemma}

\begin{pf} Let $C_1$ be the set of critical points such that for $c\in C_1$,
$z_{-n}\in \omega(c)$, $n=0,1,\ldots$, and $C_2$ be the complementary set
of critical points. Without loss of generality we can assume that already
$z_0$  does not belong to  an attracting or parabolic
cycle, nor to the closure $\text{cl}(\text{orb}(c))$ for any
$c\in C_2$.

By the assumption, $C_1$ consists of non-recurrent
points. Hence there is an $\epsilon>0$ such that 
dist$(z_{-n}, C_1)\geq\epsilon,$ $n=0,1,\ldots$. For $\delta>0$ let
$U_0=D(z_0,\delta)$,  
and $U_{-n}$ be the pull-back of $U_0$ along $z_{-n}$. By 
Ma\~n\'e's Theorem (\cite{mane:fatou} and \S\ref{dynamics for
geometers}), there is  
a $\delta>0$ such that diam$\;U_{-n}<\epsilon$. Hence $U_{-n}$ does not
hit the critical points of $C_1$. 

Moreover, if $\delta$ is sufficiently small, the orbits of the critical
points $c\in C_2$ clearly don't meet $U_0$. Hence $U_{-n}$ don't hit these
critical points either, so that the pull-back $\{U_{-n}\}$ is univalent. 
\end {pf}

When we refer to an attracting/parabolic etc.
cycle in $\NN_f$, we mean the invariant lift of the corresponding
cycle in $\bar \C$.  Let us recall from Appendix 2 that $C_r$ denotes
the set of recurrent critical points in the Julia set.

\begin{lemma}{closure} The closure of the set $\NN_f\ssm \RR_f$
 of irregular points in $\NN_f$ coincides with the invariant lift
$\hat\omega(C_r)$ together with attracting and parabolic cycles. 
\end{lemma}

\begin{pf} If $\hat z\not\in \hat\omega(C_r)$, nor is an attracting or
parabolic periodic point, then 
it follows from lemma \ref{nonrec1} that
$B(D, \hat z)\subset \RR_f$ for sufficently small neighborhood $D\ni z$.
Thus $\hat z\in \inter \RR_f$.

Vice versa, let  $\hat z\in \hat\omega(C_r)$. Let  $D$ be a neighborhood 
of $z_0$, $N>0$ be any integer, and
$B_0=\cl B( D, \hat z, N)$ be a closed neighborhood of $\hat z$.
We should  show that $B_0$ contains an irregular point.

Since $D_{-N}\cap \omega(C_r)\not=\emptyset$, there is a critical point
 $c\in C_r$ such that $f^{n_1} c\in D_{-N}$ for some $n_1>0$.
Let  $\hat z^{(1)}$ be any backward orbit with $z^{(1)}_{-(N+n_1)}=c$,
and $$B_1\equiv \cl B( D, \hat z^{(1)}, N+n_1)\subset B_0.$$
Then all leaves of $B_1$ over $D$ are at
least double branched.

Let us now consider a neighborhood base $D\equiv D^1\supset D^2\supset\ldots$
of $z$,  and let $D^2_{-(N+n_1)}\ni c$ be the pullback of $D^2$ along the
orbit $\{f^k c\}_{k=0}^{N+n_1}$.
Since $c$ is recurrent, there is an $n_2$ such that 
$f^{n_2}c\in D^2_{-(N+n_1)}$.
Take any backward orbit $\hat z^{(2)}$ with
 $z^{(2)}_{-(N+n_1+n_2)}=c$, and cosider
the closed box
 $B_2=\cl B( D^2, \hat z^{(2)}, N+n_1+n_2)\subset B_1$. All leaves
of $B_2$ over $D^2$ are at least triple branched.

Proceeding in this way, we will construct a nest
$B_0 \supset B_1 \supset B_2\supset\ldots$ of closed boxes, such that
all leaves of $B_n$ are  at least $n$ times branched over $D^n$. Hence
the intersection of these boxes consist of irregular points.
\end{pf}

Let us call a map $f$ {\it critically  non-recurrent} 
if all its critical points on the
Julia set are non-recurrent. The following fact was proved by
Carleson, Jones and Yoccoz \cite{carleson-jones-yoccoz}
(in different language).

\begin{corollary}{nonrec2}
A map $f$ is  critically non-recurrent if and only if
 $$ \RR_f=\NN_f\ssm \{\text{attracting and  parabolic cycles}\}.$$
\end{corollary}

Let us call a map $f$ {\it persistently recurrent} if
any backward orbit $U_0, U_{-1}, \ldots$ of a neighborhood $U_0$
along $\omega(C_r)$ hits a critical point.
 In other words, all points of $\hat\omega(C_r)$
are irregular. Lemma \ref{closure} also yields the
following criterion of openness of the regular leaf space.

\begin{corollary}{persistent} The regular leaf space $\RR_f$ is open
in $\NN_f$ if and only if $f$ is either critically 
non-recurrent or persistently
recurrent. In the latter case
 $$\RR_f=\NN_f\ssm
     (\hat\omega(C_r)\cup \text{parabolic and attracting cycles})$$.
\end{corollary}

\subsection{The Julia and Fatou sets.}
\label{Julia and Fatou sets}
Let us consider the pull-backs $\jr_f\equiv \jr=\pi^{-1}J\cap \RR_f$ 
and $\fr_f\equiv \fr=\pi^{-1}F$
of the Julia set $J$ and the Fatou set $F$
to the space $\RR_f$.

Note first that $\fr$ is obtained from 
the pullback of $F$ to $\NN_f$ just by removing the 
attracting cycles. Also, if we remove from  $\fr$ the
invariant lifts of Siegel disks
and Herman rings, then we obtain a Riemann surface lamination.
Indeed, if $U$ is compactly contained in the Fatou set, and a backward
trajectory $U_0,U_{-1},\ldots$ eventually does not meet 
either attracting cycles, Siegel disks or Herman rings, then
there is an $N$ such that $f^{-k}U_{-N}$ does not meet the critical points
for $k\ge 0$. In particular, the boxes $B(U_0,\hat z,N)$ have a product
structure if $\hat z\in \fr\setminus\{\text{Siegel disks and Herman
rings}\}$ and $N$ is large. 

Note further that $\hat f$ acts properly
discontinuously on $\fr$ with Siegel disks and Herman
rings removed.  Indeed for any $\hat z\in \fr$ which is not in a rotation
domain, either $z_{-n}$ lie in an attracting or parabolic basin
and pull back toward its boundary, or eventually end up in preimages of
a periodic  domain. Thus there is a neighborhood $V$ of $z_{-N}$ for
some $N$ such that all further pullbacks of $V$ accumulate onto $J$.
It follows that $\hat f^{-n}B(V,z_{-N})$ eventually escapes every
compact subset of $\fr$.

Thus, $\fr/\hat f$ is a Hausdorff topological space, and in fact a
Riemann surface lamination, since it inherits its local structure from
$\RR$.

So to each basin of the Fatou set we can associate a Riemann surface 
lamination. These play the role of the Riemann surfaces associated
to a Kleinian group. 

In \cite{sullivan:universalities,sullivan:bounds} Sullivan considered
the natural extension of the attracting 
basin of infinity for a polynomial, and obtained a ``solenoidal
Riemann surface lamination'', called $\SS$ (see Appendix 1). A similar
object appears as a
subset of $\fr/\hat f$ in general.
Let us consider the topological structure of these laminations in
somewhat greater detail.

\subsection*{Attracting domains.}
Consider a cycle of basins 
$U_1\overset{f}{\rightarrow} \cdots \overset{f}{\rightarrow} U_m
\overset{f}{\rightarrow} U_1$ 
for an attracting (or super-attracting) cycle, and let $\GG$ denote
the subset of $\fr$ consisting of orbits $\hat z$ that are attracted
(in forward time) to this cycle. This sublamination divides naturally
into two pieces: let $\GG_1$ contain orbits which stay in $\cup U_i$
for all time, and let $\GG_2$ consist of orbits which, before some
time, lie outside the $U_i$. 

Suppose that all of the
domains are simply connected.
We claim that $\GG_1/\hat f$
is Sullivan's solenoidal Riemann surface lamination (of appropriate degree),
and $\GG_2/\hat f$ is a finite union of copies of
(plane domain)$\times$(Cantor set), which accumulates onto the solenoidal
part. The full quotient $\GG/\hat f$ is, in particular, compact. 

We can study $\GG_1/\hat f$ by considering just the return map 
$f^m$ to $U_1$, and the quotient of the set of orbits of this map that stay in
$U_1$ for all time. On a neighborhood of $\boundary U_1$, $f^m$
is topologically conjugate to $z\mapsto z^d$ acting on a neighborhood of the
boundary in the unit disk $\D$,
and every orbit in $\GG_1$ accumulates in backward time onto
$\boundary U_1$ (note that we omit the orbit which remains on 
the attracting periodic cycle, since it does not lie in $\RR_f$).
It follows that the quotient $\GG_1/\hat f$ is 
homeomorphic to the quotient of the Fatou domain of  0  (or
$\infty$) for the lamination  of $z\mapsto z^d$, namely Sullivan's
solenoidal Riemann surface lamination. 

Now consider an orbit $\hat z$ which escapes $\union U_i$ in backward
time. Let $\til z$ denote the full orbit
$(\ldots,z_{-1},\BrkOK z_0,z_1,\ldots)$.  There is a finite list
$V_1,\ldots,V_p$ of preimages of $U_1$ such that no $V_i$ contains a
post-critical point, and every full orbit $\til z$ with $\hat z\in
\GG_2$ passes through a unique $V_i$.
 Let $q$ denote the smallest integer for which $z_q \in
\union V_i$.  Since $\GG_2/\hat f$ is just the space of these full
orbits modulo shift, we can identify it with $(\union V_i)\times
\Sigma$, where $\Sigma$ is a Cantor set, so that the $\union V_i$ component
is $z_q$ and the $\Sigma$ component specifies the preimages of
$V_i$ which contain the preimages $z_{q-n}$,
$n=1,2,\ldots$.

It remains to see that the closure of $\GG_2/\hat f$ is in $\GG_1/\hat
f$. Let $A$ denote a fundamental annulus in $U_1$. This is a compact
annulus, surrounding the fixed point of $f^m$, through which every
full orbit of $\GG_1$ passes exactly once (or twice if on the
boundary). Now if we consider $\hat z$ in $\GG_2$, such that $z_q\in
V_i$, we see that $z_{q+N}$ passes
through $A$ where $N>0$ gets larger as $z_q$ approaches $\boundary
V_i$. Thus $\hat z$ is very close to some orbit $\hat w\in \GG_1$
which agrees with $\hat z$ for all moments $n$ where $z_{n}\in \union
U_i$. It follows that $\GG_2/\hat f$ 
accumulates on $\GG_1/\hat f$, and in fact that all of $\GG_1/\hat f$
is obtained this way.

If the domains $U_i$ in the cycle are not simply connected the
topological structure of the quotient is more complicated and we shall
not describe it here. However let us sketch an argument showing that
it is compact. Let $D$ be a small closed disk
around the attracting fixed point for $f^m$ in $U_1$, so that $D$ maps
univalently to $f^m(D)\subset D$. For any orbit $\hat z$ attracted to
the cycle there is a first moment $q\in\Z$ when $z_q$ lies in $D$. 

Let $A=D\setminus int(f^m(D))$; this is the same fundamental annulus
described above.
Let $\DD_1$ denote all orbits $\hat z\in \GG$ for which $z_0\in A$.
Let $\DD_2$ denote all orbits $\hat z\in \GG$ for which $z_0\in f^m(D)$
and $z_{-n}\notin D$ for $n>0$. Then modulo the action of $\hat f$
every orbit is uniquely represented in $\DD_1\union \DD_2$, except for
some identifications on the boundaries. Since both $\DD_1$ and $\DD_2$
are compact, it follows that $\GG/\hat f$ is compact.

\subsection*{Leau (parabolic) domains.}
For a cycle of domains with a parabolic periodic point, 
the quotient of the corresponding lamination is not compact. 
One should think of these as obtained from the solenoidal Riemann
surface laminations by a ``pinching'', but 
we will not try to elaborate on this case in this paper.

\section{The Type Problem and affine structure on the leaves.}
\label{affine structure}

By Lemmas \ref{leaves are planes} and \ref{not elliptic}
every leaf of $\RR_f$ is either a parabolic (affine) or hyperbolic
plane, except possibly for (invariant lifts of) Herman rings, which are
hyperbolic annuli. Siegel disks are the only
example we know of hyperbolic planes. 

\proclaim Type Problem.
Are there any other cases of hyperbolic  leaves 
except Siegel disks and Herman rings?

\subsection{Criteria for parabolicity of leaves}
Let us look at the type problem in some special cases.

\subsection*{Repelling fixed point.} Let $\alpha$ be a repelling fixed point
for $f$ with multiplier $\lambda$, and $\hat \alpha=(\alpha,\alpha,\ldots)$ 
be its invariant lift to $\NN_f$.   Let us consider the invariant leaf
$L(\hat \alpha)=\{\hat z: z_{-n}\to \alpha\}$ through $\hat \alpha$. 
This leaf is parabolic since the quotient of $L\ssm \{\hat\alpha\}$
by the action of $\hat f$ is a torus.
Similar reasoning applies to the case of a repelling periodic point.

\subsection*{Parabolic fixed point.} Let now $\alpha$ be a parabolic fixed
point with combinatorial rotation number $p/q$. Then $f^q$ has $s=ql$
invariant repelling petals $P_i$. Let us consider
the set $L_i=L_i(\hat\alpha)$ consisting of backward orbits 
$\hat z$ such that  the suborbit $z_{-qn},\; n=0,1\ldots ,$ 
eventually lands in $P_i$. 
(Observe that
$\hat\alpha$ itself does not belong to these leaves.)
The map $\hat f$ permutes the leaves
$L_i$ organizing them into cycles of order $q$.
 These leaves are parabolic since their quotients
by the $\hat f^q$-action are ``Ecalle-Voronin cylinders" with infinite modulus
(that is, conformally equivalent to $\C^*$).
The case of parabolic periodic points is treated similarly.

\subsection*{General conditions.}
Let us now give a couple of general conditions for a leaf to be parabolic.
Let $D(z,\epsilon)$ denote the spherical disk of radius $\epsilon$ centered
at $z$, and $\hat D(\hat z,\epsilon)$ denote the component of 
$L(\hat z)\cap \pi^{-1}D(z,\epsilon)$ containing $\hat z$.

\begin{lemma}{big disks} Let a backward orbit 
$\hat z=\{z_0, z_{-1},\ldots\} \in {\cal R}_f\ssm\text{(rotation sets)}$
satisfy the following property.
There is an $\epsilon>0$ and a subsequence $\{n(k)\}$ such that
the disk $ D(z_{-n(k)},\epsilon)$ can be univalently pulled  back along the
rest of the orbit,
$\{z_{-m}\}_{m\ge n(k)}.$ Then the leaf 
$L(\hat z)$ is parabolic.
\end{lemma}

\noindent{\bf Remark.} In terms of the natural extension the
assumption of the 
  lemma means that the 
  $\hat D_{-k}\equiv\hat D(\hat f^{-n(k)}\hat z, \epsilon)$ univalently
  project down to the sphere. 

\begin{pf} Assume without loss of generality that $n(0)=0$.
By the Shrinking Lemma,
$\diam(\pi_{-m}\hat D_0) = \delta(m)\to 0$ as $m\to\infty$. Hence for sufficiently
large $k$ the annulus $\hat D_{-k}\ssm \hat f^{-n(k)}\hat D_0$ is univalently
mapped to an annulus on the sphere containing a round annulus with
outradius $\ep$ and inradius $\delta(n(k))$.
Its modulus can therefore be estimated via
$$\mod(\hat D_{-k}\ssm \hat f^{-n(k)}\hat D_0)\geq {1\over
2\pi}\log(\sin c\ep/\delta(n(k)))\to\infty,$$ 
where the constant $c$ accounts for distortion between spherical and
Euclidean metrics.
This is equal to the modulus of its univalent image,
$\hat A_k = \hat f^{n(k)}(\hat D_{-k}) \setminus \hat D_0$, which is an annulus in
$L(\hat z)$ surrounding $\hat D_0$. Since $\mod(\hat A_k)\to\infty$, 
the leaf $L(\hat z)$ must be parabolic.
\end{pf}

Recall that $C$ denotes the set of critical points of $f$.
 The following is an immediate
consequence of Lemma \ref{big disks}.

\begin{corollary}{away from O}
If a backward orbit 
$\hat z=\{z_0, z_{-1},\ldots\} \in {\cal R}$ does not converge to
 $\omega(C)$, then the leaf $L(\hat z)$ is parabolic.
\end{corollary}

Note that the set $C$ can be replaced here by the set $C_r$ of recurrent 
critical points.

\begin{lemma}{annuli condition} Let  $\hat z \in {\RR_f}$. Assume that for some
sequence $n(k)$ there exist annuli 
$\hat A_{-n(k)}\subset L_{-n(k)}=L(\hat f^{-n(k)}\hat z)$
enclosing $ \hat f^{-n(k)}\hat z$ and a branched point of the projection
$\pi: L_{-n(k)}\rightarrow {\bar \C}$, whose moduli stay away from 0.
Then the leaf $L(\hat z)$ is parabolic.
\end{lemma}

\begin{pf} Let $B_{-n}\subset L_{-n}$ be the set of branched points for the
projection $\pi: L_{-n}\rightarrow \bar\C$. Since every branched point is
represented by a backward orbit finitely many times 
passing through a critical point,
$\hat f^{-1} B_{-n}\supset B_{-(n+1)}$, and moreover for any $\hat c\in B_0$
there is an $n$ such that $\hat f^{-n}\hat c$ is not a branched point any more.
Let $P_n=\hat f^n B
_{-n}\subset L_0$. Then
$P_0\supset P_1\supset P_2\supset\ldots$, and $\cap P_n=\emptyset$.
As $P_0$ is discrete, the sets $P_n$ escape to $\infty$. 
Thus, if the leaf $L_0$ were hyperbolic, then the modulus of any
annulus $R_n$ enclosing $\hat z$ and a point of $P_n$ would tend to 0
as $n\to\infty$, which would contradict our assumption.
\end{pf}

\begin{lemma}{derivative condition} 
Consider a backward orbit $\hat z=\{z_0, z_{-1},\ldots\}
\in {\cal R}$ which does not hit the set ${\omega(C)}$. Assume
that $\|Df^{-n}(z)\|\to 0$ as $n\to\infty$ where $f^{-n}$ is the branch
of the inverse map which sends $z$ to $z_{-n}$, and $\|\cdot\|$ means
the hyperbolic metric in ${\C}\ssm {\omega(C)}$. Then the leaf
$L(\hat z)$ is parabolic.
\end{lemma}

\begin{pf}
Let $L_{-n}=L(\hat f^{-n}\hat z)$. Then
the projection
$$\pi: L_{-n}\ssm \pi^{-1} \omega(C)\rightarrow \bar{\C}\ssm {\omega(C)}$$
is a covering map, and hence a local hyperbolic isometry (with respect to
the corresponding hyperbolic metrics).

Assume now that the leaf $L_0$ is hyperbolic. Then all 
$L_{-n}$ are also hyperbolic.
Since the inclusion $i: L_{-n}\ssm \pi^{-1}\omega(C)\rightarrow L_{-n}$
is a hyperbolic contraction, the projection $\pi$ is expanding from
the hyperbolic metric of  $L_{-n}$ to the hyperbolic metric
of $\bar{\C}\ssm {\omega(C)}$.

Note finally that $\hat f^{-n}: L_0\rightarrow  L_{-n}$ is
a hyperbolic isometry. Hence $\|Df^{-n}(z)\|\geq \|D\pi(\hat z)\|^{-1}>0$
where the last norm is measured from the hyperbolic metric on $L_0$
to the hyperbolic metric of $\bar{\C}\ssm {\omega(C)}$. Contradiction. 
\end{pf}

\bfheading{Remark:} We don't know whether the above contracting property
along the backward orbits is always satisfied (unlike the expansion propery 
along the forward orbits: see McMullen \cite{mcmullen:renormalization}, 
Theorem 3.4). See also Lemma \ref{renorm parabolic} below.

\subsection*{Axiom A case.}   
Let $f$ satify Axiom A.
Let us consider a backward 
orbit $\hat z=\{z_{-n}\}\in \RR_f$. Then this backward orbit converges
to the Julia set, and hence stays bounded distance away from $\omega(C)$.
By Corollary \ref{away from O}, all leaves of $\RR_f$ are parabolic.

\subsection*{Critically non-recurrent case}

\begin{proposition}{nonrec3}
 Assume that all critical points on the Julia set are non-recurrent.
Then $$ \an_f=\RR_f=\NN_f\ssm \{\text{attracting and parabolic cycles}\},$$
so that all regular leaves are parabolic.
\end{proposition}

\begin{pf} The second equality
$\RR_f=\NN_f\ssm \{\text{attracting and  parabolic cycles}\}$
was proved above (Corollary \ref{nonrec2}). 

In order to prove the first one, let us consider the following ordering on
the set of critical points in  $J(f)$:
 $c_1\succ c_2$ if cl(orb$(c_1))\ni c_2$. Given a $\hat z\in \RR_f$,
let $\tilde C$ denote the set of critical points belonging to $\alpha(\hat z)$.

Assume first that $\tilde C\not=\emptyset$.  Then let us take a
critical point $a\in \tilde C$ which is a maximal element of this
ordering.  Let $\epsilon>0$ be such that $z_{-n}$ stay distance at
least $\epsilon$ from all critical points $c\not\in \tilde C$. By
Ma\~n\'e's theorem (Appendix 2), there is a $\delta>0$ such that for
all $n$ all components of $f^{-n}D(a,\delta)$ have diameter at most
$\epsilon$. Hence if $z_{-k}\in D(a,\delta)$, and we pull
$D(a,\delta)$ back along $\{z_{-(k+n)}\}_n$, then we don't hit the
critical points $c\not\in \tilde C$.  Clearly we will not hit the
critical points $c$ of $\tilde C$ either (provided $\delta$ is small
enough), since their forward orbits don't accumulate on $a$. Hence
this pull back is univalent.

Select now a sequence $k(l)$ such that $z_{-k(l)}\to a$, and apply 
Lemma \ref{big disks}. (Note that the lemma applies since there can be
no rotation domains in the non-recurrent case).

If $\tilde C=\emptyset$ then take any point $a\in \alpha(\bar z)$, 
and repeat the
above argument.
\end{pf}

{\bf Remark.} By a minor modification of the above argument one can check
  that {\it the leaf $L(\hat z)$ is parabolic, provided $\hat z$ is not 
   an attracting or parabolic cycle, and
  $\alpha(\hat z)$
is not contained in $\omega(C_r)$}, where $C_r$ is the set of recurrent
critical points.

\subsection*{Invariant measures with positive characteristic exponent.}
Let $\mu$ be an invariant measure of $f$, and suppose that for
$\mu$-a.e. $z$ the characteristic exponent
$$\chi(z)=\lim{1\over n} \log |Df^n (z)|$$
exists and is positive ($|\cdot|$ means the spherical norm). 

Let $\hat\mu$ be the lift of $\mu$ to the natural extension. The 
Pesin local unstable manifolds for $\hat \mu$ are the sets 
$\hat D(\hat z, \epsilon(\hat z))\subset L(\hat z)$ 
which univalently project down
to the sphere and whose backward orbits shrink exponentially.
Moreover, $\epsilon(z)>0$ $\hat\mu$-a.e. 

Let $X_{\epsilon}=\{\hat z: \epsilon(z)>\epsilon\}.$
 It follows from the Poincar\'{e}
recurrence theorem that for $\hat \mu$-a.e. $\hat z$ there is an $\epsilon>0$
such that the backward orbit $\hat f^{-n}\hat z$ infinitely many times
visits $X_{\epsilon}$. By Lemma \ref{big disks} the leaves $L(\hat z)$
are parabolic for  $\hat \mu$-almost all $\hat z$
(compare \cite{bedford-lyubich-smillie}, \cite{wu}).

\subsection*{Infinitely renormalizable quadratics}
We refer to the papers of Douady and Hubbard \cite{douady-hubbard:polylike} and
 McMullen \cite{mcmullen:renormalization}, \cite{mcmullen:renormfiber}
 for the background
in holomorphic renormalization theory. Here we will briefly recall the
basic concepts.

Let $U'$ and  $U$ be two topological disks such that $\cl U'\subset U$.
A double branched covering map $f: U'\rightarrow U$  is called
{\it quadratic-like}. We assume that its critical point is located at the
origin 0. The set $K(f)=\{z: f^n z\in U', n=0,1,\ldots\}$ is 
called the {\it filled Julia set}; its boundary is called the Julia
set $J(f)$.  
The Julia set is connected iff the critical point 0 is non-escaping, that is,
$0\in K(f)$.

Any quadratic polynomial can be viewed as a quadratic-like map with 
$U$ being a round disk of sufficiently large  radius, and $U'$ being its
pullback. By the Straightening Theorem of Douady and Hubbard any 
quadratic-like map $f: U'\rightarrow U$ 
is quasi-conformally conjugate to some quadratic 
polynomial $z\mapsto z^2+c$. Moreover, if $\mod(U'\ssm U)\geq\epsilon>0$
then there is a conjugacy with dilatation bounded by $K(\epsilon)$.  

We can specify a distinguished fixed point of $f$ as follows: 
Take an arc $\gamma\subset U\ssm K$ with endpoints $a$ and $f(a)$. 
Choosing appropriate pullbacks of this arc
by $f$ we obtain a curve $\Gamma\supset \gamma$ such that
$f(\Gamma\setminus\gamma)=\Gamma$. It turns out that if $J(f)$ is connencted
then this curve lands at a 
specific fixed point of $f$, usually denoted by $\beta$. This point is
repelling for any quadratic polynomial $z\mapsto z^2+c$ except $c=1/4$.

A quadratic-like map $f$ is called {\it renormalizable}  under the following
circumstances:
\begin{itemize}
\item
Some iterate $g=f^p$ $(p>1)$ restricted to
  an appropriate topological disk $U'\ni 0$ is quadratic-like with connected
  Julia set;
\item
The sets $f^k K(f)$, $k=1,\ldots, p-1$, do not 
  touch $K(f)$ except perhapsfor the $\beta$-fixed point of $g$. 
\smallskip
\end{itemize}
Under these circumstances the map $g$ 
is called a {\it renormalization } of $f$. If there is a sequence 
of renormalizations $g_n: U_n'\rightarrow U_n$ with increasing periods $p_n$,
the map $f$ is called {\it infinitely renormalizable}. If this sequence can
be selected in such a way that  the ratios
$p_{n+1}/p_n$ are bounded,  then one says that $f$ is of {\it  bounded type.} 

 We say that $f$ is
an infinitely renormalizable map with {\it a priori bounds} if
there is a sequence of renormalizations as above and an  $\epsilon>0$ such that
$$\mod (U_n\setminus U_n')\geq\epsilon,\; n=0,1,\ldots$$

Let us say that a map is {\it Feigenbaum-like} 
if it is infinitely renormalizable
of bounded type with  a priori bounds. Any infinitely renormalizable
real quadratic of bounded type is {\it Feigenbaum-like}: complex a priori
bounds were established by Sullivan (see \cite{sullivan:bounds},
 \cite{demelo-vanstrien}).

For a Feigenbaum-like map the set $\omega_f(0)$ is a Cantor set of bounded
geometry, and $f|\omega_f(0)$ 
is an invertible 
minimal dynamical system (conjugate to a translation on a group).
In particular, $f$ is persistently recurrent and hence, by
Corollary  \ref{persistent}, 
$\RR_f=\NN_f\ssm\hat \omega_f(0)$.

\begin{lemma}{renorm parabolic}
Let $f$ be a Feigenbaum-like  quadratic
polynomial. Then all leaves 
of the lamination ${\cal R}_f$ are parabolic.
\end{lemma}

\begin{pf} 
 Let $\hat z=\{z_0, z_{-1},\ldots\}\in {\cal R}$.
Then this orbit eventually stays out of the set
$\omega(0)$, so we can assume that $z_0\in {\C}\ssm \omega(0)$.

Let $g: U'\rightarrow U$ 
be  a renormalized map, mod$(U\ssm U')>\epsilon>0$.
Let $\beta$ be its distinguished fixed point, as above.
Since $f$ is of bounded type, $g$ is $K(\epsilon)$-quasi-conformally
 conjugate to
a polynomial $z\mapsto z^2+c$ with $|c-1/4|>\delta>0$ (with $\delta $ depending
on the type and a-priori bound).  Hence  
\begin{equation}
\label{Dg bound}
\|Dg(\beta)\|\geq \lambda>1.
\end{equation}

Because of Corollary \ref{away from  O}, we can assume that
$z_{-n}$ converge to $\omega(0)$.
Let $\omega_g(0)$ be the closure of the
postcritical set of $g$. Then there is a backward orbit $\zeta_{-l}$
of $g$ converging
to $\omega_g(0)$ which is a part of the backward orbit $\hat z$.
 Let us take the
second element  $\zeta_{-1}$ of this backward orbit. 
Clearly $\zeta_{-1}\in U'\ssm U''$, where $U''$ is the $g$-pullback of $U'$.

The set of Feigenbaum-like maps is compact in the
Caratheodory topology (see McMullen \cite{mcmullen:renormalization}).
Hence there is a path $\gamma$ in
$U\ssm\omega_g (0)$ 
joining $\zeta_{-1}$ and $\beta$ of bounded hyperbolic length, such
that the analytic continuation of $g^{-1}$ which fixes $\beta$ carries
$\zeta_{-1}$ to $\zeta_{-2}$. It follows from this and (\ref{Dg bound}) that
$\|Dg(\zeta_{-2})\|\geq\theta>1$ for some $\theta$ depending on $\epsilon$,
$\lambda$ and the hyperbolic length of $\gamma$ only.

Hence 
$\|Df^{-n}(z)\|\to 0$ as $n\to\infty$ (where $f^{-n}z=z_{-n}$),
 and Lemma \ref{derivative  condition} yields the desired result.
\end {pf}

\bfheading{ Remark:} The above way to get hyperbolic contraction is, 
modulo the details, 
due to Curt McMullen (compare \cite{mcmullen:renormfiber}, Proposition 5.9). 
It is actually possible to 
weaken the assumptions of the lemma: 
McMullen has an argument showing that his notion 
of ``robustness'' suffices to give the desired contraction.

\subsection{Affine structures on the leaves and linearization.}
\label{affine coordinate}
Being unable to resolve
the type problem in full generality, let us define a new leaf space
$\an_f$ by throwing away from  ${\cal R}_f$ all hyperbolic leaves.
All leaves in $\an_f$ are conformally equivalent to the complex plane
${\C}$ and hence possess a unique affine structure compatible with their
conformal structure. 

We can express this affine structure as a limit
of rescalings of backward branches of $f$:

\begin{lemma}{formula} 
Let $f$ be a rational map with $\infty$ a critical point.  Given
$\hat z\in \an_f$, there exists a sequence of similarities
$A_n(w) = \alpha_n w + \beta_n$ such that
the maps
$$ \phi_n = A_n\circ \pi \circ \hat f^{-n}: L(\hat z)\rightarrow {\C}$$ 
converge (uniformly on compact sets) to a conformal isomorphism
$\phi:L(z)\rightarrow {\C}$. 
\end{lemma}

\bfheading{Remark:} The condition that $\infty$ is critical can always
be arranged by conjugation with an appropriate M\"obius
transformation. For polynomials it is automatic.

\begin{pf}
Take a disk neighborhood $\hat U=(U_0, U_{-1},\ldots)$
of $\hat z$ in $L(\hat z)$ with compact closure.
Since the leaf $L(\hat z)$ is parabolic, for any $M>0$
$\hat U$ is contained in a disk $\hat V = \hat V(M)$ with modulus
mod($\hat V\ssm \hat U) =  M$.

Let $l = l(M)$ be such that $\pi_{-n} = \pi \circ \hat f^{-n}$
is univalent on $\hat V$
for $n\ge l$ (possible by definition of $\RR_f$).
Thus for $n>l$, $V_{-n}=\pi_{-n}(\hat V)$ contains no critical points and in
particular lies in $\C$. Choose the similarity $A_n$ such that
$A_n(z_{-n}) = 0$ and $A'_n(z_{-n}) = (\pi'_{-n}(\hat z))^{-1}$.
Therefore $\phi_n = A_n\circ\pi_{-n}$ have been normalized by 
$\phi_n(\hat z) = 0$, $\phi'_n(\hat z) = 1$.  (In order for these
derivatives to make sense on $L$ we should fix some local coordinate
chart).

For $n>l$ and $k>0$, we can write $\phi_{n+k} = \phi_n \circ G_{n,k}$,
where 
$G_{n,k} = A_{n+k} \circ f_n^{-k} \circ A_n^{-1}$, with $f_n^{-k}$
denoting the branch of $f^{-k}$ taking $V_{-n}$ to $V_{-n-k}$.
Note that $G_{n,k}$ is defined and 
univalent on $A_{n}(V_{-n}) = \phi_n(\hat V)$, and is normalized so
that $G_{n,k}(0) = 0$, $G'_{n,k}(0) = 1$. 
By the Koebe 1/4 theorem $\phi_n(\hat U)$ contains a disk of definite
radius $\delta>0$. Since $\phi_n$ is univalent the modulus of
$\phi_n(\hat V)\setminus \phi_n(\hat U)$ is  $M$, so by the
Koebe distortion theorem (see appendix 2) the 
nonlinearity of $G_{n,k}$ on the $\delta$-disk around 0 is small, and
goes to 0 as $M\to\infty$, independently of $k$.

Letting $M$, and therefore $l$ and $n$, go to $\infty$, 
it follows that $G_{n,k} \to \text{id}$ uniformly on a $\delta/2$ disk
around $0$ as $n\to\infty$. Thus $\{\phi_n\}$ form a Cauchy sequence,
and so converge uniformly on a neighborhood of $\hat z$. Since
$\{\phi_n|_{\hat U}\}$ is a normal family, they must converge on all
of $\hat U$. 

Applying this argument to a sequence of disks $\hat U_m$ exhausting
$L(\hat z)$, we conclude that $\phi_n$
converge uniformly on compact sets to a global map $\phi:L(\hat z)\to
\C$, which is univalent. Since $L(\hat z)$ is parabolic its image
must be all of $\C$, so $\phi$ is an isomorphism.
\end{pf}

In the particular case where $\pi$ is already univalent on a
leafwise neighborhood of $\hat z\in\RR_f$ (i.e. no $z_{-n}$ is a
critical point for 
$n>0$), we can identify this neighborhood with a neighborhood
of $z_0$ and obtain this local formula for the affine chart:
\begin{equation}
\label{chart}
\phi_{\hat z}(\hat\zeta)=\lim_{n\to\infty} (f^n)'(z_{-n})(\zeta_{-n}-z_{-n})
\end{equation}
(if $f$ is appropriately normalized, e.g., if $\infty$ is critical).
In the case of a leaf corresponding to a repelling 
fixed point this exactly corresponds to the classical formulas 
for the linearizing coordinate.  Note however that {\em uniform}
expansion is not necessary for this formula to hold. 

Namely, if $\alpha$ is a repelling fixed point, then
the affine map $\phi: L(\hat\alpha)\rightarrow \C$
is given by the classical K\"{o}nigs linearizing function:
$$\phi(\hat \zeta)=\lim \lambda^{-n} (\zeta_{-n}-\alpha). $$
Note that $\phi(\hat f^{-1}\hat \zeta)=\lambda^{-1}\phi(\hat\zeta),$
$\zeta\in L(\hat\alpha)$, so that $\phi$ conjugates $\hat f^{-1}$ on the leaf
to the linear map $z\mapsto \lambda^{-1}z$.

Let now $\alpha$ be a parabolic fixed point with combinatorial rotation
number $p/q$. An explicit
affine map from the associated leaves  $L_i(\hat\alpha)$ to $\C$ 
is given by the Leau-Fatou linearizing function:
$$\phi(\hat\zeta)=\lim (h({1\over (\zeta_{-nq}-\alpha)^s})-n),  $$
where $s$ is the number of petals at $\alpha$, and 
$h$ is an appropriate local chart at a sectorial region at $\infty$
(compare Milnor \cite{milnor:intro}, \S 7) .
This function conjugates $\hat f^{-q}$ on the leaf
to a translation $z\mapsto z+a$.  This corresponds to a variation on the
construction in Lemma \ref{formula}, where the rescaling map $A_n$ is
precomposed with a fixed local chart in $\bar\C$, in this case
$w\mapsto h(1/(w-\alpha)^s)$.

In general, affine structure on the leaves
of $\an_f$ can be viewed as a {\it simultanuous linearization} of 
the dynamics along the backward orbits. Indeed, $\hat f$ becomes an affine map
between the leaves. In the affine local charts (\ref{chart}) these maps
become just multiplications by the derivative at the base point:

\begin{equation}
\label{linear}
\phi_{\hat f\hat z}(\hat f\hat\zeta)= f'(z)\cdot \phi_{\hat z}(\hat\zeta),
\end{equation}
provided no $z_{-n}$ is critical for $n>0$. 

\subsection{Density of leaves} 
Let us say that a leaf space  $X$ is {\it minimal}
if all leaves are dense in $X$.

\begin{lemma}{affine leaves are dense}
Any parabolic leaf $L$ is
dense in $\NN_f$. Thus the leaf space  $\an_f$ is minimal.
Moreover, $\jr\cap L$ is dense in the pullback $\pi^{-1}J$ of the Julia
set to $\NN_f$.
\end{lemma}
\begin{pf}
Since $\pi_{-n}$ is a non-constant analytic map on the parabolic leaf
$L$, it can miss at most two points in $\bar\C$.  Now consider any
$\hat z\in\NN_f$, and large $n>0$. Since $\pi_{-n}(L)$ is dense, there
is some $\hat w\in L$ with $w_{-n}$ as close as we like to
$z_{-n}$. If it is sufficiently close then the spherical
$\dist(w_{-j},z_{-j})$ will be small for all $0\le j \le n$. Thus $L$
is dense.  If $\hat z\in \pi^{-1}J$ then clearly $\hat w$ can be
selected from $\jr$.
\end{pf}

We remark that it seems plausible that $L$ is dense even if it 
is hyperbolic, provided that it is not a rotation domain.

Let $\jn_f$ denote $\jr_f \intersect \an_f$, the Julia set in the
affine leaf space. 
\begin{corollary}{Julia compact}
The Julia set $\jn_f\subset\an_f$ is compact if and only if $f$ is
critically non-recurrent.
\end{corollary}

\begin{pf} If $f$ is critically non-recurrent then by Proposition \ref{nonrec3}
$\jn_f = \jr_f = \pi^{-1}(J)$, which is 
a closed subset of $\NN_f$, and thus compact.

Otherwise, by Lemma \ref{closure}, $f$ has irregular points  on 
$\pi^{-1}J$. On the other hand, by Lemma \ref{affine leaves are dense},
$\JJ=\jn_f$ is dense in  $\pi^{-1}J$. Hence $\JJ$ is not closed in $\NN_f$,
thus not compact.
\end{pf}    

\subsection{Local leaves on a global leaf.}
Let $\hat\alpha$ be a repelling periodic point, and $L=L(\hat\alpha)$
the leaf of $\hat\alpha$. Let $D\subset\bar\C$ be a topological disk
which does not 
contain $\alpha = \pi(\hat\alpha)$, and let $\Delta$ be a topological
disk compactly contained in $D$. Let $\hat D_i$ be the connected
components of $\pi^{-1}(D)\intersect L$ which univalently project down
onto $D$, and let $\hat\Delta_i\subset \hat D_i$ be the corresponding
components of $\pi^{-1}(\Delta)\intersect L$.

The following lemma is a ``natural extension'' of the Shrinking Lemma
(see Appendix), and will be applied in \S\ref{pcf orbifolds}.

\begin{lemma}{rel dist}
The size of the $\hat \Delta_i$ shrinks relative to their distance to
$\hat\alpha$: 
\begin{equation}
\label{relative size}
{\diam_L \hat\Delta_i \over \dist_L(\hat\Delta_i,\hat\alpha)}
\to 0, \ \ \text{as $i\to\infty$,}
\end{equation}
where $\diam_L$ and $\dist_L$ are measured in any uniformizing chart
$\phi:L\to\C$. 
\end{lemma}

\bfheading{Remarks.} 1. Clearly the ratio in (\ref{relative size}) does not
depend on the choice of uniformizing map $\phi$. 

2. The result is still valid if we take all
components $\hat D_i$ with some uniform bound on their branching over $D$.

\begin{pf}
Clearly we can assume that $\hat\alpha$ is fixed. Note then that the
$\hat\Delta_i$ escape to infinity in $L$ (that is, eventually don't
intersect any given leaf-compact subset in $L$), since they have
disjoint collars $\hat D_i\setminus\hat \Delta_i$ of definite modulus.

Let $\hat U\subset L$ be a leafwise neighborhood of $\hat\alpha$ such
that $\cl \hat U\subset \hat f (\hat U)$, and $\hat f\hat U$ is
univalently projected down onto $U\subset \bar\C$. Let $n_i$ be the
first positive integer such that $\hat f^{-n_i}(\hat\Delta_i) \subset
\hat U$. As $\Delta_i$ escape to infinity in $L$, $n_i$ goes to
$\infty$. 

The disks $\pi \hat f^{-n_i}(\hat\Delta_i)$ are univalent
pullbacks of $\Delta$ which are not contained in a rotation domain, so
by the Shrinking lemma their (spherical) 
diameters go to 0 as $i\to\infty$. 
As $\pi|_{\hat U}$ has bounded distortion (from the affine structure
on $\hat U$ to the spherical structure on $U$), we have
$$
{\diam_L (\hat f^{-n_i}\hat\Delta_i) \over
\dist_L (\hat f^{-n_i}\hat\Delta_i,\hat\alpha)} \to 0 
\ \ \text{as $i\to\infty$.}
$$
But since $\hat f$ preserves the affine structure on $L$, the ratio in
the last equation is equal to the ratio in (\ref{relative size}).
\end{pf}

\section{Post-critically finite maps.}
\label{pcf orbifolds}
\newcommand{\topl}{\tau_{\bold{\ell}}}
\newcommand{\topn}{\tau_{\bold{n}}}

The affine leaf space $\an_f$ which we have constructed so far is not, in
general, a lamination. The missing ingredients are both topological --
the lack of a local product structure -- and geometric -- non-continuity
of the affine structures in the transverse direction, even where there
is a product structure. As we shall see, these two problems are
related. 

In this section we will give an explicit rearrangement of $\an_f$ --
a change of topology and the addition of new leaves -- in the special
case of post-critically finite rational maps. This should serve as a
motivating example, an indication of the kind of structure that
arises, and a demonstration of how orbifold leaves appear in a natural
way.

In \S\ref{universal construction}, we will give a completely general
construction of an affine orbifold lamination for any rational map, 
emerging naturally from the affine  group action on a
space of meromorphic functions. Thus one could read that section
without first reading this one, but the reader may find that the
explicit examples given here help to illuminate the more abstract approach.

\medskip

We will first construct the orbifold lamination topologically, and
then discuss continuity of affine structures. 

\subsection{Topological orbifold lamination.}
To fix ideas, assume for the
moment that $f$ is a post-critically finite quadratic polynomial,
and moreover that the critical point is actually pre-fixed:
there is a fixed point $\alpha$ such that $f^l c=\alpha$ for some $l>1$. (We
will discuss the general postcritically finite case in \S\ref{general
pcf case}). It is standard that
$\alpha$ is a repelling fixed point (see discussion at the end of
\S\ref{general pcf case}).

As usual, let $\hat\alpha=\{\alpha,\alpha,\ldots\}$ denote the
invariant lift of $\alpha$ to $\NN_f$, and
let $L\equiv L(\hat\alpha)$ denote the $\hat f$-invariant leaf of 
$\hat\alpha$ in $\NN_f$.

Recall that $\an_f$ and $\RR_f$ are both equal to $\NN_f \ssm \hat\infty$
(Lemma \ref{nonrec2}). Our orbifold lamination $\ac_f$ will consist of
$\an_f$, with the leaf $L$ replaced by two copies named
$L^r$ and $L^s$. The topology $\topl$, and orbifold structure, are
described as follows.

Let $q:\ac_f \to \RR_f$ be the map that re-identifies $L^r$ and
$L^s$. Let us consider the pull-back topology 
$q^{-1} \topn$, where $\topn$ is the natural topology of
$\RR$ as a subset of $\NN$. Note that $\hat f$ is naturally lifted to a 
homeomorphism $\tilde f$ of $\ac_f$ with this topology.
However the pull-back topology is not Hausdorff since it does not
separate the leaves $L^r$ and $L^s$.
The actual topology $\topl$
 will be the minimal strengthening of the pullback 
topology  $q^{-1}\topn$, which separates these leaves, keeps 
$\tilde f$ as a homeomorphism, and gives $\ac_f$ the structure of an
orbifold lamination.

Let $\hat z=(z_0,z_{-1},\ldots)\in \NN_f$,
$D$ be a topological disk containing $z_0$ and at most
one postcritical point $f^k c$, $1\leq k\leq l$, and let $N$ be a
natural number. 
Let $B(D, \hat z,  N)$ and
$B(D,z_0)=\pi^{-1} D \equiv B(D, \hat z, 0)$
be the $\topn$ box neighborhoods defined as in
(\ref{box}), and recall that $D_0,D_{-1},\ldots$ are the pullbacks of $D$
along $\hat z$.

If $D_{-N}$ does not intersect the postcritical set then $B(D, \hat z, N)$
 has a natural product structure
$T\times D_{-N}$. Moreover, the projection 
$\pi: B(D, \hat z, N)\rightarrow D $ is either univalent 
or two-to-one branched covering on all leaves. The latter occurs when
$D$ contains a postcritical point $f^ k c$, and then 
this point  is the projection of
the branched point on any leaf. 

This situation always occurs if $\hat z\not=\hat \alpha$, and $N$ is 
sufficiently high. It is more complicated for $\hat z=\hat \alpha$.
In this case some of the
leaves are univalent and some are branched, so that 
$B(D,\hat\alpha, N)$ does not have a natural box structure.

Let us call a backward orbit $\hat z$ with $z_0 = \alpha$ {\em
singular} if it contains $c$ (i.e. it is a branch point of $\pi$), and
{\em regular} if it does not contain $c$, and is not equal to 
$\hat \alpha$.

Given a topological disk $D\ni\alpha$ (not containing other postcritical
points),
let $B^r (D, \hat\alpha, N)$
consist of the union of local
leaves in $B(D,\hat\alpha, N)$ containing regular orbits, and
$B^s (D, \hat\alpha, N)$ be the union of
local leaves containing singular orbits. These are disjoint open sets
in $\NN_f$ with a natural product structure. Moreover, together with the
local leaf $\hat D=\hat D(\hat\alpha)$
containing the fixed point $\hat\alpha$,
they make up all of $B(D,\hat\alpha, N)$.
We set $B^{\mu}(D)\equiv B^{\mu}(D,\hat\alpha, 0)$, where $\mu$ stands for
$r$ or $s$.

Given a set $X\subset \NN_f$, let $\tilde X\subset \ac_f$ denote $q^{-1} X$.
Let also  $\tilde D^{\mu}$ denote the
component of $q^{-1} (\hat D)$ lying in the corresponding leaf $L^{\mu}$.
The similar meaning is given to a point $\tilde z^{\mu}\in \tilde D^{\mu}$ 
corresponding to $\hat z\in \hat D$.

Let 
\begin{equation}
\label{Q-box}
Q^{\mu}(D,\hat\alpha, N)
 = \tilde B^{\mu}(D,\hat\alpha, N) \union \til D^{\mu}\quad
{\rm and}\quad Q^{\mu}(D)\equiv Q^{\mu}(D,\hat\alpha, 0). 
\end{equation}
These sets are going to be neighborhood bases for points $\tilde\alpha^{\mu}$.

Let us now define the topology $\topl$ as the minimal strengthening
of the pull-back topology $q^{-1}\topn$ 
for which the sets $\tilde f^n (Q^{\mu}(D))$
 are open for all $n\geq 0$.  

\begin{lemma}{Q topology works}
With the new topology, $\ac_f$ is an orbifold lamination with one
 singular point, $\tilde\alpha^s$. The projection
$q:\ac_f\to \NN_f$ is 
continuous, and $\tilde f$ acts homeomorphically.
\end{lemma}

\begin{pf} Given a $\hat z\in \hat D$ and a topological disk
$\Delta\subset D$ containing $z_0$, let 
\begin{equation}
\label{intersec}
Q^{\mu}(\Delta,\hat z,  N) = Q^{\mu} (D )\cap \tilde B(\Delta,\hat z, N).
\end{equation}
When $\hat z=\hat \alpha$, we go back to the sets $Q^{\mu}(D, \hat\alpha, N)$
introduced above.

Let  $\tilde\BB$ be the family of
all sets  $\tilde B(\Delta,\hat z, N)$ for $\hat z\in\RR_f$ and any disk
$\Delta\ni z_0$. 
Let $\QQ$ be the family of sets
 $Q^{\mu}(\Delta,\hat z, N)$, where $\hat z\in \tilde D^\mu$ and
$\Delta\subset D$ contains $z_0$.
Let $\TT=\tilde\BB\union\bigcup_{n\geq 0}\tilde f^n\QQ$.
We claim that $\TT$ is a neighborhood basis for the topology $\topl$.

All elements of $\TT$ are open in $\topl$, by definition.
We need to check that, for any $U,V$ in $\TT$ and 
$x\in U\intersect V$ there
is some $W\in \TT$ such that  $x\in W \subset U\intersect V$.

Clearly the sets $\tilde B(\Delta,\hat z, N)$ form a basis for the pullback
topology $q^{-1}\topn$. 
Also, restricting $\Delta$ or increasing $N$ without changing other
parameters clearly makes a set from $\TT$ smaller.
Taking additionally into account (\ref{intersec}), we conclude 
that it is enough to check the case when
$U=\tilde f^m Q^{\mu}(D)$ and $V=\tilde f^n Q^{\nu}(D)$ (Here $\mu$ and $\nu$ 
are independently either $r$ or $s$). By pulling back, we may assume that $m=0$.

Assume that $x$ does not belong to  the local leaf $\tilde D^{\mu}$
of $\tilde\alpha^{\mu}$
in $Q^{\mu}(D)$.  Note that
$ \Omega=Q^{\mu}(D)\ssm \tilde D^{\mu}$  is open in the pullback topology.
Hence $\tilde f^{-n}\Omega$ contains a basic $X\ni \tilde f^{-n} x$ of family
$\tilde\BB$. By (\ref{intersec}) $Q^{\nu}(D)\cap X\in \QQ$. Thus
  $\tilde f^n(Q^{\nu}(D)\cap X)\in \TT$ is a desired set $W$.  

Assume now that $x\in \tilde D^{\mu}$. We can select the basis of disks
$D$ in such a way that $\hat f \hat D^{\mu} $
overflows $\hat D^{\mu}$.  Then $\tilde f^{-n} x\in \tilde D^{\mu}$,
and hence  $\nu=\mu$. Moreover, the component $X$  of 
$ Q^{\mu}(D)\cap \hat f^{-n}  Q^{\mu}(D)$ containing $\tilde f^{-n }x$
is just $Q^{\mu}(D,\hat\alpha, n)$, a set of family $\QQ$.
Now the desired statement follows.

It is clear that $\ac_f$ is Hausdorff -- the doubled points have
separating neighborhoods, by the construction. Note that, away from the
postcritical points, the topology has been changed only in the fiber
direction, where a dense set of fibers has been doubled. 

Let us now check that $\tilde f: \ac_f\rightarrow \ac_f$ is a homeomorphism.
 Obviously $\tilde f^{-1}$ is 
continuous. To verify that $\tilde f$ is continuous, it is enough to check that
$\tilde f^{-1} Q^{\mu}(D)$ are open. Let $f^{-1}D=D_0\cup D_1$
where $D_0\ni \alpha$ while $D_1\ni f^{l-1}c$. Then
  $$ \tilde f^{-1}Q^{\mu}(D)= Q^{\mu}(D_0)\cup B(D_1, f^{l-1}c).$$

Finally let us check that all sets of the basis $\TT$ are  orbifold boxes.
Indeed,  all sets $B(D,\hat z, N)\approx D\times T$ are regular
lamination boxes. Hence the sets  $\tilde B(D,\hat z, N)\approx D\times\tilde T$
are also regular boxes with $\tilde T$ obtained from $T$ by
doubling points coresponding to the leaf $L$. 

The sets  $Q^r(D, \hat\alpha, N)$ are also  regular boxes $D\times T$
with the transversal $T$ consisting of all backward orbits
$\alpha,\ldots,\alpha,\ldots$ (at least $N$ $\alpha$'s)
which never pass through $c$.

Let us now consider the set $K$ of singular backward orbits 
$\alpha,\ldots,\alpha,\ldots$ (at least $N$ $\alpha$'s)
together with the point $a=\bar \alpha$ (in the natural extension topology). 
Then $Q^s(D, \hat\alpha, N)$ is homeomorphic to the
 orbifold box with transversal $(K,a)$ described in Example \ref{deg 2 example}.

Finally if  $\hat z\not=\hat\alpha$ and $D\not\ni\alpha$
 then the sets $Q^{\mu}(D, \hat z, N)$ are regular boxes 
with the same transversal $K$.

Thus $\ac_f$ is indeed an orbifold lamination.
\end{pf}

\subsection{Orbifold affine structure.}  By Corollary \ref{nonrec2} 
all leaves of the  lamination $\ac_f$ are parabolic. Let us supply
all leaves except $L^s$ with their unique affine structure. As to the
leaf $L^s$, let us consider a branched double covering 
$p: \Lambda^s\rightarrow L^s$ with a single branched point over 
$\tilde\alpha^s$. Then $\Lambda^s$ is a parabolic plane which hence
has a unique affine structure. Pushing this structure down to $L^s$
we obtain an orbifold affine structure on $L^s$ with one singular
point at $\tilde\alpha^s$.

There is no ambiguity in the above construction as the double covering
$p$ is uniquely defined up to pre- and post-compositions with  affine maps.
So after appropriate selection of the affine coordinates $z$ and
$\zeta$  on $L^s$ and $\Lambda^s$ correspondingly, $p$ just
becomes the quadratic map $z=\zeta^2$. But $z$ is a linearizing
coordinate on the leaf $L^s$ (see \S \ref{affine coordinate}).
Thus the {\it orbifold affine coordinate $\zeta$ on $L^s$
can be viewed as the square root of the linearizing coordinate.} 

Let $S_N$ denote the family of affine structures on the leaves of
$\an_f$, and let $S_L$ denote the family of orbifold affine structures
on the leaves of $\ac_f$.
Let us also consider the pullback affine structures $q^{-1}S_N$
on the leaves of $\ac_f$. They coincide with $S_L$ on all leaves
except the singular leaf $L^s$.

\begin{lemma}{orb affine str}
The orbifold affine structures $S_L$
 on the leaves of $\ac_f$ make it an affine
orbifold lamination.
\end{lemma}

\begin{pf} We need to check that the affine structure depends continuously 
 on the leaf.  We will use the box basis of $\ac_f$ described above
and the explicit formula for the affine coordinates of 
\S\ref{affine coordinate}.

Take an $x\in \ac_f$ with $q x=\hat z=(z_0,z_{-1},\ldots)\in \NN_f$. 
Let us first assume that $\hat z$ does not lie on
the invariant leaf $ L=L(\hat\alpha)$.
 Then  there is a subsequence $z_{-n(k)}$
staying  distance at least an $\epsilon>0$ from the postcritical set.

Take now a neighborhood $D\ni z_0$ containing at most one point of
the postcritical set. Let $D_{-k}$ denote the pullback of $D$ along
$\hat z$. Let us consider  boxes $B_n=B(D, \hat z, n)$.
Take a big $k$ and let $\hat\zeta\in B_{n(k)}$.
By Lemma \ref{formula} the affine structure on the leaf $\hat D(\hat\zeta)$
 of such a box is given by rescaling $\pi_{-m}=\pi\circ \hat f^{-m}$ and 
passing to limit. But for $m>n(k)$, $ \pi_{-m}=f^{-(m-n(k))}\circ \pi_{-n(k)}$
 for
an appropriate branch of the inverse function. Since 
$\diam (D_{-n})$ is small for $n$ sufficiently large, and $f^{-(m-n(k))}$
allows analytic extension in the $\epsilon $ neighborhood of $z_{-n(k)}$,
by the Koebe Distortion Theorem
 it is almost linear on $D_{-n(k)}$ uniformly in $\hat\zeta$.
It follows that the variation of the affine structure on the leaves of
$B_{n(k)}$ is small, provided $k$ is sufficiently large.

Hence 
 the variation of the pullback structure $q^{-1}S_N$ on
the leaves of the box $\tilde B_{n(k)}\equiv \tilde B(D, \hat z, n(k))$
is also small  for big $k$. Thus the variation of the structure
$S_L$ on all regular leaves of  $\tilde B_{n(k)}$ is small
as well.

Let now $y\in \tilde B_{n(k)}\cap \tilde L^s$, and   $\tilde D^s(y)$ be
the singular local leaf of $y$ in $\tilde B_{n(k)}$.  Let  
 $\phi: (\tilde L^s, \tilde \alpha)\rightarrow (\C, 0)$
 be a regular uniformization of $\tilde L^s$.
Then according to the discussion preceeding this lemma,
 an orbifold affine chart on $\tilde D^s(y)$
is given by $\sqrt{\phi}$. But the image
 $\phi(\tilde D^s(y))$ escapes to $\infty$
in $\C$ when $y\to x$. 
 Moreover, by Lemma \ref{rel dist} 
its size relative to the distance to the origin 
is vanishing.
    Hence the non-linearity of the square root map on
this set goes to zero.  Thus the affine structure $S_L$
 on $\tilde D^s(y)$ is close to the pullback structure $q^{-1}S_N$
on this local leaf. Consequently, it is close to the affine structure
on the leaf $\tilde D(x)$ when $y$ is close enough ot $x$.
  We are done with the case when $\hat z\not\in L$.

Let now $\hat z\in L$, so that $x\in \tilde L^{\mu}$
for $\mu=r$ or $\mu=s$. We wish to check that  the affine structures 
$S_L$ on leaves $\tilde\Delta(y)$ 
of a box  $f^n Q^{\mu} (\Delta, \hat z, N)$
defined by (\ref{intersec})
approach the affine structure on $\tilde\Delta(x)$. 
By pulling back and enlarging the box, we see that it is enough to
check this for $x=\tilde\alpha^{\mu}$ 
and boxes $Q^{\mu}(D)$ defined in (\ref{Q-box}). 

Let first $x= \tilde\alpha^r$. Let us consider a regular orbit 
$$\hat\zeta=(\alpha,\ldots,\alpha, \zeta_{-(N+1)}, \ldots)\in
B^r(D,\hat\alpha,N).$$ Then 
the inverse branches of 
$f^{-(n-N)}: (D_{-N},\alpha)\rightarrow  (D_{-n},\zeta_{-n})$ along
$\hat\zeta$ allow a uniform $\epsilon>0$-enlargement, and hence
have small non-linearity for big $N$. It follows that the affine 
structure $S_N$ on the local leaf $\hat D(\hat\zeta)$ is close to the
regular affine structure on $\hat D^r(\hat\alpha)$ (given by the linearizing
coordinate near $\alpha$). Now we can pass to the orbifold structures
on $Q^r(D)$ in the same way as in the above case $\hat z\not\in L$.

Finally let $x= \tilde\alpha^s$. Let us now consider a singular orbit
$$\hat\zeta=(\alpha,\ldots,\alpha, \zeta_{-(N+1)}, \ldots, \zeta_{-n}=c, 
\ldots)\in B^s(D,\hat\alpha,N),$$ where $n=N+l$. Then  
for sufficiently large $N$ the
rescaled branch $f^{-(n-1)}: D\rightarrow D_{-(n-1)}\ni fc$ is close to the
linearizing coordinate near $\alpha$. 
 The next inverse branch
$f^{-1}: D_{-(n-1)}\rightarrow D_{-n}$ is almost the square root map 
(since $D_{-n}$ is small), while all further
inverse iterates are almost linear on $D_{-n}$.
It follows
that the affine coordinate on the local leaf $\hat D^s(\hat\zeta)$
is close to the square root of the linesrizing coordinate, which is exactly
the orbifold affine coordinate on $\hat D^s(\hat\alpha)$. 

Now we  again can pass from the box $B^s(D)$ to $Q^s(D)$
 in the same way as above. 
\end{pf}  

\subsection{General post-critically finite construction}
\label{general pcf case}
Let  $f$ be an arbitrary post-critically finite map. If a critical
point lands in a cycle then the cycle is either repelling or
super-attracting (contains a periodic critical point) -- see e.g.
\cite[Thms. 1.4,1.6, Prop. 1.11]{lyubich:topological}. In the latter
case this cycle is omitted from $\RR_f$. Thus we need only consider
repelling cycles. 
It also follows that 
there is a uniform bound on the branching index of $\pi$ at all points in
$\RR_f$, since a backward orbit in $\RR_f$ can only hit the critical set a
bounded number of times. 
Given a postcritical repelling periodic point $\alpha$, let 
us consider all occurring branching indices 
$1=d_1(\alpha)< \ldots< d_{l(\alpha)}(\alpha)$ of the leaves over $\alpha$.  

A general construction of the orbifold lamination for a
post-critically finite map has the following
differences as compared
with the previous particular case:

\begin{itemize}
\item
Make $l(\alpha)$
copies of the post-critical 
periodic leaf $L(\hat\alpha)$. 

\item
Supply these copies with orbifold structures of
degrees $d_i(\alpha)$.

\item
Organize the leaves of the lamination over 
$\hat\alpha$ into the boxes according to their
branching indices and then compactify them by adding the corresponding 
orbifold leaves. These boxes will be open in the new topology.
\end{itemize}

\subsection{Structure of the Chebyshev and Latt\`es laminations}
\label{chebyshev and lattes}
Let us consider the quadratic Chebyshev polynomial $p: z\mapsto 2z^2-1$,
$J(p)=[-1,1]$.
Let $T(z)=z^2$, and 
$\phi(z)=1/2(z+1/z)$. Then $\phi\circ T=
p\circ\phi$, so that $p$ is conformally equivalent to $T$ on
the quotient space of $\C^*$ by the involution  $\sigma: z\mapsto 1/z$.

Then the natural extension $\NN_p$ is  the quotient  of the natural extension
$\NN_T$ modulo the involution $\hat \sigma: (z_0, z_{-1},\ldots)\mapsto
(\sigma z_0, \sigma z_{-1}, \ldots)$.  The only invariant leaf of this
involution is the invariant leaf  $L=L(\hat 1)$ of $\hat p$. 
Since $\an_T = \RR_T$ is a regular affine lamination, we obtain a natural
orbifold affine lamination structure on $\an_p$ 
(with one singular leaf $L$).  
The orbifold lamination $\ac_p$ constructed above is obtained from this
one by adding an isolated copy of $L$ (with regular affine structure). 

The situation for the higher degree Chebyshev polynomials is completely
analogous.

 Similarly, the regular leaf associated with the post-critical
fixed point of a Latt\`es map is isolated. Proposition \ref{C and L}
shows that these are the only postcritically finite maps with 
isolated leaves (see Proposition \ref{minimality} for a more general
statement). 
  After removing this leaf, the  lamination becomes the
quotient of the ``torus solenoid'' (that is, the natural extension of 
the torus endomorphism) modulo an involution.

\section{Hyperbolic 3-laminations.}
\label{3d lams}

\subsection{Affine extensions in the abstract}
\label{axiomatic}
In this section, let us forget the specific construction of Section
\ref{pcf orbifolds} and take an ``axiomatic'' approach to what we call
affine extensions. The general construction of  Section \ref{universal
construction} will yield objects of this type.

Let $f:\bar \C \to \bar \C$ be a rational map. An {\em affine
extension of $f$} is an affine (orbifold) 2-lamination $\ac$ with
simply connected leaves,
together with a homeomorphism $\hat f:\ac \to \ac$, acting by 
conformal automorphisms on leaves, and a projection
$\pi:\ac\to \bar\C$, such that
\begin{enumerate}
\item $f \circ \pi = \pi \circ \hat f$
\item $\pi$ is continuous, and restricted to any leaf is non-constant
and complex-analytic.
\end{enumerate}
Condition (1) immediately implies that $\pi$ factors through a map
$p:\ac\to\NN_f$, given by 
$$
	p(\bold{z}) = (\pi(\bold{z}),\pi \hat f^{-1}(\bold{z}),\pi
					\hat f^{-2}(\bold{z}),\ldots). 
$$
Let $\pi_{k} = \pi\circ\hat f^k$, as usual.

In fact $p$ is continous by (2), and 
we immediately see that 
$p(\ac)$ is contained in $\RR_f$: on any leaf $L$, $\pi$ factors
through $f^n$ for any $n>0$ and so the pullbacks $\pi\circ \hat
f^{-n}(U)$ for any disk $U\subset L$ with compact closure are
eventually unbranched. Thus $p$ restricted to each leaf is a 
complex analytic map to a leaf of $\RR_f$, and we further conclude
that the leaf must be parabolic. Hence $p:\ac \to \an_f$.

The construction in \S\ref{pcf orbifolds} yields just such
an object, and as in that case
the map $p$ need not be injective: it re-identifies the leaves which
we separated in our construction.

\subsection{Extending to three dimensions}
Even before we consider the action on $\ac$ we can 
associate to it a naturally defined 
${\Hyp}^3$-lamination $\hc$,
by attaching a copy of hyperbolic 3-space, realized as its upper
half-space model, to (a finite cover of) every leaf.
Since transition maps for affine charts on the leaves of $\ac$ are
affine, they extend naturally to isometries on the hyperbolic
3-spaces.

In particular given two affine charts $\phi,\phi':\C\to L$, the
corresponding transition map from $\Hyp^3$ to $\Hyp^3$ multiplies
heights by the norm of the derivative of $\phi^{-1}\circ\phi'$. Thus
we can consider a copy of $\Hyp^3$ for each chart $\phi$, and
define the leaf $H_L$ attached to $L$ as the identification of
all these copies via the transition maps. However, we prefer to make
the following definition, which will be easier to
work with:

Consider the group
$\Aff$ of complex-affine maps $A:\C\to\C$ (henceforth just ``affine'').
We can identify the
complex plane $\C$ and the hyperbolic space 
$\Hyp^3\equiv \C\times (0,\infty)$  (with a preferred point at
$\infty$)  as  
homogeneous spaces for  $\Aff$, namely  $\C\homeo \Aff/\C_*$ and
$\Hyp^3\homeo \Aff/S^1$. 
In other words, consider the projections
$p_1:\Aff \to \C$ and $p_2:\Aff\to\Hyp^3$ given by
\begin{equation}\label{p1 def}
p_1:g \mapsto g(0) \in \C
\end{equation}
and
\begin{equation}\label{p2 def}
p_2:g \mapsto g(0,1)=(g(0),|g'|) \in \Hyp^3.
\end{equation}
Fibres of $p_1$ are orbits of the {\em right action} of the subgroup $\C_* =
Fix(0)\subset \Aff$, that is $\{z\mapsto \alpha z: \alpha\ne 0\}.$
Fibres of $p_2$ are orbits
of the right action of $S^1$, that is the group
$\{z\mapsto \lambda z: |\lambda|=1\}.$
The left-action of $\Aff$ on itself projects to complex-affine maps on
$\C$, and to hyperbolic isometries on $\Hyp^3$.

Now suppose first that $\ac$ has no orbifold leaves, and for a leaf $L$
consider the set
$\{\phi:\C\to L\}$ of affine isomorphisms (``charts'') from $\C$ to
$L$, which admits a fixed-point-free 
right-action by $\Aff$. We may identify $L$ with $\{\phi\}/\C_*$
by taking $\phi$ to $\phi(0)$. The space $\{\phi\}/S^1$ is  naturally
identified with $\Hyp^3$ as above, and we call this the hyperbolic 
leaf $H_L$ associated to $L$. 

Thus, the total space $\hc$ may be defined as $\{\phi:\C\to
\ac\}/S^1$, where the maps $\phi$ vary over all charts for leaves of $\ac$.
This clearly inherits the structure of a hyperbolic 3-lamination.
We will usually write $[\phi]$ for an equivalence class of charts
modulo rotation in $S^1$.

The same construction works for the orbifold leaves, with the
charts replaced by finite coverings.  Thus
an orbifold affine 2-lamination extends to 
a hyperbolic 3-orbifold lamination.

One should think of a chart $\phi:\C\to L$ as determining a point and
a choice of {\em scale} for the leaf $L$. Changes of scale correspond
to vertical motion in the upper half-space model. Indeed, let $e^\R$
denote the subgroup of $\C_*$ acting by scaling without rotation. 
The $\R$-action induced on $\ac$ by the right-multiplication 
$r:[\phi] \mapsto [\phi\circ e^r]$ is simply the vertical geodesic flow in
each leaf, where $r$ measures arclength and increasing $r$ corresponds
to increasing heights in each leaf (as is evident from (\ref{p2 def}).

Finally, we remark that this extension of an orbifold affine
lamination to an orbifold $\Hyp^3$-lamination is unique, in the sense
that if $\hc'$ is another orbifold $\Hyp^3$-lamination with a
projection $\hc'\to\ac$ such that on each leaf, fibres of points are
geodesics with a common endpoint at infinity, then $\hc$ and $\hc'$ are
related by an isomorphism fixing $\ac$.

\subsection{Proper discontinuity of actions}
The action $\hat f$ on $\ac$ (even without assuming that it projects
to a rational map) extends naturally to an action, 
which we also call $\hat f$, on $\hc$ by 
hyperbolic isometries, namely
\begin{equation}
\label{action on H}
\hat f : [\phi] \mapsto [\hat f\circ \phi].
\end{equation}

It is useful to note that we now have two commuting actions on $\ac$: 
a $\Z$-action generated by $\hat f$ on the left, and an $\R$-action,
the vertical geodesic flow, generated by $e^\R$ on the right. 
These actions have a certain {\em coherence}:  forward iterates
of $\hat f$ tend to increase heights, as a result of the general
expansive properties of the rational map $f$. Let us make this precise
with the following statement:

\begin{lemma}{coherence}
Let $(\ac,\hat f,\pi)$ be an affine
extension of a rational map $f$, and let
$\hc$ be the hyperbolic 3-lamination associated to $\ac$.
For any two points $p,q\in \hc$ there are neighborhoods $U_p,U_q$ for
which the following holds: if $n_i, r_i$ are  sequences such that
$$
(\hat f^{n_i}\circ U_p) \intersect (U_q \circ e^{r_i}) \ne \emptyset
$$
then $n_i \to +\infty$ if and only if  $r_i \to +\infty$, and 
$n_i \to -\infty$ if and only if  $r_i \to -\infty$.
\end{lemma}

In other words, whenever high forward/backward iterate of $\bz\in U_p$ is comparable
with $\bzeta\in U_q$ in the sense that these points lie on
 the same vertical  geodesic, the former point is much higher/lower than the latter.

\begin{pf}
Represent $p$ by a chart $\phi:\C\to \ac$ and $q$ by a chart
$\psi:\C\to \ac$. Recall that $\pi\circ\phi$ is analytic, and hence
its image misses at most two points in $\bar\C$. Thus there exists an
open set $W$ which meets the Julia set $J_f$, and a disk $D\subset \C$
around 0 such that $W \subset \pi\circ\phi(D)$.
Let $U_p$ be small enough that for any $[\phi']\in
U_p$, $\pi(\phi'(D))$ contains $W$. This is possible since
$\pi$ is continuous in $\ac$. In addition choose $U_p$ small enough
that there is some upper bound on the degree of $\pi\circ \phi'$ 
in $D$ (here by
``degree'' we mean the maximal degree over any point in the image).

Now we can see that $\pi\circ\hat f^n\circ\phi'(D)$ will tend to blow up as
$n\to\infty$, and down (in diameter) as $n\to -\infty$. Indeed, 
there is some $n_0$ such that $f^n(W)$ contains all of $J_f$ for
$n>n_0$, and as $n\to\infty$ the degree of $f^n$ on $W$
increases without bound -- hence the same is true for $\pi\circ\hat
f^n\circ\phi'$ on $D$, for any $[\phi']\in U_p$. 

To see what happens to $\pi\circ\hat f^{-n}\circ\phi'(D) = \pi_{-n}(\phi'(D))$ 
as $n\to \infty$, we may invoke the Shrinking
Lemma given in Appendix 2, once we
observe two things: (1) the degree of $f^n$ on $\pi_{-n}\phi'(D)$ is 
bounded by the degree of $\pi$ on $\phi'(D)$, and hence uniformly over $U_p$.
(2) Since every leaf of $\ac$ is affine,
$\pi_{-n}\phi'(D)$  is eventually outside the closure of the  rotation domains,
so that
 $\diam (\pi_{-n}\phi'(D)) \to 0$ as $n\to \infty$, uniformly for
all $[\phi']\in U_p$.

Now let us choose  $U_q$ such that for $[\psi']\in U_q$ the
degree of $\pi\circ\psi'$ on $D$ is uniformly bounded.
Suppose that we have $\hat f^{n_i} \circ \phi'_i
\equiv \psi'_i \circ e^{r_i}$ in $\ac$ for $[\phi'_i]\in U_p$ and
$[\psi'_i]\in U_q$. Then if $n_i\to\infty$ then  
$r_i\to\infty$ as well, so that the degree of $\pi\circ\psi'_i$ on
$e^{r_i}(D)$ can go to infinity. 
Conversely, suppose that $n_i \to -\infty$. Then the above Shrinking
Lemma argument implies that $\diam \pi \circ \psi'_i(e^{r_i}D)\to 0$,
and so $r_i$ must go to $-\infty$. 
The other two implications are similar.
\end{pf}

Recall that a group action on a space $X$ is {\em proper} if, for any
two points $p,q\in X$,  there exist
neighborhoods $U\ni p$ and $V\ni q$ so that the subset of group
elements $g$ such that $gU \intersect V \ne \emptyset$
has compact closure in the group. If the group has the discrete
topology this set must be finite, and we say the action is {\em
properly discontinuous}.

A consequence of the Shrinking Lemma, as used in Lemma \ref{coherence},
is the following central fact: 

\begin{proposition}{action discontinuous}
Let $(\ac,\hat f,\pi)$ be an affine
extension of a rational map $f$, and let
$\hc$ be the hyperbolic 3-lamination associated to $\ac$.
The induced action of $\hat f$ on $\hc$ is properly discontinuous.
Similarly, the vertical geodesic flow on $\hc$ is a proper $\R$-action.
\end{proposition}

\begin{pf}
Given $p$ and $q$ in $\hc$, choose
$U_p$ and $U_q$ as in lemma
\ref{coherence}. Then in particular 
$(\hat f^{n}\circ U_p) \intersect U_q = \emptyset$ for all but
finitely many $n$. 
Geometrically, we say that forward iterates of $U_p$ cannot continue to
intersect $U_q$, since their heights are going to infinity whenever
comparison is possible.

Similarly, $U_p \intersect (U_q \circ e^r) = \emptyset$ for all but a
bounded set of $r$.
\end{pf}

At this moment we can conclude
 that  the quotient $\hc/\hat f$ is a Hausdorff space which inherits the 
structure of
the  hyperbolic orbifold 3-lamination. On the other hand, the  quotient via the
flow action recovers the affine lamination $\ac$. This 
 duality between $\hc/\hat f$ and $\ac$  will be useful
in what follows (see Proposition \ref{convex cocompactness2}).

\subsection{Convex Hulls}
\label{convex hull of lamination}
In analogy with the situation in Kleinian groups, we denote by
$\CC(\jc)$ the {\it convex hull} in $\hc$ of the lift $\jc=\pi^{-1}(J)$ of the
Julia set to $\ac$.  This is simply the union of the convex hulls of 
$(\jc\intersect L)\union\{\infty\}$ in $H_L$ for every leaf $H_L$ of
$\hc$ bounded 
by a leaf $L$ of $\ac$.   The quotient $\CC(\jc)/\hat
f$ can be called the {\em convex core} of $\hc_f/\hat f$.

Using lemma \ref{variation of convex hulls} we can obtain the
following. Let $\CC_\delta$ denote the leafwise $\delta$-neighborhood
of $\CC(\jc)$.
\begin{corollary}{boundary is lamination}
For $\delta>0$, $\CC_\delta$ inherits the structure of a 3-lamination
with boundary. In fact $\CC_\delta$ is homeomorphic to
$\hc\union\fc/\hat f$, where 
$\boundary\CC_\delta$ is taken 
to the ``boundary at infinity'' $\fc/\hat f$.

Except when the Julia set of $f$ is smooth, the above holds for
$\delta=0$, and we note also that $\boundary\CC$ inherits a metric
 from $\hc_f$ which makes it into a hyperbolic 2-lamination. 
\end{corollary}

\begin{pf}
Since $\jc_f$ is the pullback of $J_f$ by $\pi$ and $\pi$ varies
continuously in the transverse direction, for any product box $T\times
D$, if $T$ is sufficiently small then the intersections of $\jc_f$
with the local leaves $\{t\}\times D$ are close to each other in the
Hausdorff topology in $D$. The same applies to the finite covers of
orbifold boxes, and hence for any (large) closed disk on a leaf we can
take a small transversal neighborhood so that the Julia sets vary only
slightly in the Hausdorff topology.

It follows, applying lemma \ref{variation of convex hulls},
that any point $x\in \CC_\delta$ has an (orbifold) box neighborhood in
$\hc$ which intersects $\CC_\delta$ in a set of the form $T\times C$, up
to bilipschitz homeomorphism.
Here $C$ 
is the intersection of $\CC_\delta$ with a leafwise neighborhood of $x$.

The homeomorphism from $\CC_\delta$ to $\hc\union \fc/\hat f$ comes
directly from the leafwise homeomorphism discussed in section
\ref{geometry for dynamicists}. The case where $J_f$ is smooth
corresponds exactly to the case where $\jc\intersect L$ is contained
in a straight line in $L$ (again, the Shrinking Lemma), and this is
the only case where the discussion fails for $\delta = 0$.
\end{pf}

\subsection*{The $z^2+\ep$ case.}
The simplest possible quadratic polynomial is $f(z) = z^2 + \ep$ where
$\ep$ is small (more precisely, let $\ep$ lie in the main cardioid of the
Mandelbrot set, so that $f$ has one attracting fixed point).
In this case $J$ is a quasi-circle, the Fatou domain lifts to two
components 
in $\RR_f$, each of which has quotient homeomorphic to Sullivan's
solenoidal Riemann surface lamination $\SS$ (see Appendix 1),
and in
fact $\RR_f$ is already an affine lamination (Proposition \ref{nonrec3}). 

In each leaf, $\jc$ is a quasi-line separating the plane into two
components, where
one component projects to the outside and one to the inside of $J$ in
$\C$. We claim that the convex core $\CC$ is (for $\ep\ne 0$) simply a
product $\SS\times[0,1]$.
This makes concrete the analogy between $z^2+\ep$ and a
quasi-Fuchsian group.

To prove this, or the equivalent fact that $\hc\union\fc/\hat f \homeo
\SS\times [0,1]$, make the following leafwise construction. Let $H$ be
a leaf of $\hc$ bounded by $L$. 
Foliate each component $D$ of $\fc\intersect L$ by Poincar\'e
geodesics coming from infinity in $L$ (vertical geodesics in the upper
half-plane uniformization). Above each such ray $r$ lies a ``curtain''
in $H$, bounded by the vertical line above the point $r\intersect
\jc$. Let $l$ be the union of two rays on opposite sides meeting at
$\jc$. The curtain above $l$ is, in the induced metric, isometric to
$\Hyp^2$, and the vertical line $v$ above $l\intersect \jc$ is a geodesic.
Use the orthogonal projection to $v$ in this surface to define a
product structure. This varies continuously with the lines $l$ in $L$,
and varies continuously in the transverse direction of the lamination.
Thus it gives  a product structure for the entire lamination, which is
also preserved by $\hat f$, since it is clearly affinely invariant.

Note in particular that the convex core is compact. In section
\ref{convex cocompactness section} we will discuss this phenomenon
more generally. 

\subsection{The scenery flow}
In Bedford-Fisher-Urbanski \cite{bedford-fisher-urbanski} a
construction called the ``scenery flow'' is discussed, which is
related to the constructions of this paper, in the case of
an axiom A rational map $f$. 

The scenery flow is, roughly, the set of all ``pictures'' of the
Julia set at small scales. That is, one considers all (complex) affine
rescalings (and rotations) of $J$ in $\C$ and takes limits in the
Hausdorff topology. The resulting collection of subsets of $\C$ is
indexed by backward orbits of $f$, using the linearization formula
(\ref{chart}): for each backward orbit $\hat z$ we consider
the Hausdorff limit $J(\hat z)$ of the sets $J_n = A_n(J)$ where 
$A(z) = (f^n)'(z_{-n})(z-z_{-n})$. The natural action of $f$ on such a
set is $\hat f(J(\hat z)) = J(\hat f(\hat z)) = f'(z_0)\cdot J(\hat
z)$. The flow on the set of pictures is defined by $J(\hat z)\mapsto
e^t J(\hat z)$. 

Translating this into our terminology, given $\hat z\in \jc$ a point
in $\hc$ lying above $\hat z$ can be written as
$[\phi]$ where $\phi:\C\to L(\hat z)$ is a chart such that $\phi(0) =
\hat z$. 
The corresponding picture $J(\hat z)$ is given by 
$\phi^{-1}(\jc\intersect L(\hat z))$, and
the scaling flow is exactly the vertical geodesic flow (this
interpretation of rescaling as geodesic flow was part of the original
motivation for \cite{bedford-fisher-urbanski}) .
Thus the scenery  flow is taken to the ``curtain'' above the lift of
the Julia set.

\section{Universal orbifold laminations}
\label{universal construction}

In this section we introduce the machinery for a general construction
of an affine orbifold 2-lamination and accompanying hyperbolic
orbifold 3-lamination, for any rational map. 

In the original construction we were faced with the following issue: 
A small disk $D$ in the affine part could be approached (in $\NN_f$) by a
sequence of disks $D_i$ in such a way that the projections $\pi$ on
$D_i$ do not converge in any sensible way to the projection on $D$.
For example there could be branching on the $D_i$ whereas $D$ projects
univalently. We resolved this in the post-critically finite case 
by creating new leaves and redefining
the topology so as to sort out the different branching possibilities. 

 From the point of view of this section, the projection maps
themselves from the affine leaves to the sphere will be the basic
objects, so that the topology will automatically include convergence
of the maps.  Thus we will consider a space of meromorphic functions,
with an associated action by affine transformations of the domain
which gives rise to a leaf structure (that is, we can think of a leaf
as a set of choices of basepoint for a meromorphic function, with
precomposition with affine maps giving the change of basepoint). This
will be our ``universal'' orbifold foliation, and any rational map
will act naturally on it and give rise to an invariant lamination in
which our original affine space $\an_f$ will be a subset.  The new
topology $\al_f$ is induced from this space, and a final closure step
will yield the added leaves.

\subsection{ Leaves of the affine action in the Universal space }

Let $\tilde\uu$ denote the space of  meromorphic functions on
$\C$, with the topology of uniform convergence on compact sets,
and let $\uu$ denote the  open subset of non-constant functions.
Since $\tilde\uu$ is a complex vector space, $\uu$ can be viewed as
an infinite-dimensional complex analytic manifold (anlyticity amounts to
analytic dependence on Taylor coeffiecients).

The space $\uu$ admits two natural commuting
analytic  actions:
a {\em left-action}
$\psi\mapsto f\circ\psi$ by the semigroup of rational maps
$f:\bar\C\to\bar\C$, and a {\em right-action} 
$\psi\mapsto \psi\circ A$ by the group $\Aff$ of 
complex-affine maps $A:\C\to\C$.

Let us first consider the structure of the individual orbits $\psi\circ\Aff$ of
the right-action of $\Aff$ on $\uu$, 
and later show that they fit together into a foliation.
On each orbit we place the
{\em leafwise topology}, in which open neighborhoods are sets of the
form $u\circ V$ where $u\in\uu$ and $V$ is an open set
in $\Aff$. Note that this may be a stronger topology
(more open sets) than the induced topology from $\uu$, since a leaf
may accumulate on itself in $\uu$. 

The map $\Aff \to \psi\circ\Aff$ is locally non-singular -- that is, the
derivative map $D_\psi: T_{id}(\Aff) \to T_\psi(\uu)$ is
non-singular, as one may check by explicit computation.
Note that the tangent space $T_\psi(\uu)$ can be
identified with the space $\tilde\uu$. 
It follows that, for $h$ sufficiently close to but
not equal to the identity, $\psi\circ h \ne \psi$. Thus
the  isotropy subgroup $\Gamma_\psi
= \{\delta\in\Aff: \psi\circ\delta = \psi\}$ is discrete in $\Aff$.
We may therefore make the identification
$$
	\psi\circ\Aff \homeo \Gamma_\psi \backslash \Aff.
$$
which is a homeomorphism if $\psi\circ\Aff$ is taken with the leafwise
topology. (The
quotient is on the left since $\Aff$ acts on the right, so that 
$\psi\circ g = \psi \circ h \iff g=\delta\circ h$ for $\delta\in\Gamma_\psi$).

Note also that $\Gamma_\psi$ must in fact consist of isometries of $\C$ since
a non-constant meromorphic function cannot be invariant under a dilation.

Now as in Section \ref{3d lams}, we may think of
$\C$ and $\Hyp^3$ as the quotients
$\C\homeo \Aff/\C_*$ and $\Hyp^3\homeo \Aff/S^1$, with associated
left-action of $\Aff$.

Since right and left actions commute, we may form the quotients
$$
L_{aff}(\psi)\equiv \psi\circ\Aff/\C_* \homeo \Gamma_\psi\backslash\Aff/\C_*
\homeo  \Gamma_\psi\backslash\C,
$$
which is a Euclidean 2-orbifold, and
$$
L_{hyp}(\psi)\equiv \psi\circ\Aff/S^1 \homeo \Gamma_\psi\backslash\Aff/S^1 
\homeo \Gamma_\psi\backslash\Hyp^3,
$$
which is a hyperbolic 3-orbifold. Note that the singularities, always
arising from rotations in $\Gamma_\psi$, are  cone axes.

The natural projection $L_{hyp}(\psi)\ra L_{aff}(\psi)$ is a one-dimensional fiber
 boundle whose leaves are the orbits of the {\it vertical flow}
$V_r: \phi\ra \phi\circ e^{r}$ (this flow is well-defined since $\C_*$ is commutative). 

As an example, consider $\psi (z) = z^m$, so that $\Gamma_\psi$ is the
cyclic group generated by $\xi:z\mapsto e^{2\pi i/m} z$.
The leaf $\psi\circ\Aff/\C_*$ is then the orbifold
$\langle\xi\rangle\backslash\C$ with one order $m$ cone point. More
interesting examples are Chebyshev polynomials associated with trigonometric functions
and  Latt\`es maps associated with elliptic functions.

A local (orbifold) affine
chart on a leaf $L_{aff}(\phi)$ near $\phi$ is given by translations $t\mapsto \phi(z+t)$, where
$t\in \C$ is small. In these coordinates
the map $f: L_{aff}(\phi)\ra L_{aff}(f\circ\phi)$ becomes the
identity. Thus $f$ is affine on the leaves.  
Hence it is automatically a covering.

Similar statements are valid for the hyperbolic leaves of $L_{hyp}(\phi)$, with local charts
$(t,e^r)\mapsto \psi(e^r z+t)$.

\subsection{Foliation structure}
\label{foliation structure}
With this point of view on individual leaves, let us consider how they
fit together into the total space $\uu$, and its quotients.

\begin{lemma}{really foliation}
The $\Aff$ action supplies  the 
space $\uu$ with  an analytic foliation 
with two complex dimensional leaves. 
\end{lemma}
\begin{pf}
This is  a generality about any non-singular analytic Lie group action.
However, rather than using deep Implicit Function Theorems
(see \cite{H}), we can check the statement directly.

Let $\phi\in \uu$. Without loss of generality we can assume that 
$\phi'(0)\not=0$. Let 
$$T=\{\theta\in \uu: \theta(0)=\phi(0); \theta'(0)=\phi'(0)\}.$$
We will show that $T$ is a local transversal to the action of $\Aff$.
Indeed take a $\psi\in \uu$ near $\phi$ and a $\gamma\in \Aff$ near $\id$,
$\gamma(z)=az+b$. The condition that $\psi\circ\gamma\in T$ amounts
to the following system of two equations for $a$ and $b$:
\begin{equation}\label{Hurwitz}
\psi(b)=\phi(0),\;\; a\psi'(b)=\phi'(0).
\end{equation}
If $\psi$ is close to $\phi$ then the first equation has a unique root $b$
near 0 by the Hurwitz theorem.  Thus the second equation
has a unique root $a$ near 1.

It follows that the $\Aff$ action has a local product structure near $\phi$
given by  the map
$T\times \Aff\rightarrow \uu$,
$(\theta,\gamma)\mapsto \theta\circ\gamma$ near $(\phi, \id)$.
This structure is analytic since this map is so. 
The inverse map is also analytic as the solutions of (\ref{Hurwitz})
analytically depend on the Taylor  coefficients of $\psi$ 
(by the Implicit Function Theorem).
\end{pf}

Now we may form the quotients $\ua=\uu/\C_*$ and $\uh=\uu/S^1$, and we claim
that they are orbifold 2- and 3-foliations, respectively.

This follows from the following general fact. 
Suppose $\LL$ is a lamination with
finite-dimensional smooth leaves and 
a Lie group $G$ acts on $\LL$ (say from the right) preserving leaves.
We call the action smooth if its leafwise derivative exists, and is continuous 
in $\LL$ (in the transverse direction as well).

\begin{lemma}{general orbifold}
Let $\LL$ be a lamination with
finite-dimensional leaves, admitting a nonsingular proper smooth action
by a Lie group $G$.
Then $\LL/G$ is an orbifold lamination, where the leaves have
dimension equal to the codimension of $G$-orbits in the leaves of
$\LL$. If $\LL$ is actually an analytic foliation and the action of $G$
is anlytic then $\LL/G$ is an analytic orbifold  foliation as well. 
\end{lemma}

\begin{pf} Note first that the properness of the action ensures that 
the quotient $\LL/G$ is Hausdorff.
 
Let now $p\in B\subset U$ where $B$ is  a product box $B=T\times V$, with
$V$ a leafwise neighborhood. Write $p = (t,v) \in B$. Then because the
$G$-action  is smooth and non-singular, we can find a continuous
family of transversals 
$K_s$ to the $G$-orbits in each local leaf $\{s\}\times V$. The union
$K$ is a transversal to the $G$ action, which itself has a product box
structure. 

The subgroup $G_p$ fixing $p$ is discrete by the non-singularity
assumption, and is finite by the properness assumption. We now get a
``first return'' action of $G_p$ on a small enough neighborhood of $p$
in the transversal $K$, and the quotient of this neighborhood by this
action is our orbifold box in the quotient $\LL/G$.
To see this, note that if $K'$ is a
sufficiently small neighborhood of $p$ in $K$ then for any $q\in K'$
and $g\in G_p$, $qg$ is in the original neighborhood $B$, and hence
can be uniquely pushed to $K$ along its $G$-orbit in $B$. Thus each
element of $G_p$ induces a map $K'\to K$ fixing $p$ and altogether we
obtain a finite group action on the union of images of $K'$.

Finally, it is obvious that if the lamination $\LL$ and the action of $G$ have
some transversal regularity (e.g., analytic), then the quotient
lamination inherits  it.
\end{pf}

Let us summarize the above discussion:

\begin{corollary}{universal foliations}
The quotient $\uh=\uu/S^1$ is a hyperbolic orbifold  3-foliation
The quotient $\ua=\uu/\C_*$ 
is an affine orbifold 2-foliation, and the
projection $\uh\ra \ua$ is a fiber bundle.
On the leaves of $\uh$ it is identified with the
vertical projection in each half-space to the bounding plane. 
\end{corollary}

{\bf Remark:} The projection $\uu\to \uh$ is similar to a Seifert
fibration. A function admitting rotational symmetries around $0$ gives
rise to a singular fiber: its $S^1$ orbit is finitely covered by the
$S^1$ action, whereas for nearby functions without the symmetry the
orbit is an injective image of $S^1$. However, note that singular
fibers are not isolated, as they are in Seifert-fibred three-manifolds.

\begin{pf}  
To apply Lemma \ref{general orbifold}
 we need to check that the actions
of $\C_*$ and $S^1$ are proper. For $S^1$ this is clear since it is
compact. For $\C_*$  we just have to consider the vertical flow
 $\psi \mapsto \psi\circ e^r, r\in\R$ (as in Section \ref{3d lams}). But
it is easy to see 
that as $r$ goes to $\pm\infty$ $\psi\circ e^r$ diverges in $\uu$ --
it becomes a constant in one direction, and blows up at every point in
the other. In fact for a small enough neighborhoods $U\ni\psi$ and
$V\ni\phi$ there is
a fixed $R$ so that for $|r|>R$ the rescaling $u\circ e^r$ is outside $V$
for any $u\in U$. This proves that $\C_*$ acts properly.

Note finally the local affine charts are transversally analytic, so
that we obtain 
an affine foliation. Indeed taking an analytic transversal $K$ to the 
foliation $\ua$, the map $(\psi, t)\mapsto \psi(z+t)$, where $\psi\in K$,
 $t\in \C$ is small, provides us with an orbifold affine box.
Similarly, the hyperbolic structure on $\uh$ is transversally analytic.
\end{pf}

We remark that it is easy to see that $\uu$ is metrizable, and in fact
one can give it a complete metric which is invariant under the right
$\C_*$-action. However, we shall not need this explicitly.

\subsection{Characteristic laminations}
Now given a rational map $f:\bar\C\to\bar\C$, we extract from our
universal space $\uu$ the characteristic orbifold laminations for
$f$. First consider the  ``global attractor"
$$
\k_f\equiv \k  = \bigcap_{n\ge 0} f^n(\uu)
$$
which is the maximal invariant subset for which $f:\k \to \k$ is
surjective. Note also that $\k$ is naturally a
sublamination since it is leafwise saturated. 

Let us show that $\k$ is closed in $\uu$. It is enough to check that
for any rational map $g$, $g(\uu)$ is closed. Let $g\circ \phi_n\to \psi$.
Then $\{\phi_n\}$ is a normal family. Indeed, given any point $a\in \C$,
consider two neighborhoods $U\Supset  V\ni a$. 
Then eventually for all $n$, $\phi_n(V)\subset g^{-1}\circ\psi(U)$.
Take $U$ so small
that the complement of $g^{-1}\circ\psi(U)$ has non-empty interior.
By Montel's Theorem,   the family  $\{\phi_n\}$ is normal on $V$. As normality
is a local property, $\{\phi_n\}$ is normal.
Let $\phi: \C\rightarrow \bar\C$ be any limit function. Then $\psi=g\circ \phi$,
and we are done.

It is not necessarily true that $f|_\k$ is 
injective (see remark below).
Thus we take the natural extension, or inverse limit, of the system
$ \k \underset{f}{\leftarrow} \k \underset{f}{\leftarrow} \cdots$.
Call this new system $\hat f:\hat\k_f \to \hat\k_f$. Elements in $\hat
\k\equiv \hat\k_f$ are simply sequences $\hat \psi=\{\psi_n\in\uu\}_{n\le 0}$ such that
$\psi_{n+1} = f\circ\psi_n$.

Note that $\hat\k$ is still naturally a leaf space: 
$\k$ is invariant under the right action of $\Aff$, 
which then extends to $\hat\k$ via $\{\psi_n\}\circ A
= \{\psi_n\circ A\}$, so that the leaves project down to (orbifold)
cover leaves
in $\k$. We want to check that $\hat\k$ is in fact a lamination,
i.e. that there is a local product structure.

\begin{lemma}{K hat lamination}
$\hat\k_f$ is a lamination  whose leaves are the right $\Aff$-orbits.
The projection from $\hat\k_f$ to $\k_f$ is
an orbifold covering on leaves.
Similarly, $\hat\ka_f\equiv \hat\k_f/\C_*$ 
and $\hat\kh_f\equiv \hat\k_f/S^1$ are orbifold affine 2- and
hyperbolic 3-laminations, respectively. 
\end{lemma}

\begin{pf}
Fix $\hat\psi = \{\psi_n\}$ in $\hat\k$. We will describe the
structure of a neighborhood of $\hat\psi$ as follows.
Let $D$ be some disk in $\C$ on which $\psi_0$ is univalent, and let
$D'' \Subset D' \Subset D$ be nested open disks.
Let $U$ be an open neighborhood of
$\psi_0$ in $\k$ for which any $u\in U$ is univalent in $D'$ and 
such that $\psi_0(D'')\subset u(D')$. 
Note that we may assume $U$ is a product
neighborhood of the form $T\times V$ where $V$ is a neighborhood of
the identity in $\Aff$. 

The preimage $\hat U$ of $U$ in $\hat\k$ consists of sequences
$\hat u = \{u_n\}$ for which $u_0\in U$. For any such $\hat u$, notice
that, since $u_0 = f^n\circ u_{-n}$ for any $n\ge 0$, 
$u_{-n}$ is univalent on $D'$ 
and $f^n$ is univalent on
$u_{-n}(D')$. 
Thus $\hat u$ determines an infinite sequence of univalent pullbacks
of $u_0(D')$, and hence of $W_0 \equiv \psi_0(D'')$. Conversely $\hat u$ is
determined by this sequence and $u_0$, using analytic continuation.

Let $\Sigma$ denote the set of all possible infinite pullback sequences
$W_0, W_{-1},\ldots$ where
$f$ is univalent at each step.
With the natural topology, 
this is a closed subset of the set of all possible pullback sequences
for $W_0$, and hence a closed subset of a Cantor set.
We thus have an injection of $\hat U$ into $\Sigma\times U$.
Let us show it is also a topological
embedding. Suppose for a sequence $\hat u^i$ that the images 
$((W^i_{-n}),u^i_0)$ in
$\Sigma\times U$ converge. Then $u^i_0$ converge as meromorphic
functions, and for each $n>0$ eventually $W^i_{-n}$ are constant.
The functional equation
$u^i_0 = f^n \circ u^i_{-n}$ then implies that 
$u^i_{-n}$ converge locally on $D''$ as $i\to\infty$, and
in fact form a normal
family so that they converge
globally. Thus $\hat u^i$ converge. 
The other direction is easy.

The subset of $\Sigma\times U$ obtained is saturated in the leaf
direction, since any 
composition $u_0\circ A$ lying in $U$ (with $A\in \Aff$)
can be pulled back along the
same sequences as $u_0$. Hence there is some subset $Q\subset \Sigma\times T$
so that we may identify $\hat U$ with $Q\times V$. This is the desired
product box.

The fact that $\hat\k/\C_*$ and $\hat\k/S^1$ are orbifold
laminations now follows by another application of Lemma 7.2.
\end{pf}

\noindent{\it Remark.}
We expect that in most cases the natural extension step is
unnecessary; that is, $f$ is already injective on $\k$.
Counterexamples are maps with symmetry: for example, 
the leaf of $\psi_0(z) = e^z$ in $\uu/\C_*$ is a cylinder $\C/{2\pi i}$, 
and if $f(z) = z^d$ then $f$ is a $d$-fold cover from this leaf to
itself. It follows that
the lift of this leaf to $\hat\k_f/\C_*$ is
a solenoidal Riemann surface lamination (in fact it is just the
original $\an_f$ in this case).
We conjecture that non-injectivity only happens when $f$ or a
power of $f$ has a M\"obius symmetry. 

\subsection{Completion}
There is an  equivariant inclusion of $\an_f$ into our new object 
$\hat\ka$,
as follows. Let $\hat z$ be a point on a leaf $L$ of $\an_f$, and let
$\phi:\C\to L$ be an isomorphism such that $\phi(0) = \hat z$.
Then for $n\le 0$,  $\pi_n\circ\phi$ is an element of $\k$ (compare
\S\ref{3d lams}). The
choice of $\phi$ was determined only up to precomposition by $\C_*$,
so that $\hat z$ determines a well-defined sequence in 
$\ka$, which gives an element $\iota(\hat z) \in \hat\ka$.

The map $\iota$ takes leaves to leaves, since another element of $L$
can be written as $\phi(A(0))$ with $A\in\Aff$. 
$\iota$ is
injective, since at least one of the coordinates $z_n$ must differ
for different points on $\an_f$.

On the other hand, $\hat\ka$ is also an affine orbifold extension
of $f$, in the sense of Section \ref{axiomatic}, and hence there is 
also a continuous, equivariant projection 
$p:\hat\ka \to \an_f$. That is, 
for any $[\hat\psi]=\{[\psi_n]\}\in\hat\ka$, let $p([\hat\psi])$ be
the backward orbit $\{\psi_n(0)\}$.

It is immediate from the definitions that $p\circ \iota$ is the
identity, but we note that the opposite is false, since in fact
$p$ is not injective and hence $\iota$ is not surjective. Indeed, let
$g:\C\to\C$ be any non-affine entire function and $L$ a leaf of $\an_f$
with chart $\phi:\C\to L$. Then
the sequence $\{\pi_n\circ\phi\circ g\}$ is on a leaf of
$\hat\k/\C_*$ which projects to $L$ but is different from $\iota(L)$.

Note that the topology on $\an_f$ induced from
$\hat\ka$ is in general stronger than its own topology, induced from 
$\NN_f$ (so that the inclusion $\iota$ is discontinuous). This is in
fact the main point of the construction.
Let $\al_f$ denote $\iota(\an_f)$, with the topology induced from $\hat\ka$.
We also think of $\an_f$ and $\al_f$
as being the same underlying space, with different topologies, which
we call ``natural'' and ``laminar''.

Our final step is to take
the closure, in $\hat\ka$, of $\al_f$, obtaining
automatically an affine orbifold extension of $f$ (in the sense of
\S\ref{3d lams}) which we call $\ac_f$. We think of $\ac_f$ as a {\em
completion} of $\al_f$.

Going through  the same construction replacing $\C_*$-action with $S^1$-action,
we obtain the hyperbolic 3D-extension 
$\hc_f\subset \hat\k/S^1\equiv \hat\kh$, 
with the hyperbolic  action of $\hat f$.  

\smallskip\noindent
{\bf Remark.} The laminated space $\ac_f$ inherits from the universal space
$\uu$ the quality of a metrizable separable space.
Moreover, it has a natural uniform structure coming from the linear
structure of $\uu$, and complete with respect to it. 
However, $\ac_f$ may presumably inherit from $\uu$ also the bad fortune
of not being locally compact.

\subsection{Induced topology}
Let us give a dynamical description of the new laminar topology
$\al_f$ on the leaf space $\an_f$.

By a local leaf $L_{loc}(\hat z,V)$ over a domain $V\subset
\bar\C$ containing $\pi(\hat z)$ we mean the connected component 
of $L(\hat z)\cap\pi^{-1}V$ containing $\hat z$. 

\begin{proposition}{convergence}
A sequence  $\hat z^n\in \al_f$ converges to $\hat \zeta\in\al_f$
if and only if the following hold:
\begin{itemize}
\item[(i)] $\hat z^n\to \hat \zeta$ in the natural topology,
\item[(ii)] For any $N$ and a neighborhood $V$ of $\zeta_{-N}$, if the
local leaf $L_{loc}(\hat f^{-N}\hat\zeta,V)$ is univalent over $V$, then for
sufficiently large $n$, $L_{loc}(\hat f^{-N}\hat z^n,V)$ is univalent
over $V$ as well. 
\end{itemize}
\end{proposition}

We remark that convergence to a point in $\ac_f\sm\al_f$ is more
subtle to characterize in general. Proposition \ref{dynamical
description} does this in the post-critically non-recurrent case. 

\begin{pf} 
Assume first that conditions (i) and (ii) are satisfied.

Represent $\hat\zeta$ as a sequence $\hat\psi = \{\psi_j\}$ in
$\hat\k_f$, and each $\hat z^n$ as $\hat\phi^n = \{\phi^n_j\}$, in
particular noting $\psi_j(0) = \zeta_j$ and $\phi^n_j(0) = z^n_j$.

The statement that $L_{loc}(\hat f^{-N}\hat\zeta,V)$ is univalent
over $V$ is equivalent to saying that $\psi_{-N}$ is univalent in the
component $W$ of $\psi^{-1}_N(V)$ containing $0$ (henceforth we say
``$\psi_{-N}$ is locally univalent over $V$''),
and similarly for 
$\hat f^{-N}\hat z^n$ and $\phi^n_{-N}$. 

Whenever, for some $n,N,V$,  both $\psi_{-N}$ and $\phi^n_{-N}$ are
locally univalent over $V$, there is a unique univalent map
$h^n:W\to \C$ satisfying $\psi_{-N} = \phi^n_{-N}\circ h^n$ on $W$.
Note that, applying $f$ a finite number of times, we have
\begin{equation}
\label{h defn}
\psi_{-j} = \phi^n_{-j} \circ h^n
\end{equation}
on $W$ for any $j\le N$.
Thus if we increase $V$ or change $N$ (but preserve the local
univalence), we obtain $h^n$ equal to the original on the original
domain, or in other words $h^n$ is locally independent of $N$ and $V$.
Choose the normalization of each $\hat\phi^n$ (mod $\C_*$) so that
$(h^n)'(0) = 1$. 

Because $\hat \zeta\in\al_f$, for any disk $D_r$ around $0$ there is
some $N(r)$ for which $\psi_{-N}$ is univalent on $D_r$ whenever
$N>N(r)$. Let $V = \psi_{-N}(D_r)$. 
For sufficiently large $n(r)$, by (i) we have that
$\phi^n_{-N}(0) = z^n_{-N}\in V$, and by (ii) that $\phi^n_{-N}$ is
locally univalent over $V$. Thus we have $h^n$ defined as above on
$D_r$ if $n>n(r)$.

If we let $x^n$ be the preimage of $z^n_{-N}$ in $D_r$ by $\psi_{-N}$
(note that $x^n$ is independent of $N$ if $N>N(r)$), we note that 
$h^n(x^n) = 0$, and by (i), $x^n\to 0$ as $n\to\infty$.

Thus, the sequence of functions $h^n$ now has these properties:
$(h^n)'(0) = 1$, $h^n(x^n) = 0$ where $\lim_{n\to\infty}x^n= 0$, and $h^n$ is
eventually defined on any compact set in $\C$. It is an application of
the Koebe distortion lemma now to show that $h^n$ converges to the
identity on compact sets, and indeed that the image of $h^n$
eventually contains any compact set in $\C$ so that $(h^n)^{-1}$
converges to the identity on compact sets as well. 

Applying (\ref{h defn}) for any $j$, we conclude that 
$\phi^n_j\to\psi_j$ on compact sets for all $j$. Thus $\hat z^n\to\hat
\zeta$ in $\al_f$.

\smallskip

Conversely, let $\hat z^n\in\al_f$ converge to
$\hat \zeta\in \al_f$. 
Assertion (i) is obvious; it is just the statement that $p:\hat\k_f \to
\an_f$ is continuous, which we have already observed. 

For (ii), let $V$ be a neighborhood of $\zeta_{-N}$ such that 
$L_{loc}(\hat f^{-N}\hat\zeta,V)$ is univalent over $V$, and let
$\psi_{-N}:\C \to \bar\C$ be as above. If $W$ is the component of
$\psi^{-1}_{-N}(V)$ containing $0$ we then have $\psi_{-N}$ univalent
on $W$. 

Because a slight enlargement $V'$ of $V$ (so that $V\Subset V'$)
pulls back along the rest of
$\hat \zeta$ with bounded branching (by definition of $\an_f$), it
follows that $W$ has compact closure in $\C$. Let $W\Subset
W'\Subset W''$ be a pair of 
enlargements of $W$, also with compact closure. By definition of
convergence in $\al_f$, there are representatives
$\hat\varphi^n = \{\varphi^n_j\}$ of $\hat z^n$ in $\hat\k_f$ such that
$\varphi^n_{-N}$ converges on $W''$ to
$\psi_{-N}$. It follows that for large enough $n$, $\varphi^n_{-N}$ is
univalent on $W'$ and $\varphi^n_{-N}(W')$ contains $V$, and thus (ii) holds.
\end{pf}

\subsection{Uniqueness}

Let us now consider, for an abstract affine orbifold extension $\ac$
of $f$ in the sense of Section \ref{axiomatic}, what properties force
it to be equal to our universal construction $\hat\ka$.

There is a natural map $I:\ac \to \hat\ka$, defined similarly to $\iota$:
for any $\bold z\in\ac$, let $\phi:\C\to L(\bold z)$ be (the inverse
of) any affine chart for the leaf of $\bold z$ that takes $0$ to $\bold
z$. Then the sequence $\{[\pi_n\circ\phi]\}$ gives a well-defined
element of $\hat\ka$, where $\pi_n$ are the projections of $\ac$
to $\bar\C$. The difference between $I$ and $\iota$ is that
$I$ is automatically continuous,
because of the transverse continuity of the affine structures in $\ac$. 

We now observe that $I(\ac)$ is equal to $\hat\ka$
if the following conditions hold:
\begin{enumerate}
\item The map $I$ is an embedding,
\item $\al_f$ is dense in $I(\ac)$, and
\item $I(\ac)$ is closed
\end{enumerate}
In particular, the first
condition reduces to checking that $I$ is both injective and proper:
i.e. that an element of $\ac$ is determined uniquely by the sequence
of functions $\pi_n\circ\phi$, and that convergence in $\ac$ follows from 
convergence of the sequence of functions.

For the construction of Section \ref{pcf orbifolds}
of orbifold laminations for post-critically finite maps, these
properties evidently hold, and therefore the general construction
produces the same object.

\subsection{Minimality}
Let us show that the laminations
we constructed are minimal. Note that this does not follow from
Lemma \ref{affine leaves are dense} 
since topology of $\ac_f$ is stronger than that of $\NN_f$.
 
\begin{proposition}{minimality}
The laminations $\ac_f$ and  $\hc_f$  are  minimal except for the
Chebyshev and Latt\`es examples. In those cases the lamination
becomes minimal after removing the isolated invariant leaf associated
with a post-critical fixed point.
\end{proposition}
Hence every open set $K$ of either lamination
contains  a global cross-section for it 
(except the isolated leaves in the above special cases).

\begin{pf} Clearly it is enough to consider $\ac_f$. Since
$\al_f$ is dense in $\ac_f$, it suffices to demonstrate density of
leaves in $\al_f$.

Let us first show that any invariant leaf $L$ is dense.
Take a point $\hat z=\{z_0,z_{-1},\dots\}$ in $\al_f$, 
and a finitely branched pullback of neighborhoods 
$\{U_0,U_{-1},\dots\}$ along it. In the case when $f$ is
Latt\`es or Chebyshev assume that $\hat z$ is not 
a postcritical fixed point. Then Proposition \ref{C and L}
and the expansion property of $f$ on the Julia set
easily yield existence of 
a limit point $a\in\bar\C$ for $\hat z$ and $\eps>0$
such that one of the local leaves $L_{loc}$ over 
$D(a,\eps)$ is not branched. 

For sufficiently large $N$, $U_{-N}$ pulls back univalently along the
rest of $\hat z$, and 
by the Shrinking Lemma, there is a sequence $N_i\to\infty$ 
such that $U_{-N_i}\subset D(a,\eps)$.

Thus $L_{loc}$ is 
univalent over $U_{-N_i}$. 
Let $\hat b^i$ be the point on this local leaf which projects to
$z_{-N_i}$. Then by Proposition \ref{convergence}, the sequence
$\hat f^{N_i} \hat b^i\in L$ converges to $\hat z$ in the $\ac_f$ topology
as $i\to\infty$, which proves density of $L$.   

Replacing $f$ by its iterate, we conclude that every periodic leaf is dense
in $\ac_f$. 

Let us now show that every leaf $L(\hat z)\subset \al_f$ 
accumulates on some periodic leaf.
 To this  end take five periodic points $\alpha_k$ and associated periodic
leaves $L_k\equiv L(\alpha_k)$. Select five disjoint topological discs
$D_k\ni \alpha_k$. By Ahlfors' Five Islands Theorem (see \cite{tsuji},
 Theorem VI.8), for any $n$, each $\hat f^{-n} L(\hat z)$ has a univalent 
local leaf   over one of the domains $D_k$. Take a $k$ for which this
happens for infinitely many $n$'s. Then by the same argument as above 
$L(\hat z)$ accumulates on the periodic leaf $L_k$.
\end{pf}

\section{Convex-cocompactness, non-recurrence and conical points}
\label{convex cocompactness section}
Define the Julia set $\jc_f$ in $\ac_f$ to be the pullback of $J_f$ by
$\pi:\ac_f \to \bar\C$.
Let $\jl_f$ denote $\jc_f \intersect \al_f$. Note that $\jl_f$ and
$\jn_f$ have the same underlying set and different topologies, and
that $\jc_f$ is the closure of $\jl_f$.

We say that $f$ is {\em convex cocompact} if the quotient
$\CC(\jc_f)/\hat f$ of the convex hull is compact.
In this section we prove several criteria for convex
cocompactness. The main criterion is the following:

\begin{theorem}{convex cocompactness1} A rational map $f$ is convex
cocompact if and only if it is postcritically non-recurrent and has no
parabolic points. 
\end{theorem}  

\noindent{\bf Remark.} This criterion is closely related to the 
``John domain criterion'' given by Carleson, Jones and Yoccoz
for polynomials \cite{carleson-jones-yoccoz}.   See also McMullen
\cite{mcmullen:john} for the connection between convex-cocompactness
and the John condition in the setting of Kleinian groups.

\subsection{Convergence and compactness}
For a critically non-recurrent map $f$ without parabolics, we can give
a dynamical criterion for convergence in $\ac_f$ (note that
Proposition \ref{convergence} only applied to convergence within
$\al_f$. This criterion includes the possibility  that a bounded
amount of branching
persists in the limit and yields a point outside $\al_f$).
Let $p: \ac_f\ra \al_f$ denote the natural projection.

\begin{proposition}{dynamical description} Let $f$ be critically
non-recurrent without parabolics. 
A sequence of points $\hat z^n\in \al_f$ 
converges to $\bzeta\in \ac_f$, with $p(\bzeta)=\hat \zeta$,
 if and only if
\begin{enumerate}
\item[(i)]
$\hat z_n\to \hat \zeta$ 
   in the natural topology and
\item[(ii)]
For any $N$ and a neighborhood $V$ of $\zeta_{-N}$, if the
local leaf $L_{loc}(\hat f^{-N}\hat\zeta,V)$ is univalent over $V$,
then the following holds:

There is a finite set of points $\{c_k\} \subset V$ such that for any
neighborhood $\Omega$ of $\{ c_k\}$ there exists $M=M(\Omega)$
so that, if $n>M$, the local leaf $L_{loc}(\hat z^n,V\sm\Omega)$
covers $V\sm\Omega$ without branching, and for any $n,m>M$ the
coverings are topologically equivalent.
\end{enumerate}
Moreover, the projection $L(\bzeta)\ra L(\hat \zeta)$ is
a finitely branched covering with uniformly bounded degree.  
\end{proposition}

\begin{pf}
By Ma\~n\'e's Theorem,  there is a neighborhood $W$ of $J_f$, and 
$\eps_0>0$ and $K_0$ with the following property:
for any backward trajectory
$\hat z=\{z_0,z_{-1},\dots\}\in\NN_f$ with $z_0\in W$, 
the pullback of the disk $D(z_0, \eps_0)$
along $ \hat z$ branches at most $K_0$ times. 
(Compare the proofs of Lemma \ref{nonrec1} and Proposition \ref{nonrec3}).

Furthermore,  for any $\hat z$ which is not an attracting cycle, there
is an $N_0(\hat z)$ such that $z_{-n} \in W$ for $n>N_0$.

Assuming that (i) and (ii) hold, represent $\hat\zeta$ using a
sequence $\{\psi_{-N}\}\in\hat\k_f$. For any disk $D\subset \C$, for
large enough $N>N_0$ the map $\psi_{-N}$ is univalent in $D$ and has
image in $W$, and in fact in $D(\zeta_{-N},\ep_0)$. For sufficiently
close $\hat z^n$ to $\hat\zeta$, 
$z^n_{-N}$ is also in $W$, and hence the pullback of $D(z^n_{-N},\ep_0)$
along the rest of $\hat z^n$ has uniformly bounded branching. 

Condition (ii) now gives a branched cover of $D$ which is
conformally equivalent to 
the coverings of 
$L_{loc}(\hat f^{-N}\hat z^n,\psi_{-N}(D))\to \psi_{-N}(D)$
away from a small neighborhood of the critical points, for large
enough $n$. This branching is uniformly bounded no matter how large
$D$ is taken, so we obtain a polynomial $h:\C \to \C$. The
sequence $\{\psi_{-N}\circ h\}\in\hat\k_f$ will represent
 the limit $\bzeta$ of the
$\hat z^n$ in $\ac_f$, by an argument similar to that in Proposition
\ref{convergence}, where condition (ii) keeps the branching
consistent. 

More precisely, let $D'=h^{-1}D$ and assume that $D'$ is large enough
that the (finite) set $C$ of critical points of $h$
is separated from 
$\boundary D'$ by an annulus of modulus $M>0$. For each $n$ represent 
$\hat z^n$ by a sequence of functions
 $\phi_{-N}^n:\C\to\bar \C$, normalized so its
1-jet agrees with $\psi_{-N}\circ h$ at a fixed non-critical point
$w\in D'$.
Let $Y$ be a neighborhood of $C$ so that $D'\sm Y$ contains an annulus
of modulus $M$ around each puncture. 
Then condition (ii) gives, for large enough $n$, a univalent map 
$u_n : D'\sm Y\to \C$ such that $\psi_{-N} \circ h = \phi^n_{-N}\circ u_n$,
and $u_n(w) = w, u'_n(w) =1$. Note that $u_n$, once defined on $D'\sm
Y$, remains the same there as we enlarge $D'$, shrink $Y$ and
increase $N$, and that $\psi_{-N}\circ h = \phi^n_{-N} \circ u_n$ wherever
it is defined. Again using Koebe distortion (this time on a multiply
connected domain), we have $u_n \to \id$ on compact subsets of 
$\C\sm C$. It follows that for every $N$,
$\phi^n_{-N} \to \psi_{-N} \circ h$, 
so that $\hat z^n\to \bzeta$ as $n\to \infty$.

Moreover, $p: L(\bzeta)\ra L(\hat\zeta)$ is a finitely branched covering
with bounded degree since $h$ is.  

Conversely, suppose that the sequence $\hat z^n$ converges in $\ac_f$.
Since by the same discussion the branching over each disk 
$D(z,\ep_0)$, $z\in J_f$,
 is eventually uniformly bounded, there must be some
subsequence of the $\hat z^n$ for which the branching converges in the
sense of (ii), and so the limit is equal to the limit defined in the
previous paragraph. It follows that the same holds for any
subsequence, so that in fact (ii) holds for the whole sequence. 
\end{pf}

\begin{corollary}{compactness}
Let $f$ be critically non-recurrent without parabolics. 
A set $K\subset \al_f$ is pre-compact in $\ac_f$ if and only if 
its closure in the natural extension $\NN_f$ does not contain 
attracting cycles. 
\end{corollary}

\begin{pf}  If a sequence $\{\hat z^n\}\subset\al_f$ does not accumulate
on  attracting cycles
then Mane's Theorem easily yields existence of a subsequence
satisfying (i) and (ii) of the previous proposition. 
\end{pf}

\subsection{Proof of Theorem \ref{convex cocompactness1}}
Let us split the proof  into two steps represented by the following
two criteria. 

\begin{lemma}{compactness of the Julia set}
The Julia set $\jc_f$ is compact if and only if $f$ is critically
non-recurrent and has no parabolic points.
\end{lemma}

\begin{pf}  
If $\jc_f$ is compact then $\jn_f$ is compact in $\an_f$, since
$\jn_f = p(\jc_f)$ where $p:\ac_f \to \an_f$ is the natural continuous
projection. 
Hence by Corollary \ref{Julia compact} $f$ is critically
non-recurrent without parabolic points. 

Vice versa, if $f$ is critically non-recurrent without parabolic
points then compactness of $\jc_f$  
follows from Corollary \ref{compactness}.
\end{pf}

\begin{proposition}{convex cocompactness2}
A rational map $f$ is convex
cocompact if and only if the Julia set $\jc_f$ is compact. 
\end{proposition}

\begin{pf}
Let $\VV=\VV(\jc)$ denote the ``curtain'' over $\jc=\jc_f$ in $\hc$.
That is,
the union of vertical geodesics over points of $\jc$. We will first
show that $\VV/\hat f$ is compact if and only if $\jc$ is
compact. 

Observing that $\jc$ is just the quotient $\VV/e^\R$ by the vertical
geodesic flow, we may view this equivalence in slightly generalized
terms:

Let $X$ be a Hausdorff space admitting commuting actions by two closed
non trivial subgroups $G$ and $H$ of $\R$. Let $G$ act on the left and
$H$ on the right, for clarity. Suppose $G$ and $H$ both act properly,
and that they are {\em coherent} in the sense of Lemma
\ref{coherence}: any $x,y\in X$ are contained in neighborhoods $U_x, U_y$
for which $g_iU_x \intersect U_yh_i \ne \emptyset$ only if $g_i,h_i$
both remain bounded, both go to $+\infty$ or to $-\infty$.
We claim that $G\backslash X$ is compact if and only if $X/H$ is
compact.

Suppose without loss of generality that $G\backslash X$ is compact,
and let $K\subset X$ be a compact fundamental domain, i.e. $GK = X$.
Let $x\in K$ and consider the {\em positive return time} $g_+(x)$ for
the orbit $xH$ to return to $KH$ under $G$. That is, let $g_+$ be the
smallest positive element of $G$ such that $g_+ xH \intersect KH \ne 0$.
We claim this is bounded for $x\in K$. Choose $h<0$ in $H$
sufficiently far from $0$ that $Kh\intersect K = \emptyset$ (by the
proper action of $H$), and that $gK \intersect Kh \ne \emptyset$ only
for $g<0$ (this is possible by coherence, after covering $K$ with a
finite number of neighborhoods $U_p$). 

Thus the point $xh$ is not in $K$ so there is some $g_+>0$ such that
$g_+ xh \in K$. Thus $g_+Kh \intersect K \ne \emptyset$ for each
$g_+$, so that fixing $h$ we have an upper bound for $g_+$ independent
of $x$, by the proper action of $G$.

Reversing the signs in the argument we also obtain a bounded {\em
negative} return time for every $x\in K$. We conclude that, in the
action of $G$ on $X/H$, every point has a bounded negative and
positive return time to the projection $KH$ of $K$. Since $X/H$ is
covered by $G$-translates of $KH$, it follows that there is a bounded
subset $I$ of $G$  such that $IKH$ covers $X/H$. Thus $X/H$ is
compact.

\medskip

In our situation the groups are $\Z$ and $\R$, and we conclude that
$ \VV/\hat f$ is compact if and only if $\VV/e^\R = \jc$ is
compact (note: it would be more consistent to denote the first quotient 
$\hat f\backslash\VV$). It remains to check that compactness of
$\VV/\hat f$ is equivalent to compactness of the convex core quotient. Since
the curtain is closed in the convex core, one implication is clear. 
Conversely, if we know that $ \VV/\hat f$ is compact, 
we need only to observe that the convex core
lies in a bounded neighborhood of the curtain. That
is, let $p\in \CC$ be some point, represented as $(z,t)$ in a
half-space model of the leaf $L_{hyp}(p)$ of $\hc$.
 If $z_0$ is the nearest point
to $z$ in the local Julia set $\jc\intersect L$, then $t>|z-z_0|$
because otherwise $p$ lies in a hemisphere over $z$ disjoint from
$\jc$, and therefore outside the convex hull. It follows that the
hyperbolic distance from $p$ to $(z_0,t)$, which lies in $\VV$,
is less than 1.

It is easy to check that a leafwise 1-neighborhood of a compact subset
of a hyperbolic 3-lamination is itself compact, so this concludes the proof.
\end{pf}

\subsection{Conical points}
Given a point $\bold{z}\in\jc$, 
let  $\gamma_{\bold{z}}$
be the vertical geodesic in $\hc_f$  terminating at $\bold{z}$.
By analogy with Kleinian groups, let us say that $\bz\in \ac$ is 
a {\it conical point} if the projection of the geodesic $\gamma_\bz$ 
to the quotient lamination $\hc/\hat f$ 
does not escape to infinity (which means that there is a sequence of 
points $p_n\in \gamma_\bz$ tending to $\bz$ whose 
 projection to $\hc/\hat f$ converges).
Note that in this definition
the vertical geodesic can be replaced by any geodsic terminating at $\bz$
since all of them are asymptotic in the hyperbolic metric.

Equivalently, $\bz\in\ac_f$ is conical iff its forward orbit $\{\hat
f^n\bz\}_{n=0}^\infty $  
is non-escaping in $\ac_f$, that is, the $\omega$-limit set
$\omega(\bz)$ is non-empty. Indeed, given
two proper commuting group actions $G$ and $H$ on a space $X$, 
the $G$-orbit of a
point $x\in X$ is non-escaping in the quotient by $H$ if and only if its
$H$-orbit is non-escaping in the quotient by $G$ (since either is
equivalent to non-escaping of the double orbit $GxH$ in $X$). 
In our situation we have
a $\Z$-action by $\hat f$ on $\hc_f$, and the $\R$-action of the
vertical geodesic flow (as in Section \ref{3d lams}). The 
directionality of our statement (forward $\hat f$-orbits accumulate in
$\hc_f/e^{\R} = \ac_f$ if 
and only if backward $\R$-orbits accumulate in 
$\hc_f/\hat f$) comes directly from the coherence of the
actions, lemma \ref{coherence}.

Let $\Lambda=\Lambda_f$ denote the set of conical points.

\medskip

We further note that the property of being conical depends only on the
projection to $\bar\C$.
Let us say that a  set $ X\subset \ac$ is {\it fiber saturated} if
$X = \pi^{-1}(\pi(X))$. The reason is that the fibers play the role 
of local stable manifolds for $\hat f$ (the proof below makes 
precise the sense of this statement).

\begin{proposition}{conical is saturated} The set of conical points
 is fiber saturated.
\end{proposition}

\begin{pf} Let us show that  
the $\omega$-limit sets of $\bold{z}$ and $\bold{w}$ are equal, up to
finite branched cover.
Represent $\bold{z}$ and $\bold{w}$ in $\hat\ka_f$ by sequences of meromorphic
functions $\{\phi_n\}$ and $\{\psi_n\}$ such that $\phi_0(0) =
\psi_0(0) = \pi(\bold{z})$. 
To compute the $\omega$-limit sets it suffices to consider just the 
first coordinate functions, $\phi=\phi_0$ and $\psi=\psi_0$.
Suppose first that $\pi$ is non-singular at $\bold{z}$ and $\bold{w}$,
so that we may assume
$\psi'(0)=\phi'(0)=1$. 

Now suppose that $h$ is a limit point of $f^n\circ\phi$ in $\ua$.
This means that for some sequence $n_i$, and $\lambda_i\in \C_*$,
$f^{n_i}\circ\phi\circ \lambda_i$ converges to $h$ on compact subsets
of $\C$. By Lemma \ref{coherence}, we know that $|\lambda_i| \to 0$.
Now fixing a disk $D\subset \C$ around $0$, we see that for $i$
sufficiently large, $\phi$ and $\psi$ are both invertible in
$\lambda_i D$, and by the Koebe distortion theorem, the combined map
$(\psi\circ\lambda_i)^{-1}(\phi\circ\lambda_i)$ converges to the 
identity on $D$. It follows that $f^{n_i}\circ\psi\circ\lambda_i$ also
converges to $h$. 

If, on the other hand, $\psi$ and/or  $\phi$ have branched points at
$0$, say with degrees $k$ and $m$ respectively, let $d=lcm(k,m)$ and write
$\til\psi = \psi\circ b_{d/k}$ and
$\til\phi = \phi\circ b_{d/m}$, where $b_j(z) = z^j$. Now $\til\psi$
and $\til\phi$ both have degree $d$ at $0$, and for small
$|\lambda_i|$ we still make sense of 
$(\til\psi\circ\lambda_i)^{-1}(\til\phi\circ\lambda_i)$ as a univalent
map. Hence the Koebe distortion argument goes through and we may
conclude that $\{f^n\circ\til\psi\}$ and $\{f^n\circ\til\phi\}$ 
have the same $\omega$-limit points in $\ua$.
\end{pf}

Let $\Delta=\pi\Lambda\subset J_f\subset\bar\C$.
 By the above proposition, it is justified
to call the points of this set conical as well. 
Let us show that it is trapped in between two well-studied sets. 

First, let $\Delta_1$ denote the set of points $z\in J_f$
such that there is an $r>0$ and a sequence $n_i\to \infty$
(depending on $z$)
such that the multi-valued inverse branch 
$f^{-n_i}: D(f^{n_i}z, r)\ra U_i\ni z$ has a bounded degree
(compare \cite{lyubich:typical}).

The second set, $\Delta_0$ is the union of all expanding 
subsets of the Julia set (a compact invariant set $X\subset\bar\C$
is called {\it expanding} if $f: X\ra X$ is surjective and
some iterate $f^n|X$ has spherical derivative strictly greater than 1).   

\begin{proposition}{trap}
$\Delta_0\subset\Delta\subset\Delta_1$.
\end{proposition}

\begin{pf} Let us start with the right-hand inclusion.
Let $z=\pi(\bz)$ for $\bz\in\Lambda_f$. Then there exists a
sequence $n_i\to\infty$ such that $\hat f^{n_i}\bz\to \bzeta\in\ac_f$.
Translation of this to the language of meromorphic functions
 provides us with a desired family of inverse branches with bounded
degree.  

For the left-hand inclusion,  take a point $z$ in an expanding set
$X\subset\bar\C$. First notice that by Lemma \ref{big disks}
any backward orbit $\hat z$
in the invariant lift  $\hat X\subset\NN_f$ belongs to a parabolic leaf.
Then, take any convergent subsequence $f^{n_i}\hat z\to \hat \zeta\in \hat X$
in the natural topology and apply Proposition \ref{dynamical description}
to see that it is convergent to the same point in the laminar topology as well
(the local leaves in condition (ii) of this proposition can be
selected univalent).  
\end{pf}

\begin{proposition}{cocompactness and conical pts}
If $f$ is convex cocompact then all points of the
Julia set $\jc_f$ are conical.

Conversely, if the lamination $\ac_f$ is locally compact and all points
of the Julia set $\jc_f$ are conical then $f$ is convex cocompact.
\end{proposition}

\begin{pf} Assume $f$ is convex cocompact, that is, the 
convex core $\CC_f/\hat f$ is compact. Since
$\gamma_\bz\subset\CC_f$ 
for any $\bz\in\jc_f$, the conical property of $\bz$ follows. 

\smallskip 

For the converse,
suppose the lamination $\ac_f$ is locally compact.
Then there  is a compact set $K$ with non-empty interior.
By Proposition \ref{minimality}, $K$ meets every leaf of the
lamination. Since the set $K\intersect L_{hyp}(p)$ is
closed in the intrinsic leaf topology, for any $p\in \hc$, 
there is a length minimizing geodesic $\Gamma_p$ joining
$p$ to $K$. Let $\dist(p,K)$ denote the 
hyperbolic length of this geodesic.
It can be also defined as the infimum of lenths of all curves joining
$p$ and $K$.

Given a set $X\subset \hc$, let $N(X,r)=\{p\in\hc: \dist(p,X)< r\}$
denote the leafwise $R$-neighborhood of $X$.
Then any compact set $Q\subset  \hc$ is covered by some $N(K,R)$.
Indeed, for every $q\in Q$ there is a curve $\gamma$ joining $q$ with
a point in the interior of 
$K$. If the length of this curve is $r$ then all
sufficiently nearby points
can be joined with $\inter K$ by a curve of length less than $r+\eps$.
Now compactness of $Q$ yields the statement.

Note also that the space $\GG$ 
of one-sided geodesics beginning in $K$ is parametrized by the unit
tangent bundle over $K$ and hence is compact. 

Assuming that the convex core $\CC_f/\hat f$ is
not compact let us construct in it an escaping geodesic.
Consider a sequence of points $q_n\in \CC/\hat f$
escaping to $\infty$ and the corresponding minimizing geodesics
$\Gamma_n\equiv \Gamma_{q_n}$. 
 By compactness of $\GG$, there is a limit
geodesic $\Gamma$ beginning at $K$. 
Let us show that this geodesic escapes to $\infty$.

Indeed, otherwise there is a compact set $Q\subset \CC/\hat f$ 
which  $\Gamma$ does not escape. Let us consider
the leafwise $1$-neighborhood $N(Q,1)$ of $Q$. 
Its closure is compact and hence is contained in some 
leafwise neighborhood $N(K,R)$ of $K$.

Since the $\Gamma_n$ accumulate on $\Gamma$, for some $n$ there are
two points $a,b\in \Gamma_n\cap N(Q,1)$ such that the distance berween
them along $\Gamma_n$ is greater than $R$. On the other hand, there is
a curve from $b$ to $K$ of length less than $R$ which contradicts the
minimality of $\Gamma_n$.
\end{pf}

Let us say that a set $X\subset\ac$ is {\it locally fiber saturated} if
for any point $p\in X$ there is a box neighborhood $U\ni x$  such that
if $q\in U\cap X$ then the whole fiber $\pi^{-1}(\pi q)\cap U$ belongs to $X$. 
We can then say that such a set $ X$  is measurable and has  ``zero",
``positive" or ``full" measure if   
the corresponding property is satisfied leafwise, that is for its
intersection with every leaf.
Note that these notions are well defined on the
affine leaves though the Lebesgue measure is not. Note also that they
don't require any transversal measure.
 
Given a measurable locally fiber saturated 
$\hat f$-invariant set $X\subset \ac$, we say that $\hat f|X$ is {\it ergodic}
if every measurable locally fiber saturated $\hat f$-invariant subset
$Y\subset X$ has either zero or full measure. 

An {\it invariant line field} on $\ac$ is a measurable real
one-dimensional distribution 
in the tangent bundle $T\ac$ over a set of positive measure, which is 
transversally continuous in measure and invariant under $\hat f$.
We say that the line field is {\it constant} if it is constant in the
affine chart on  
any leaf. Note, if we are considering an orbifold leaf then this must
take place in a finite cover -- this allows the case of an orbifold
point of order two,  and a line field with a simple pole singularity. 
This is exactly what occurs for the deformable Latt\`es example.

Given a measurable set $X$ and a set of positive  Lebesgue measure $Y$
on an affine leaf $L$,  
let $\dens(X|Y)=\meas(X\cap Y)/\meas(Y)$ (note that this is a well
defined quantity).  
Let us formulate some general ergodic properties of the conical set:

\begin{proposition}{ergodic propertis}

\begin{itemize}
\item
The set $\Lambda_f$ of conical points has either zero or full Lebesgue measure.
\item
In the latter case $f$ is ergodic, except for the Latt\`es
examples.
\item
Any  invariant line field on  $\Lambda_f$ is constant, except for the isolated
leaves of Latt\`es examples.
\end{itemize}
\end{proposition}

\noindent{\it Proof 1.} This proof demonstrates how the blow-up  method works
in the lamination context.

Take any invariant locally fiber saturated set   $X\subset \Lambda_f$ of
positive measure.  Then
$X\intersect L$ has positive measure for any leaf $L\subset \ac$. 
Take a leaf $L$ and a density point $\bz$ of $X$ in $L$.
Since $\bz$ is conical, there is a convergent sequence $\hat
f^{n(k)}\bz\to\bzeta$.   
Take an arbitrary round disc $D\subset L(\bz)$ 
and a box neigborhood  $D\times T$  of $\bz$. Let
 $\bzeta=(\zeta,\tau),\quad \bz_{n(k)}=(z_{k},t_k)\in D\times T$. 

Let us consider round discs $\Delta_k=\hat f^{-n(k)}(D\times t_k)$ on
the leaf $L(\bz)$.  By the Shrinking Lemma, they shrink to $\bz$ and
hence $\dens(X|\Delta_k)\to 1$.  Since $\hat f$ is leafwise affine,
$\dens(X|(D\times t_k))\to 1$ as $k\to \infty$.  Since $X$ is 
fiber saturated, $\dens(X|(D\times\tau))=1$.  Since the disc $D$ is
arbitrarily big, $X$ has full measure on the leaf $L(\bzeta)$.

Since the leaf $L(\bzeta)$ is dense in $\ac_f$ (except the isolated
leaves in the Latt\`es examples) and $X$ is locally fiber saturated, it has
full measure on every leaf. This proves the first two statements,
except for the Latt\`es examples.

Finally, the first statement holds for the Latt\`es examples since
 $\Lambda_f=\ac_f$ 
by Theorem \ref{convex cocompactness1} and Proposition \ref{cocompactness and
conical pts}. The second statement fails for the trivial reason
that the isolated leaf is an invariant locally fiber saturated subset of $\ac_f$.
However, the previous argument shows that this leaf and its complement are
the only subsets like this.

If now $X$ supports an invariant line field $\mu$, take  $\bz$ to be a Lebesgue continuity
point for this field on the leaf $L(\bz)$, so that $\mu$ is almost constant on
$\Delta_k\cap X\sm Y$ where $\dens(Y|\Delta_k)\to 0$ as $k\to\infty$.  It follows that
$\mu|(D\times t_k)$ accumulates in measure on constant line fields. Since  $\mu$
is transversally continuous in measure, $\mu|(D\times\tau)$ is constant almost everywhere,
and hence almost everywhere
 on the leaf $L(\bzeta)$.  As this leaf is dense in $\ac$, except for the isolated leaves
of Latt\`es examples, the last statement follows as well.  \QED

\smallskip\noindent {\it Proof 2.} This proof (for first two statements only)
exploits Ahlfors' harmonic extension method. Namely, 
let $X\subset\Lambda$ be a locally fiber saturated set of positive measure.
Then  we construct a harmonic function on $\hc/\hat f$ by solving a
Dirichlet problem on each leaf. That is, given an affine leaf $L\subset \ac$
and the attached hyperbolic leaf $H_L\subset\hc$,
 construct
harmonic $h:H_L\to\R_+$ whose boundary values are 0 on $\fc$ and $1$
on $\Lambda\cap L$.  It is perhaps best to think of $h(x)$ for $x\in H_L$  as
the area of $\Lambda\intersect L$ as measured in the ``visual metric'' at
$x$. That is, we map $H_L\union L$ by M\"obius transformation to the unit ball
taking $x$ to $0$, and measure the area of the image of $\jc\intersect
L$ on the unit sphere. This is the same as integrating the Poisson
kernel against the characteristic function of $\Lambda\intersect L$.

One must check that $h$ is continuous (in the transverse direction).
But if we fix a box neighborhood $T\times D$ for a large $D$ in $L$ 
(in the orbifold case this should be the finite cover of a box
neighborhood), then for points near $x$ (in the transverse direction)
the visual measure induced on 
the leaves near $L$ changes continuously on $D$ (and if we choose $D$
large enough, is very small on the complement of $D$ in each leaf).
The intersections of $\Lambda$ with nearby leaves 
is a continuous family of analytic branched covers. It
follows that area measure on $\Lambda$ varies continuously in the
transverse direction, and therefore so does its integral with respect
to the visual measure.

Consider now a density point  $\bz\in X$ and the geodesic $\gamma_\bz\subset L(\bz)$
terminating at this point. Then $h(p)\to 1$ as $p\to\bz$ along $\gamma_\bz$,
 since the visual are of $X$ as seen from $p$ is going 1.   

Observing also that $h$ is invariant by $\hat f$, we obtain a
continuous leafwise harmonic function $g$ on the quotient $\hc/\hat
f$. Since $\bz$ is conical, the projection of $\gamma_\bz$ to
$\hc/\hat f$ has a limit point $q$. By continuity, $g(q)=1$. By the
Maximum Principle, $g$ is identically equal to 1 on the whole leaf
$L(q)$.  By Proposition \ref{minimality}, this leaf is dense, except
for the Latt\`es examples, and thus $g$ is identically equal to 1 on
the whole lamination. It follows that $X$ has full measure.  \QED

\noindent{\bf Remark.} Given Proposition \ref{trap}, the results of the above
proposition are not really new (compare \cite{lyubich:typical},
\cite[Lemma 10]{BL}, \cite[Theorem 3.9]{mcmullen:renormalization}).
However, the laminations give a new insight on them, and strengthen 
the connection to the corresponding results for Kleinian groups
\marginpar{refs?}.

\begin{corollary}{no line fields}   
If $f$ is not Latt\`es, then
there are no invariant line fields on $\Lambda_f$  
which come from the sphere $\C$.
\end{corollary}

\begin{pf} It is easy to see that one can  always find  two
leaves $L(\hat z_1)$ and $L(\hat z_2)$, 
with $\pi(\hat z_1)=\pi(\hat z_2)\equiv z$ such that $L(\hat z_1)$ is branched
at  $\hat z_1$ while $L(\hat z_2)$ is regular at $\hat z_2$.  Then the push-forward of
a constant line field from $L(\hat z_1)$ has a singularity, while the 
push-forward from $L(\hat z_2)$ does not. 

The only case when this does not lead to a contradiction is when one of
the above leaves is isolated, so that the invariant line field is not necessarily
constant on it. 
But this may happen only for the Latt\`es examples.
\end{pf}

\subsection{Elliptic structure of the Latt\`es examples}\label{Lattes}
Let us show in conclusion how the invariant line field imposes the ``elliptic
structure" of the Latt\`es examples. We have seen that the invariant line field
may exist only if there is an isolated leaf $L^r$. But then there should exist a
non-isolated orbifold leaf $L^s$
with an orbifold-constant line field. 

Considering the projection $\pi: L^s\ra\bar\C$ we see that
the line field on $\bar \C$ is locally (a.e.) the image of the constant
line field under a branched cover. It follows that the branching of
$\pi$ can be at most degree 2, and that the line field on $\bar\C$ can
only have isolated index $-1/2$ (pole) singularities. By the index
theorem on line fields, 
 there must be exactly four of these. Thus
$\bar\C$ has the structure of an orbifold with four order-2 singular
points, the (2,2,2,2) orbifold (this is exactly Thurston's orbifold
for this map). 

Let $X\subset \bar\C$ 
denote the above set of four singular points. It is clearly forward
invariant under $f$.
The property that the leaf $L^r$ is isolated means that all backward
orbits $\hat z$ with $z_0\in X$ eventually escaping $X$ hit a critical point. 
In other words, $\pi: L^s\rightarrow \bar C$ 
is double branched at all points  of $L^s\cap\pi^{-1} X$,
 except the singular periodic point.
Thus this map  is an 
orbifold cover. (See e.g. Thurston \cite{wpt:textbook} or Scott
\cite{scott:survey} for a
discussion of orbifolds and orbifold covers).  

Let $q:  \til L^s\rightarrow L^s$ be the double covering associated to
the orbifold structure of  $L^s$, $\til L^s\approx \C$. It follows that 
$\pi\circ q: \til L^s\rightarrow (\bar\C, X)$
is an orbifold universal cover.  
The group of deck translations for such a cover is generated
by a lattice of translations and the involution $z\mapsto -z$.

Let $m$ be a period of the leaf $L^s$.
Note that $\hat f^m: L^s\rightarrow L^s$ lifts 
(in two ways, because of choice of sign) to
a multiplication map $g: z\mapsto nz$ on $\til L$. 
The constant $n$ must be real,
since $g$ preserves the line field.
On the other hand $g$  
commutes with $\pi\circ q$, so it preserves the lattice.
Hence $n$ is an
integer. In other words the original map $f$ is the projection of an
integral torus endomorphism, i.e. a deformable Latt\`es example.

\section{Quasi-isometries and rigidity}
\label{pcf rigidity section}

\subsection{Rigidity}
In this section we will use the convex-cocompactness of the quotient 
3-lamination to prove  rigidity of 
critically non-recurrent maps without parabolic points,
which extends Thurston's rigidity theorem (see
\cite{douady-hubbard:pcf}). 

\begin{theorem}{pcf rigidity}
Let $f$ and $g$ be two critically non-recurrent rational maps
without parabolic periodic points. 
\begin{enumerate}
\item If $f$ and $g$ are topologically conjugate  then they are  
quasi-conformally conjugate. 
\item If the conjugacy is equivariantly homotopic to conformal on the
   Fatou sets, then $f$ and $g$ are M\"{o}bius conjugate, 
  except for the Latt\`es examples.
\end{enumerate}
In particular, the second case holds automatically when 
the Julia sets of $f$ and $g$ coincide with the whole sphere.
\end{theorem}

\subsection*{Remarks.}
Thurston's proof  of rigidity for post-critically finite maps 
used a contraction principle on a
Teichm\"uller space,
which is another aspect of the connection between
rational maps and Kleinian groups (see \cite{mcmullen:iter,mcmullen:icm90}).

Our proof uses another familiar scheme from both dynamics
and hyperbolic geometry, which is roughly as follows.
In step one, a topological conjugacy is promoted to a quasi-conformal
conjugacy, using some geometric information. In step two, the
quasi-conformal conjugacy is found to be conformal by an
ergodic reasoning, because it induces
an invariant line field on the Julia set. 

In the convex cocompact case,
the topological conjugacy is almost immediately quasi-conformal,
because it gives rise to a homeomorphism on compact sets (the convex
cores), which is automatically a quasi-isometry of the 3-laminations.
This is directly analogous to the proof of Mostow's rigidity theorem
in the case where the Fatou domain is empty, and to Marden's
isomorphism theorem otherwise. 

The second step, absence of invariant line fields, follows from
the properties of the conical limit set given in the previous sections.

\subsection*{Proof}
Let $\ac_f$ and $\ac_g$ be the affine orbifold laminations constructed
 from the natural extensions of $f$ and $g$, and let $\hc_f$ and
$\hc_g$ be the hyperbolic orbifold 3-laminations built over $\ac_f$
and $\ac_g$.

Let $\Phi:\bar \C\to\bar \C$ be the homeomorphism conjugating $f$ to $g$. 
Let $\hat\Phi: \NN_f \to \NN_g$ denote the natural extension of $\Phi$, 
which conjugates the action of $\hat f$ to that of $\hat g$.  
This map admits a continuous extension to a homeomorphism, 
which we also call 
$\hat\Phi$, from $\ac_f$ to $\ac_g$, again conjugating $\hat f$ to $\hat 
g$, and preserving orbifold affine structure. Indeed, 
Proposition \ref{dynamical description} describes 
convergence in $\ac_f$ in dynamical terms which are
respected by topological conjugacy.
(Note: this is not obvious and maybe not true for critically recurrent
maps.)

\comm{ We need to check that
if a sequence $\hat z^n\in \an_f$ converges in $\ac_f$ then its image
under $\hat\Phi$ converges in $\ac_g$. To this end let us use the dynamical
criterion of convergence, Proposition \ref{dynamical description}.
It is topologically invariant except the condition (i) of big modulus.
However, in the critically non-recurrent case without parabolic points,
 the domains $V_{-N}$ in Proposition \ref{dynamical description}
can be selected as disks $D_{-N}\equiv D(z_{-N}, \eps)$, with a sufficiently
 small $\eps>0$ independent of $N$. 
Indeed, 
by the Ma\~n\'e Theorem  the  pullback of $D_{-N-l}, \; l=0,1,\dots$, of 
 $D_{-N}$ along $\{z_{-N-l}\}$ does not branch, provided $N$ is big enough.
It follows that for sufficiently large $l$, $D_{-N-l}\subset V_{-N-l}$, so that
all properties from Proposition \ref{dynamical description} valid for $V_{-N-l}$
are also valid for $D_{-N}$. 
The domains $D(z_{-N},\eps)$ have inner radius bounded from below, while
the domains $ \Phi U_{-N}$ shrink (By the Shrinking Lemma for $g$).  Hence
$\mod(\Phi D(z_{-N},\eps)\setminus \Phi U_{-N})\to \infty$ as $N\to\infty$,
so that the modulus condition is satisfied. 
end comm}

We may assume that $\Phi$ is quasi-conformal
on the Fatou set
$F(f)$, possibly after applying an equivariant homotopy. 
 Let us give a sketch of this well-known procedure.
Let $\bar a$ be an attracting cycle.
If it is not superattracting, 
 we may choose a fundamental annulus around one of its
points $a$. On this annulus we may homotope $\Phi$, fixing it on the
post-critical points, to some $C^1$ diffeomorphism which conjugates $f$ to
$g$ on the boundary. This homotopy can then be transported by the
action of $f$ and $g$ to the rest of  the attraction basin  $B$ 
of $\bar a$.
 By a Poincar\'e length
argument the tracks of the homotopy have vanishing  Euclidean length near
$\boundary B$, so that it can be extended as the identity to 
$\boundary B$. Finally, Man\`e's Theorem implies that the diameters of the
Fatou components tend to 0, so that the homotopy can be extended 
as the identity to the rest of the sphere.

If $\bar a$ is superattracting,
the B\"ottcher coordinate provides us with
an invariant  circle foliation in a punctured neighborhood of $a$.
Moreover, this foliation is affine (that is, there is a canonical affine
structure on the leaves), as the B\"{o}ttcher coordinate is
unique up to scaling and rotation. 
Select now a fundamental annulus, with the affine circle foliation inside 
and marked post-critical points. There is a homotopy of $\Phi$ in the
fundamental annulus to some diffeomorphism, which respects
this extra structure, and conjugates $f$ and $g$ on the boundary. 
By means of dynamics this homotopy can be spread 
around the whole basin $B$. By the same reason as above it can be
extended to the rest of the sphere as the identity.

In the post-critically finite case the action
of a power of $f$ on the immediate
 basin of $a$ (that is, the component of $D$ containing $a$)
 is conjugate to $z\mapsto
z^d$, and similarly for $g$ (see \cite{lyubich:topological}, Theorem 1.6).
 Then $\Phi$ can be homotope  in the
fundamental annulus to a diffeomorphism which is linear in the
B\"ottcher coordinates. Hence it is conformal on the basin,
and we are in case (2) of the theorem.

\medskip

We next extend $\hat\Phi$ to a conjugacy of the 3-laminations, using
the following elementary fact:

\begin{lemma}{continuous extension}
For any homeomorphism $\phi:\C\to\C$ there is a homeomorphism 
$e(\phi):\Hyp^3\union\C\to\Hyp^3\union\C$ which restricts to $\phi$ on $\C$,
such that the following are satisfied:  
\begin{enumerate}
\item The extension is affinely natural: If $\alpha,\beta$ are (complex) affine 
maps of $\C$ then $e(\alpha)$ and $e(\beta)$ are the unique possible
similarities of $\Hyp^3$, and 
$$e(\alpha\circ\phi\circ\beta) = e(\alpha)\circ e(\phi)\circ e(\beta).$$
\item $e(\phi)$ depends continuously on $\phi$, in the compact-open 
topology on maps of $\C$ and $\Hyp^3$.
\item $e(\phi)^{-1}$ depends continuously on $\phi$ or, 
equivalently, on $\phi^{-1}$. 
\end{enumerate}
\end{lemma}

\begin{pf}
The definition of $e(\phi)$ is the following: 
$$
 e(\phi)(z,t) = \left(\phi(z), \max\limits_{|w| = t} 
 |\phi(z+w)-\phi(z)| \right).
$$
Note in particular that the vertical line over each $z\in\C$ is mapped 
homeomorphically to the vertical line over $\phi(z)$, since the $\max$ 
is monotonic in $t$ as a result of the assumption that $\phi$ is a 
homeomorphism. Hence the map is a homeomorphism. The other properties 
follow easily. Note that part (3) is not completely automatic since 
$e(\phi)^{-1}$ is not in general equal to $e(\phi^{-1})$. 
\end{pf}

As a corollary, we can
extend $\hat\Phi$ leafwise to a map 
$\hat E:\hc_f\to \hc_g$, which is a homeomorphism on every leaf. The extension 
is well-defined because it is affinely natural. 
Note that, on 
the orbifold leaves, we must apply the lemma to the appropriate branched 
cover of the leaf. Since the map back to the orbifold leaf is quotient 
by rotations,  
the affine naturality of the extension implies that the extension is 
well-defined downstairs.

Continuity of $\hat E$ follows from part (2) of lemma \ref{continuous
extension}, applied to a local trivialization, i.e. a product-box (or
orbifold-box) neighborhood in $\hc_f$ and in $\hc_g$. 
Continuity of $\hat E^{-1}$ follows from the same argument, using part (3) of 
lemma \ref{continuous extension}. Thus $\hat E$ is a homeomorphism. 

Again 
the affine naturality of the extension and the fact that $\hat f$ and 
$\hat g$ act  by affine isomorphisms on the leaves imply that $\hat E$ 
conjugates $\hat f$ to $\hat g$. We conclude that it projects to a 
homeomorphism
$$
E:\hc_f/\hat f \to \hc_g/\hat g.
$$

We next show that the $E$ can be deformed to a quasi-isometry:

\begin{lemma}{quasi-isometry} 
There exist $K,\delta>0$ and a map $\hat 
E':\hc_f \to \hc_g$, which agrees with $\hat E$ on $\ac_f$ and is a 
$(K,\delta)$-quasi-isometry on each leaf. 
\end{lemma} 

\begin{pf}
Note that to show a map $h:\Hyp^3\to\Hyp^3$ is a quasi-isometry it 
suffices to show that there exist $\ep_1,\ep_2$ such that for all balls 
$B$ of radius $\ep_1$, $\diam(h(B))\le \ep_2$, and similarly for 
$h^{-1}$. Let us call this property quasi-Lipschitz, so that 
quasi-isometry is equivalent to quasi-Lipschitz in both directions. 

Consider first the case that $f$ (and therefore $g$) has no Fatou domain. 
In this case the convex cores are the entire quotients, and by
Theorem \ref{convex cocompactness1} $\hc_f/\hat f$ and $\hc_g/\hat
g$ are both  
compact. If we fix $\ep_1>0$ then the function $x\mapsto \diam(\hat 
E(B(x,\ep_1)))$ is continuous in $x\in\hc_f$ -- as one can see by 
considering a local trivialization of the lamination. (Here $B(x,\ep_1)$ 
is a leafwise hyperbolic ball of radius $\ep_1$, and
$\diam$ refers to diameter measured inside a leaf.)
By compactness, then, it has a finite upper bound. Since we can do the 
same for $\hat E^{-1}$, we are done in this case. 

In the case where the convex core $\CC_f$ is not all of $\hc_f$, we first 
adjust the map so that it takes a small neighborhood of $\CC_f$ to 
$\CC_g$. 


Let $\CC_f(\ep)$ denote the closed $\ep$-neighborhood of $\CC_f$, by which we 
mean the union of leafwise $\ep$-neighborhoods. Note that 
$\CC_f(\ep)/\hat f$ is still compact. Recall the product structure on
$\hc\setminus \CC_f(\ep)$, discussed in Appendix 2 for the leafwise
case, but extended to the global lamination by virtue of the
discussion in \S \ref{convex hull of lamination} and lemma
\ref{variation of convex hulls} on continuous variation of convex hulls. 
This product structure (in particular projection along the gradient
lines) gives a $C^1$ 
identification between $\boundary \CC_f(\ep)$ and $\fc_f$, and
moreover we obtain a homeomorphism $P_f:\hc_f\union\fc_f \to 
\CC_f(\ep)$ which is the identity on $\CC_f$, and equal to $\Pi_\ep$ on 
$\fc_f$. On each leaf $P_f$ is the map $h^{-1}_{\ep,J}$ discussed in 
the proof of lemma \ref{variation of convex hulls}. Because the
construction is natural, $P_f$ commutes with $\hat f$. 

Letting $P_g$ denote the corresponding construction for $g$, we then have 
(fixing $\ep>0$) a map 
$$
\hat E' = P_g \circ \hat E \circ P_f^{-1}:\CC_f(\ep)\to\CC_g(\ep)
$$
which is a homeomorphism that restricts 
to a $C^1$ diffeomorphism on $\boundary\CC_f(\ep)$, and conjugates $\hat 
f$ to $\hat g$. 
We can extend this to 
a map, also called $\hat E'$, on all of $\hc_f$, using the product 
structure; that is, sending gradient lines to gradient lines at unit speed.

This map is the desired quasi-isometry. On $\CC_f(\ep)$ it is 
quasi-Lipschitz as before, by the same compactness argument on the 
quotient; and similarly for $(\hat E')^{-1}$ on $\CC_g(\ep)$. In the 
exterior, proposition \ref{product structure metric} determines the
metric up to bilipschitz homeomorphism in terms of the metric on the
boundary of $\CC_f(\ep)$ (or $\CC_g(\ep)$). It follows that
it is bilipschitz on the exterior, since
$\hat E'$ is a $C^1$ diffeomorphism on the boundary. (We are also using 
the fact that $\boundary\CC_f(\ep)/\hat f$ is compact to bound the derivatives 
of the map on the boundary).

Since $\hat E'$ is a quasi-isometry it extends continuously to a
quasiconformal homeomorphism on the boundary at infinity, namely
$\ac_f$. It remains to check that the boundary values of 
$\hat E'$ agree with the origional ones of $\hat E$, namely $\hat
\Phi$. In the  
Fatou domain this is automatic from the construction. For any point in
$\jc_f$, we note that it lies in the closure of $\CC_f$.  For any
point $x\in\CC_f$, the maps $\hat E$ and $\hat E'$ differ by an
application of $P_g$, so their leafwise distance is (again by
compactness of the quotient) uniformly bounded. It follows that the
two maps have identical boundary values on $\jc_f$.
\end{pf}

We can now complete the proof of theorem 
\ref{pcf rigidity}. Lemma \ref{quasi-isometry} implies that $\hat\Phi$ 
extends to a quasi-isometry of the 3-laminations -- that is, a map which 
is a quasi-isometry on every leaf, with uniform constants. 
and therefore (lemma \ref{qc extension}) $\hat\Phi$ is in fact a quasiconformal 
map on every leaf, with uniform constant. Since $\hat\Phi$ is just the 
lift of the original conjugacy $\Phi$, we conclude that  $\Phi$ itself 
is quasiconformal. 

This concludes step one of the proof (that topological conjugacy implies
quasi-conformal), which is case (1) of the theorem.
To finish the proof we need to show that a quasi-conformal conjugacy
which is conformal on the Fatou set is M\"obius, except for the Latt\`es
examples. But this is equivalent to the absence of invariant line fields
on the Julia set which follows from
 Proposition \ref{cocompactness and conical pts}
and Corollary \ref{no line fields}.

\section{Further program}
\label{conjectures}

Let us outline some possible directions for further development,
problems and conjectures.

{\it 1. Regular leaf space.} Study the regular leaf space $\RR_f$
in more detail. What is the behaviour of the leaves of $\RR_f$ near
irregular points? In particular, look at the Feigenbaum case.
What happens to $\RR_f$ at a parabolic bifurcation? Other than
rotation domains, are there any leaves which are not dense? (Lemma
\ref{affine leaves are dense} shows that all parabolic leaves are dense.)
Can it happen that a leaf other than a rotation domain does not
intersect the Julia set? 

{\it  2. Type Problem } (see \S 4). Are there hyperbolic leaves in $\RR_f$
  except for Siegel disks and Herman rings? It seems that  the right place to 
  look for hyperbolic leaves are maps with non-locally connected Julia set
 (Cremer points or infinitely renormalizable polynomials of highly 
  unbounded type, see \cite{milnor:localconn}). Prove that all leaves
  of a ``fake  
  Feigenbaum" quadratic (that is, a rational map which is topologically
  equivalent to the Feigenbaum quadratic)  are parabolic.
  Conjecturally there are no fake Feigenbaum maps (a special case of
  the rigidity problem), but this would be the first step of trying to apply
  the laminations to this problem.
  More generally does the topological type of the map determine the
conformal types of the leaves?

{\it 3. Uniqueness problem.}  In general, can one reconstruct $f$ from its
   3-lamination? How does the lamination detect the difference between
polynomial and polynomial-like maps?

 {\it 4. Geometric finiteness.} 
There are many definitions of geometrically finite Kleinian groups,
all equivalent for dimensions 2 and 3 (see Maskit
\cite{maskit:book}, Bowditch \cite{bowditch:geomfinite}).
The definition in terms of finite-sided fundamental domain
(see Ahlfors \cite{ahlfors:geomfinite}) seems to fail altogether in
the lamination context; it is also not equivalent to the others for
hyperbolic manifolds in higher dimensions
\cite{bowditch:geomfinite}.
The definition in terms of conical and parabolic points (Beardon-Maskit
\cite{beardon-maskit}) 
can be translated into the lamination setting.  We expect it to
pick out critically non-recurrent maps with or without
parabolic points. Thurston's definition in terms of finite volume of a
neighborhood of the convex core, or compact thick part of the convex
core (similar also to Marden's definition in \cite{marden})
seems harder to transport to laminations. Is there a good replacement
for the notions of volume and injectivity radius which would make this
translation work?

{\it 5. Deformation theory.}
Describe the space of $\Hyp^3$ laminations, or affine 2-laminations, or
just those arising from rational maps. A fundamental difficulty here is that
there is no common ``universal cover'', as there is for hyperbolic manifolds.

{\it 6. Topology of  $\hc_f/\hat f$.} 
 What is the topological structure of 
 $\hc_f$ and  $\hc_f/\hat f$?  Does  $\hc_f/\hat f$
 always have two ends for quadratic $f$?

Particular cases are the Axiom A polynomials (take $z\mapsto z^2-1$
first) and the Feigenbaum quadratic.  Is there an internal structure
to $\hc_f$ that mirrors the sequence of bifurcations going from
$z\mapsto z^2$ to $f$ (degree 2 case)?

Let us consider the following model. Let $f_c : z\mapsto z^2+c$, $c\in
[ c_0, 0]$, where $c_0$ is the Feigenbaum point, or any point
preceding it.  Let $K_c$ and $J_c$ denote the filled Julia set and
the Julia set for $f_c$.  Consider their lifts $\KK_c$ and $\JJ_c$
to $\RR_f$. 
Consider the set $\MM=\{(c, \hat z): c_0\leq c\leq 0,
\hat z\in \KK_c\}$.

There is a natural projection from $J_c$ onto $J_{c_0}$,
since $J_{c_0}$ is obtained from $J_c$ by some ``pinchings" (compare
Douady \cite{douady:pinching}).  This induces a projection $r_c:
\JJ_c\rightarrow \JJ_{c_0}$. Let us consider the quotient $\MM/\sim$
where the equivalence relation $\sim$ identifies $(c, \hat z)$, $\hat
z\in \JJ_c$ with $(c_0, r_c \hat z)$.  The map $f$ induces a self-map
$\tilde f$ of $\MM/\sim$.

Is $\tilde f: \MM/\sim \to \MM/\sim$ topologically
equivalent to $\hat f_{c_0}: \hc_{f_{c_0}}\union
\ac_{f_{c_0}}\rightarrow \hc_{f_{c_0}}\union \ac_{f_{c_0}}$?

{\it 7. Geometry of $\hc_f/\hat f$. } Give a quasi-isometric model for
 $\hc_f/\hat f$. Does topology of this lamination determine its
 geometry? (It is certainly a quite strong version of the Rigidity
 Problem).

Can one place ``pleated solenoids'' inside $\hc_f/\hat f$, and
use them in analogy with pleated surfaces in hyperbolic
3-manifolds? (In the Feigenbaum case, one can consider the pullback  of the 
little Julia set $J(R^n f)$  to $\ac_f$ (where $R$ denotes the renormalization operator),
 take the boundary of its convex hull in $\hc_f$,
and spread it around by iterates of $\hat f$).  

{\it  8. Spectral Theory.} We define 
  the three dimensional Poincar\'{e} series of $\hat f$ by taking a transversal $K$ of $\hc_f$,
 averaging $\exp (-\rho (\hat f^{-n}x, K))$ along a natural transversal  measure on $K$
 (where $\rho$ stands for the leafwise hyperbolic distance), and summing up over $n$
(see Su \cite{Su} for a  discussion of the transversal  measure).
  Is it true that the corresponding critical exponent coincides with the Hausdorff dimension
 of the conical limit set?  
  A natural further project is to
  develop a spectral theory on the lamination $\hc_f/\hat f$, and
  to study measure and
  dimension of the Julia sets from this point of view
 (compare Sullivan \cite{sullivan:density,sullivan:entropy}, Canary
\cite{canary:laplacian},
Bishop-Jones \cite{bishop-jones:dimension}, Denker-Urbanski \cite{DU}). 
The Ahlfors-type
 argument used in \S 8 of this paper is a first step in this direction.

 {\it 9. Added leaves of $\ac_f$.}
  Can it happen that $\ac_f$
 is not locally compact? This problem requires understanding of the
 added leaves 
 of $\ac_f$.  What one can say about the entire function corresponding to
 the leaf  projection $p: L_{aff}(\bz)\ra L_{aff}(p(\bz))$?    Can it have
 asymptotic values? (In the critically non-recurrent case it is polynomial.)

{\it  10.  Action of rational functions in the Universal space.}
 It would be interesting to have a general idea of this action.
What is the structure of the characteristic attractor $K_f$? 
Is a generic $f: \uu\ra\uu$ injective? More precisely, let us 
consider a functional equation $f\circ\phi=f\circ\psi$ where $\phi,\psi\in \uu$
are meromorphic. Is it true that any solution  of this equation has a form
$\phi=\gamma\circ\psi$ where $\gamma$ is a {\it symmetry} of $f$ (that is,
a M\"obius transformation such that $f\circ\gamma=f$),
or $\phi = \psi \circ \delta$ where $\delta$ is a rotation?
See Fatou \cite{fatou} and Ritt \cite{ritt} for further discussion of
 this problem (we are grateful to A. Eremenko for providing these
 references).

\section{Appendix 1: Circle and polynomial-like maps}
\label{appendix1}

\subsection{Sullivan's laminations for circle maps}

Let $f: S^1\rightarrow S^1$ be a $C^2$ expanding map of the circle
of degree $d>1$.
The expanding property means that there exist constants $C>0$
and $\lambda>1$ such that  $|Df^n(x)|\geq C\lambda^n$, $n=0,1,\ldots$
Sullivan's construction goes as follows (see Sullivan
\cite{sullivan:universalities,sullivan:bounds}, and de Melo-van Strien
\cite{demelo-vanstrien}):  

\medskip\noindent Step (i). Consider the natural extension 
  $\hat f: \NN_f\rightarrow \NN_f$. Topologically $\NN_f$ is the
  standard solenoid over the circle. Dynamically $\hat f$ is a 
  hyperbolic (in the sense of Anosov and Smale) map with 
 one-dimensional unstable leaves.

\medskip\noindent Step (ii). Supply the leaves with the affine structure
  by means of the explicit formula  (\ref{chart})
  (existence  of the limit follows from the standard distortion estimates
  for hyperbolic maps). The map $\hat f$ preserves this structure.

\medskip\noindent Step (iii). Attach hyperbolic planes to the leaves
  and extend $\hat f$ to the corresponding hyperbolic 2-lamination
  $\HH^2_f$ acting isometrically on the leaves.

\medskip\noindent Step (iv). Take the quotient $\HH^2/\hat f$. This is
 Sullivan's Riemann surface lamination associated to $f$. Topologically
  it is a solenoidal fibration over the circle.

\medskip
The main difference between this construction and the one outlined in the
Introduction is related to the critical points on the Julia set.
These tend to distort the affine structures and complicate the transversal
behavior of the leaves. Also, as we have seen,  even in the Axiom A case
the topological structure  of the 3-lamination is not at all obvious.

Sullivan  constructed  2-laminations to build up the deformation space
of expanding circle maps. We try to study rigidity phenomenon by means of
 3-laminations. This is a usual philosophical difference between dimensions
two and three.

\subsection{ Polynomial-like maps: globalization of the leaves }
\label{polynomial-like}
Polynomial-like maps are not globally defined, and certainly cannot be 
in general extended to the whole sphere. However, such a globalization
can be carried out on the natural extension level.  Lemma \ref{direct limit}
shows that it leads to the same object, provided the map was a priori
globally defined.

Let $U$ and $V$ be two open  sets of $\C$ such that $\cl U\subset V$, 
and 
$f:  U\rightarrow V$ be
an analytic branched covering. Keep in mind Douady-Hubbard
polynomial-like \cite{douady-hubbard:polylike}
maps, generalized polynomial-like maps \cite{lyubich:measure}, 
or a rational function $R$
restricted on the sphere minus an invariant neighborhood of attracting
cycles.

For such a map
we can consider the space ${\cal N}_f$ of backward orbits, and lift $f^{-1}$
to this space as the
map which forgets the first coordinate:
${\hat g}\equiv {\hat f}^{-1}: {\cal N}_f\rightarrow {\cal N}_f$. 
This map is injective but not surjective:
its image consists of the orbits which start with a $z_0\in U$.

To make it invertible, let us consider the {\it inductive (direct) limit}
of $${\cal N}\underset{\hat g}\rightarrow {\cal N}\underset{\hat g}
 \rightarrow{\cal N}\underset{\hat g}\rightarrow\ldots,$$
which is defined in the following way.
Take infinitely
many copies ${\cal N}^{m}$ of the same space ${\cal N}$. Let us embed
${\cal N}^{m}$ into ${\cal N}^{m+1}$ by means of the map
 $$i_m\equiv {\hat g}: {\cal N}^{m}=
{\cal N}\rightarrow {\cal N}={\cal N}^{m+1}.$$ In other words, we identify
a point $\hat z\in {\cal N}^m$ with the point $i_m \hat z\in {\cal N}^{m+1}$.
Thus we obtain an increasing sequence of the spaces
\begin{equation}\label{inclusions}
{\cal N}^0\hookrightarrow {\cal N}^1\hookrightarrow{\cal N}^2
\hookrightarrow\ldots
\end{equation}
Let ${\cal D}\equiv {\cal D}_f=\cup {\cal N}^m$.  To define a topology 
on ${\cal D}$, let us call a set $W\subset {\cal D}$ open if 
$W=\cup W_i$ where $W_i$ is an open set in ${\cal N}^i$.  

The map ${\hat g}: \NN^k\rightarrow \NN^k$ respects the embeddings 
$i_m: {\cal N}^m\hookrightarrow {\cal N}^{m+1}$,
and hence induces the self-map of ${\cal D}$, which we will denote by the same
letter. Moreover, ${\hat g}$ homeomorphically maps 
 ${\cal N}^m$ onto $i_{m-1}{\cal N}^{m-1}$,
$m>0$, so that it is invertible on ${\cal D}$.
We will keep the notation ${\hat f}$ for $\hat g^{-1}$.

\begin{lemma}{direct limit} 
Assume that a branched covering $f:  U\rightarrow V$
is the restriction of a rational endomorphism $R: \bar{\C}\rightarrow \bar{\C}$
such that $\C\ssm V$ is contained in the basin of attraction of a finite 
attracting set $A$.
Then $\hat f: {\cal D}_f\rightarrow {\cal D}_f$ is naturally conjugate
to $\hat R: {\cal N}_R\ssm\hat A\rightarrow {\cal N}_R\ssm\hat A$.
\end{lemma}

\begin{pf} 
Let us consider the following commutative diagram:
\begin{equation*}
\begin{array}{cccccc}
\NN^0  & \underset{i_0}{\hookrightarrow} & 
\NN^1  & \underset{i_1}{\hookrightarrow} & 
\NN^2  & \underset{i_2}{\hookrightarrow} \cdots \\
\Big\downarrow\vcenter{\rlap{id}} & &
\Big\downarrow\vcenter{\rlap{$\hat R$}} & &
\Big\downarrow\vcenter{\rlap{$\hat R^2$}} & \\
\NN & \underset{i}{\hookrightarrow} & 
\hat R\NN & \underset{i}{\hookrightarrow} & 
\hat R^2 \NN & \underset{i}{\hookrightarrow} \cdots \\
\end{array}
\end{equation*}
where $\NN\equiv \NN_f$, 
the upper line is the sequence (\ref{inclusions}) for $\hat f$, while
the lower one is the sequence of natural inclusions. It induces a
homeomorphism between $\DD_f$ and $\cup \hat R^n \NN=\NN_R\ssm \hat A$,
which is the desired conjugacy.   
\end{pf}

\section{Appendix 2: Background material}
\label{appendix}
\label{appendix2}

\subsection{Dynamics}
\label{dynamics for geometers}
We assume the following background in holomorphic dynamics:

\begin{itemize}
\item Classification of periodic 
points as {\it attracting, repelling, parabolic, Siegel} and {\it Cremer},
and the local dynamics near these points;

\item
Notions of the Julia set $J(f)$ and the  Fatou set $F(f)$;

\item Classification  of components of the Fatou
set as {\it attracting basins, parabolic basins, Siegel disks and
 Herman rings}; Siegel disks and Herman rings will be also called the
{\it rotation sets}.

\item The notion of an {\it Axiom A or hyperbolic}
rational function. There are two equivalent definitions of this property:

\begin{itemize}
\item
All critical points are in basins of attracting cycles;

\item
The map is uniformly expanding on the Julia
set, that is, there exist constants $A>0$ and $\lambda>1$ such that
for any $z\in J(f)$,
 $$\|Df^n (z)\|\geq A\lambda^n,\; n=0,1,2\ldots,$$
where $\|\cdot\|$ denotes the spherical metric.
\end{itemize}
\end{itemize}

All this material can be found in any book or survey in holomorphic
dynamics -- e.g. \cite{carleson-gamelin,lyubich:topological,milnor:intro}.

As usual, $\omega(z)\equiv \omega_f(z)$
 denotes the $\omega$-limit set of a point z.
A point $z$ is called {\it recurrent} if $z\in \omega(z)$.
Given a set $Z$, let 
$$\operatorname{orb}(Z)=\bigcup_{z\in Z}\operatorname{orb} z, \qquad\omega(Z)=\bigcup_{z\in Z}\omega(z).$$
Let $C$ denote the set of citical points of $f$, and $C_r$ the set 
of recurrent critical points.

The critical values of $f^n$ are the points of $f^k C,\; 1\leq k \leq n$.
So if a simply connected neighborhood $U$ does not meet
$\operatorname{orb} C$ then 
all inverse branches of $f^{-n} $ are well defined univalent functions in $U$.

The {\it non-linearity}, or {\it distortion} of a conformal map
$\psi: U\hookrightarrow \C$ is defined as
$$\Dis({\psi})=\sup_{z,\zeta\in U}
    \log \left|{\psi' (z)\over \psi'(\zeta)}\right|. $$

\proclaim Koebe Distortion Theorem. Let $\psi: B(a, r)\hookrightarrow \C$
be a conformal map, $k<1$. Then the distortion of $\psi$ in $B(a, kr)$
is bounded by a constant $C(k)$ independent of $\psi$. Moreover
$C(k)=O(k)$ as $k\to 0$.

Let $U\subset\bar\C$ be any domain. Let us select a base point $z\in U$,
and count its $n$-fold preimages: $z^n_i$. Let $U^{-n}_i$ denote a
component of  $f^{-n} U$ containing $z^n_i$. 
This specifies a ``multi-valued branch" $f^{-n}_i$ of the inverse map.
(The reader can think of these branches as functions living
on appropriate Riemann surfaces, or as equivalence relations, or just
as a convenient way of describing the situation).
Singular points for an inverse branch are critical values
for the direct map. There is a natural way of composing and restricting
the inverse branches (with an appropriate adjustment of the base points,
which may change only the way of counting). 

The following lemma  is a variation of a well-known fact
(compare \cite{lyubich:topological}, Proposition 1.10).
As it plays a crucial role for this paper, we will include the
proof.

\proclaim Shrinking Lemma. Let $f$ be a rational map of degree $d>1$.
Let $U\subset \C$ be a domain
which is not contained in any rotation set of $f$, and let $k$ be a
natural number. 
Let us consider a family $\{f^{-n}_i\}$
of all inverse branches  in $U$ with at most $k$ singular points
(counting with multiplicities).
Then  for any domain $W$ compactly contained in $U$,
$\diam(f^{-n}_i|W)\to 0$ as $n\to\infty$ independently of $i$ 
(where $\diam$ denotes spherical diameter).

\begin{pf} 
We first consider the case that $U$, and every
pullback $U^{-n}_i$, are disks. Let $z\in U$ be a point outside any
rotation domain of $f$.

Let 
$\Phi_{n,i}: D\to U^{-n}_i$ be a Riemann mapping taking $0$ to a preimage
of $z$, where $D$ is the unit disk. Then $\pi_{n,i} = f^n\circ
\Phi_{n,i}$ is a proper branched covering from $D$ to $U$, with at most
$k$ critical points counted with multiplicity.  (One can think of the
disk $D$ here as the 
Riemann surface over $U$ for the corresponding branch of the inverse
function.)  

Let $\alpha_1,\ldots,\alpha_k$ be a periodic cycle of $f$ of length at
least 3, not meeting some neighborhood of $z$. Then no preimage of
this neighborhood meets the cycle either. 
By normality of the family $\{\pi_{n,i}\}$, 
there must be some disk $D'$ compactly contained in $D$ such that
$\Phi_{n,i}(D')$ omits $\{\alpha_j\}$ for all $n,i$, and such that
$\pi_{n,i}(D')\ni z$. 
Thus $\{\Phi_{n,i}\}$ is  a normal family on $D'$. 

Because of the bound $k$ on the number of critical points of $\pi_{n,i}$,
there is some $\delta$ such that the disk $B=B(z,\delta)$ is contained
in $\pi_{n,i}(D')$ for all $n,i$ (one can show this for example by noting that 
$\pi_{n,i}^{-1}(U\setminus B(z,\delta))$ contains an annulus whose
modulus is bounded below
depending only on $k$ and $\delta$, and goes to $\infty$ as $\delta\to 0$).
We now claim that 
the diameters $\diam(B^{-n}_i)$ go to 0 uniformly. 

If not, we can extract a convergent subsequence
$\Phi_{n_k,i_k}|_{D'}$, and conclude 
that for the limit 
point $z_\infty = \lim \Phi_{n_k,i_k}(0)$  there is a neighborhood
$B_\infty$ whose images under arbitrarily high iterates are in $U$.
This implies in particular that $B_\infty$ (and therefore $B$) is
disjoint from the Julia set (as any neighborhood intersecting the Julia
set covers it under some iterate of $f$). By a smaller choice of
$\delta$ we may assume it is compactly contained in the Fatou set. 
Thus, either forward iterates of
$B_\infty$ under $f$
limit to an attracting/parabolic periodic cycle, or $B_\infty$ is contained
in a rotation domain. The former is impossible since
$f^{n_k}(B_\infty)$ limits onto all of $B$. 
The latter is ruled out by the choice of $z$.

It now follows that $\diam(W^{-n}_i) \to 0$ for any $W$ compactly
contained in $U$, since $\Phi'_{n,i}$ must converge to 0 uniformly on
compact sets.

To treat the general case, take a finite covering of $W$ by disks $D$
compactly contained in $U$, none of which are contained in a rotation domain.
We must consider the possibility that some of the pullbacks $D^{-n}_i$ are
not disks. 
For any $\ep>0$ there exists $N = N(D,\ep)>0$ such that, if $D^{-n}_i$
is a disk 
and $n\ge N$, then $\diam D^{-n}_i \le \ep$. For if not, we could find
a subfamily of pullbacks, all disks, whose diameters fail to shrink to
0. The previous argument applies, so this is impossible. 

Thus, let $\ep$ be less than half the distance between any two
critical values of $f$. Then the preimage of any disk of diameter less
than $\ep$ is a disjoint union of disks. It follows that, if some
$D^{-n}_i$ is not a disk then some image $D^{-m}_i$ of it, with $0\le
m \le N$, is also not a disk. That is, the transition from disk to
non-disk occurs in the first $N$ levels. Thus, if we remove from
consideration the finite number of non-disks $D^{-n}_j$ with $n\le N$,
and all their preimages, we are left with a family in which all
preimages are disks. For this subfamily, we have uniform shrinking by
the previous arguments.

For each of the finitely many non-disks $D^{-n}_j$ ($n\le N$), we can
now repeat the argument, covering $W^{-n}_j$ with disks not contained
in rotation domains, and so on. However now the bound on the number of
singular points is $k-2$, since in the transition from disk to
non-disk at least two singular points must be used. We can therefore
obtain a uniform rate of shrinking for this family, by induction on
$k$. This concludes the proof.
\end{pf}

A key result on  critically non-recurrent rational
maps is the following theorem of Ma\~n\'e \cite{mane:fatou}
closely related to the Shrinking Lemma.

\proclaim Ma\~n\'e's Theorem { \cite{mane:fatou}}.
Let $f:\bar \C\to\bar\C$ be a rational map. If a point $x\in J(f)$ is
neither a parabolic periodic point, nor belongs to the $\omega$-limit
set of a recurrent critical point then, for all $\ep>0$, there exists
a neighborhood $U$ of $x$ such that for all $n\ge 0$ every connected
component of $f^{-n}(U)$ has diameter $\le \ep$.

\noindent{\bf Chebyshev and Latt\`es examples.} Let us finally dwell
on the remarkable examples of rational functions whose dynamics
often present some special features.

The {\it Chebyshev polynomial} $p_d$ of degree $d$
 can be defined by means of the functional
equation $p_d (cos z)= cos ( dz)$. In other words, 
consider the dilation map $T_d: z\mapsto dz$ on the
 cylinder $C=\C/2\pi{\Z}$. Then 
$p_d$ is the quotient of this map via the involution $z\mapsto -z$. 

The Julia set of $p_d$  coincides with
the interval $[-1,1]$. The endpoint 1 is always fixed, while
-1 is either fixed (for odd degrees) or pre-fixed (for even degrees).
Any critical point is mapped by $p_d$ to one of the endpoints.

Similarly, the Latt\`es examples come from the functional equations
$f_d(P(z))=P(dz)$, where $P: \C\ra \bar\C$ is a Weirestrass $P$-functin,
$\deg f_d=|d|^2$ where $d$ is not necessarily integer.
They can be viewed as quotients of torus
endomorphisms. That is, let $\T= \C/\Lambda$ be a torus, where $\Lambda$
is a lattice. Then identifying $z$ with $-z$ sends $\T$ to $\bar\C$ via a
two-fold branched cover.  If $T_d(\Lambda)\subset \Lambda$ then the dialtion
 $T_d$ induces a torus endomorphism, which further
projects to a rational map  of $\bar\C$ of degree $|d|^2$. 
(This occurs for all integer $d$'s on any torus, 
but also for some special tori and special  non-real values
of  $d$: take, e.g.,  the standard lattice $\Lambda=\Z^2$ and $d=1+i$).

The Julia set of the Latt\`es examples is the whole sphere.
Like in the Chebyshev case,  every
critical point of a Latt\`es map  is pre-fixed. 

The following dynamical characterization of these examles is well-known:

\begin{proposition}{C and L} Assume that a rational map $f$ has a periodic
point $a\in J(f)$ such that every backward trajectory $a=a_0, a_{-1},\dots$
which passes through $a$ only finitely many times hits a critical point.
Then $f$ is either Chebyshev or Latt\`es.
\end{proposition}

We will see in this paper how this property manifests itself in the lamination structure.  

For integer values of $d$ the Latt\`es
 maps are  quasi-conformally deformable,
since $\Lambda$ may be varied (or, since the constant line field
on the torus is dilation invariant).
Conjectually they are the only examples which admit quasi-conformal
deformations on the Julia set. We will see a lamination reasoning
behind  this conjecture.

\subsection{Geometry}
\label{geometry for dynamicists}

\subsection*{Hyperbolic geometry and convex hulls}
We assume familiarity with the hyperbolic space $\Hyp^3$ and its
boundary at infinity the Riemann sphere. (See e.g. Beardon
\cite{beardon}, Thurston \cite{wpt:textbook}).
Most natural for us will be 
the upper half space model $\C \times \R_+$.

We recall some fundamental facts about hyperbolic convex hulls. 
Most of these facts appear in Epstein-Marden \cite{epstein-marden}, or
can be obtained from that paper with a small amount of effort. 

The convex hull $C=C(E)\subset
\Hyp^3$ of a closed set $E$ on the Riemann sphere 
$\bar\C$ is defined as the smallest convex set in $\Hyp^3$ whose
closure in $\Hyp^3\union\bar\C$ contains $E$. Equivalently, $C$ is the
intersection of all closed half-spaces in $\Hyp^3$ containing $E$ at
infinity.  Provided $E$ is not contained in a round circle, $C\union
E$ is homeomorphic to a closed 3-ball, and $\boundary C$ is a
subsurface of $\Hyp^3$, which is isometric to a complete hyperbolic
surface, using the metric of shortest paths in $\boundary C$.

The geometry of the complement $\Hyp^3-C$ is well-understood
We begin with the projection
$\Pi:\Hyp^3\to C$ assigning to $x\in \Hyp^3$ the point in $C$ nearest to
$x$, which is unique by the convexity of $C$. This projection also
extends continuously to $\bar\C - E$.

Let $d:\Hyp^3\to [0,\infty)$ be the distance
function $d(x) = d_{\Hyp^3}(x,C)$. This is a $C^1$ function in $\Hyp^3-C$,
and its gradient is the unit vector tangent to the geodesic through
$x$ and $\Pi(x)$, and pointing away from $\Pi(x)$ (lemma 1.3.6 in
\cite{epstein-marden}). In fact these geodesics are the integral lines
of this gradient field, and they foliate $\Hyp^3\ssm C$. The gradient vector
field itself is Lipschitz, with a uniform constant outside a
neighborhood $C_\ep = d^{-1}([0,\ep])$, for any fixed $\ep>0$ (see
\S 2.11 in \cite{epstein-marden}).

The level surfaces $S_\ep = d^{-1}(\ep)$ are, therefore, $C^1$
submanifolds for $\ep>0$, and are all homeomorphic via the
gradient flow. 
Since each gradient line terminates at infinity, the
level surfaces can be identified with $\bar\C\ssm E$, which we may label
$S_\infty$. Thus we have a natural product structure identifying
$\Hyp^3\union\bar\C\ssm (C\union E)$ with $(0,\infty]\times S_\ep$ for
$\ep\in(0,\infty]$. 

The identification between $S_\ep$ and $S_\infty$
is  a quasiconformal map, and in fact the following is a consequence
of Theorem 2.3.1 in \cite{epstein-marden}: 

\begin{proposition}{product structure metric}
Let $\sigma$ denote the Poincar\'e metric on $S_\infty = \bar\C-E$. 
Let $\rho$ denote the metric on $(0,\infty)\times S_\infty$ given
infinitesimally as
$$
d\rho^2 = dr^2 + (\cosh^2 r) d\sigma^2
$$ 
where $r\in (0,\infty)$ is the first coordinate.
The identification of $(\ep,\infty)\times S_\infty$  with $\Hyp^3-C_\ep(E)$ is
bilipschitz with constant $L$ depending only on $\ep>0$. 
\end{proposition}

The dependence of $C(E)$ (or $C_\delta(E)$) on $E$ is {\em continuous},
with respect to the Hausdorff topology on closed subsets of the ball
$\Hyp^3 \union \bar\C$. This is easy in our setting;
a proof for 
a more general context appears in Bowditch \cite{bowditch:hulls}.
In fact more is true: on compact sets in $\Hyp^3$, a small variation of
$E$ produces a locally homeomorphic deformation of $C_\delta$:

\begin{lemma}{variation of convex hulls}
Let there be given a closed $E_0\subset \bar\C$, a hyperbolic $R$-ball
$B(x,R)$ around 
a point $x\in\Hyp^3$, and  $\delta>0$. For each $\ep>0$ 
there is a neighborhood $U$ of $E_0$ in the
Hausdorff topology on closed subsets of $\bar\C$ such that, for any
$E\in U$, there is a 
$(1+\ep)$-bilipschitz map
$\Psi_E:B(x,R)\to \Hyp^3$ fixing $x$, 
such that $\Psi_E^{-1}(C_\delta(E))  = C_\delta(E_0)\intersect B(x,R)$.
\end{lemma}

{\bf Remarks.} (1) In particular, note that ($\delta$-neighborhoods
of) convex hulls of sufficiently nearby sets are, locally,
homeomorphic, {\em even if the sets themselves are not homeomorphic}.
(2) We take $C_\delta$ rather than $C$ itself here in order to avoid the
exceptional case when $E_0$ lies on a round circle. Then the convex
hull fails to have interior, and is not homeomorphic to convex hulls
of nearby sets. In all other situations the lemma holds for $C_0 = C$. 

\begin{pf}
We give only a sketch, and refer the reader to \cite{epstein-marden}
for a thorough treatment of the techniques.

Using the product structure on $\Hyp^3-C_\delta(E)$ discussed above,
there is a homeomorphism
$h_{\delta,E}:C_\delta(E) \to \Hyp^3 \union S_\infty(E)$, which
expands segments to gradient lines, and is the identity on $C(E)$. 
Now note that, 
for a fixed ball $B(x,R)$ and $E$ sufficiently close to
$E_0$, the image $h_{\delta,E_0}(B(x,R)\intersect C_\delta(E_0))$
misses $E$. Therefore the map $h_{\delta,E}^{-1}\circ h_{\delta,E_0}$ is
defined on $B(x,R)\intersect C_\delta(E_0)$. Extend to the rest of
$B(x,R)$, again using the product structure. 
\marginpar{is bilipschitz clear?}
\end{pf}

\subsection*{Quasi-isometries and QC maps}
We call a map $h:\Hyp^3\to\Hyp^3$ a $(K,\delta)$ quasi-isometry if the 
following holds for all $p,q\in\Hyp^3$: 
$$
{1\over K} d(p,q) -\delta \le
d(h(p),h(q)) \le Kd(p,q) + \delta.
$$

The connection (in one direction) of quasi-isometries to
quasi-conformal maps is given by the following lemma. For a proof, see
Thurston \cite{wpt:notes} or (in the more general context of
hyperbolic spaces in the sense of Gromov) \cite{c-d-p,ghys-harpe}.

\begin{lemma}{qc extension}
Given $(K,\delta)$ there exists $L$ so that any
$(K,\delta)$-quasi-isometry 
$h:\Hyp^3\to\Hyp^3$ extends continuously to an $L$-quasiconformal
homeomorphism $\til h:\bar\C\to\bar\C$.
\end{lemma}

\ifx\undefined\bysame
\newcommand{\bysame}{\leavevmode\hbox to3em{\hrulefill}\,}
\fi


\end{document}